\documentclass[a4paper,11pt,reqno]{amsart}

\usepackage[margin=1in,includehead,includefoot]{geometry}
\usepackage[utf8]{inputenc}
\usepackage[T1]{fontenc}
\usepackage{lmodern,microtype}
\usepackage{upref}
\usepackage{mathtools,amssymb,amsthm}
\usepackage[foot]{amsaddr}
\usepackage{enumitem}
\usepackage{ifthen}
\usepackage{graphicx}
\usepackage{hyperref}
\usepackage{marginnote}
\usepackage[margin=5mm]{caption}
\usepackage{todonotes}
\usepackage{tikz}
\usepackage{bookmark}
\usepackage{siunitx}
\usepackage[b]{esvect} 
\usepackage{xcolor}
\usepackage{mycommands}

\graphicspath{{./graphics/}}

\definecolor{darkblue}{rgb}{0,.25,.7}
\newcommand{\problembox}[1]{\fcolorbox{darkblue}{white}{\color{darkblue}P#1}}
\newcommand{\openproblem}[2]{\marginnote{\problembox{#1}}{\color{darkblue}#2}}
\newcommand{\solvedproblem}[2]{\marginnote{\problembox{#1} \\ \color{darkblue}{solved}}{\color{darkblue}#2}}

\theoremstyle{remark}


\makeatletter
\def\l@subsection{\@tocline{2}{0pt}{2pc}{6pc}{}}
\makeatother

\hypersetup{
  pdftitle={Combinatorial Gray codes---an updated survey}
  pdfauthor={Torsten M\"utze}
}

\title{Combinatorial Gray codes---an updated survey}

\author{Torsten M\"utze}
\address[Torsten M\"utze]{Institut f\"ur Mathematik, Universit\"at Kassel, Germany}
\email{tmuetze@mathematik.uni-kassel.de}

\thanks{Torsten M\"utze is also affiliated with the Faculty of Mathematics and Physics, Charles University Prague, Czech Republic.
He was supported by Czech Science Foundation grant GA~22-15272S, and he participated in the workshop `Combinatorics, Algorithms and Geometry' in March 2024, which was funded by German Science Foundation grant~522790373.}

\begin{document}

\begin{abstract}
A combinatorial Gray code for a class of objects is a listing that contains each object from the class exactly once such that any two consecutive objects in the list differ only by a `small change'.
Such listings are known for many different combinatorial objects, including bitstrings, combinations, permutations, partitions, triangulations, but also for objects defined with respect to a fixed graph, such as spanning trees, perfect matchings or vertex colorings.
This survey provides a comprehensive picture of the state-of-the-art of the research on combinatorial Gray codes.
In particular, it gives an update on Savage's influential survey [C.~D.~Savage. A survey of combinatorial Gray codes. SIAM Rev., 39(4):605--629, 1997.], incorporating many more recent developments.
We also emphasize the connections to closely related problems in graph theory, algebra, order theory, geometry and algorithms, which embeds this research area into a broader context.
Lastly, we collect and propose a number of challenging research problems, thus stimulating new research endeavors.
\end{abstract}

\maketitle

\tableofcontents

\section*{Version history}

This is version 2 of the manuscript.
Each item summarizes the changes compared to the previous version.

\begin{enumerate}[label=\arabic*.,leftmargin=4mm, noitemsep, topsep=3pt plus 3pt]
\item (2022-02-02) Initial version.
\item (2024-07-31) New Section~\ref{sec:01} on 0/1-polytopes; new Section~\ref{sec:otrees} on ordered trees; Figure~\ref{fig:comb3} was split into two (now Figure~\ref{fig:comb4}); new Figure~\ref{fig:counter}; problems P27, P28 and P30 were solved; added progress on problems P10, P29, P42 and P56; added new problems P61--P67; included new results and added or updated several references; added version history (added version history (added version history (added version history (... [stack overflow] ...)))).
\end{enumerate}

\section{Introduction}

In computer science and discrete mathematics, we frequently encounter different classes of combinatorial objects.
For example, we may be interested in all subsets or all permutations of some ground set, all strings over some alphabet, all binary trees on a certain number of nodes, or all spanning trees of a graph, to mention just a few examples.
A \emph{combinatorial Gray code} is a listing of the combinatorial objects in the class of interest that contains each object exactly once such that any two consecutive objects in the list differ only by a `small change'.
This concept is named after Frank Gray, a physicist and researcher at Bell Labs, whose patent~\cite{gray_1953} from 1953 describes a method for listing all bitstrings of length~$n$ such that any two consecutive strings differ in exactly one bit.
Combinatorial Gray codes have a large number of real-world applications, and we shall see that this research area has ramifications into graph theory, algebra, order theory, geometry and algorithms.

The starting point of this work is Savage's influential survey~\cite{MR1491049} from~1997.
Since its publication, the area of combinatorial Gray codes has flourished tremendously, and many new and exciting results have appeared.
In particular, several of the hard open problems mentioned in her survey are now solved, and various general techniques for constructing and arguing about Gray codes have emerged, drawing on the aforementioned connections to related mathematical disciplines, which embeds this research area into a much broader context.
This survey aims to provide a comprehensive overview over the state-of-the-art, and to provide pointers to closely related research strands.
Another goal is to collect and propose a number of challenging research problems, thus stimulating new research endeavors.
These open problems of varying difficulty are marked with~\framebox{P} on the page margins.
Some of these problems are also mentioned in the problem collection by Jackson, Stevens and Hurlbert~\cite{MR2548549}, and some in Gould's survey~\cite[Sec.~6]{MR1106528} on Hamilton cycles in graphs.
Knuth's book~\cite{MR3444818} also contains several challenging research problems on Gray codes.

This survey focuses mainly on basic combinatorial objects of the types mentioned before, i.e., bitstrings, permutations, words, trees etc.
We also discuss some objects defined with respect to a fixed graph, such as spanning trees, perfect matchings, vertex colorings etc., and objects defined with respect to a fixed poset, such as linear extensions and ideals.
These and other discrete structures and small change operations between them are subject of the broader research area of \emph{combinatorial reconfiguration}; see the surveys~\cite{DBLP:books/cu/p/Heuvel13, MR3798903, mynhardt_nasserasr_2019}.
For important classes of combinatorial objects, we often refer to the corresponding counting sequence in the OEIS~\cite{oeis} by its ID number in square brackets.
As this survey is dynamic, it allows incorporating new results and developments as they arise.
The reader is encouraged to contact the author with any such updates, corrections or pointers to missing references, and also with suggestions for additional open problems.

Undoubtedly, one of the main applications of Gray codes is \emph{combinatorial generation}, i.e., the algorithmic task of efficiently producing a listing of the combinatorial objects of interest, each object exactly once.
Gray codes lend themselves to fast generation, because of the `small change' property between consecutive objects, provided that those changes can be computed efficiently.
However, there are orderings of objects that are not Gray codes, and that can still be computed efficiently in an amortized sense.
For example, in the counting order of all bitstrings of length~$n$, which is the same as lexicographic order, the average number of flipped bits is $2-o(1)$, but in some steps, many bits (up to $n$) are flipped.
Having said that, in this survey we focus entirely on Gray code listings for combinatorial objects, and not on the broader task of efficient combinatorial generation.

There are many excellent books on combinatorial generation, which include a lot of material on Gray codes.
Specifically, the author drew a great amount of inspiration from Ruskey's book draft~\cite{ruskey_2003} and from Knuth's book~\cite{MR3444818}.
Moreover, the book by Nijenhuis and Wilf~\cite{MR0396274} is a classic, as well as Wilf's update~\cite{MR993775}.
Other books that contain material on the subject are~\cite{MR0335266,MR0471431,MR561109,MR799324,MR1116324,kreher_stinson_1999}.
Another useful survey on combinatorial generation is~\cite{wasa_2016}, and many of the results listed there are based on the well-known reverse search method of Avis and Fukuda~\cite{MR1380066}.

\emph{De Bruijn sequences} and \emph{universal cycles} are closely related to Gray codes, but beyond the scope of this survey.
We refer to the papers~\cite{MR652466,MR1197444} and to the online resource~\cite{bruijn} for readers interested in this subject.

\subsection{Outline of this paper}

In Section~\ref{sec:prelim} we introduce some terminology that is used throughout this paper, and we collect a few general observations that are valid for any Gray code problem.
We also show how Gray codes relate to several other research areas in mathematics and computer science.
Sections~\ref{sec:bits}--\ref{sec:greedy} describe the literature on combinatorial Gray codes, sorted by objects and sometimes by methods.
Some objects generalize multiple other objects, and assigning them to one particular section is somewhat arbitrary.
For example, multiset permutations generalize combinations and permutations, and they are treated in Section~\ref{sec:perm} on permutations.
This survey contains numerous figures of different Gray codes, which facilitate comparison between the different constructions, thus also stimulating new research questions.
Moreover, the figures convey the wealth and aesthetic appeal of the subject.

\section{Preliminaries}
\label{sec:prelim}

\subsection{Flip graphs}

\begin{figure}
\centerline{
\begin{tabular}{c}
\includegraphics{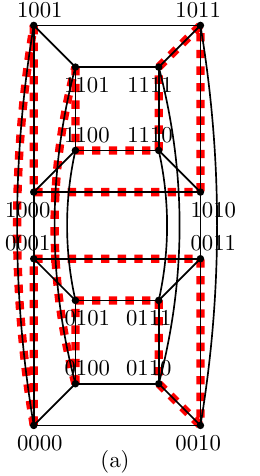}
\includegraphics{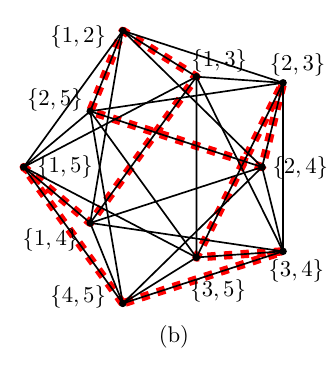}
\includegraphics{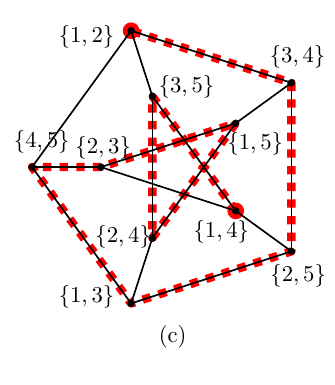} \\
\includegraphics{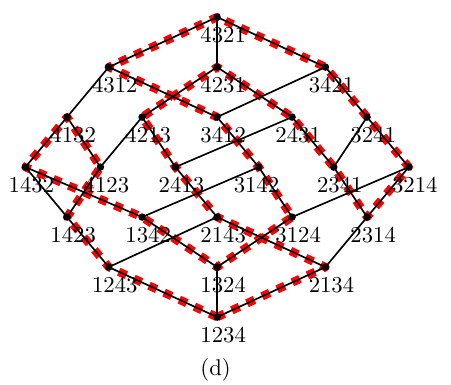}
\includegraphics{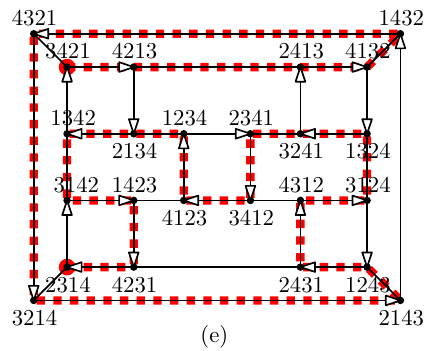} \\
\includegraphics{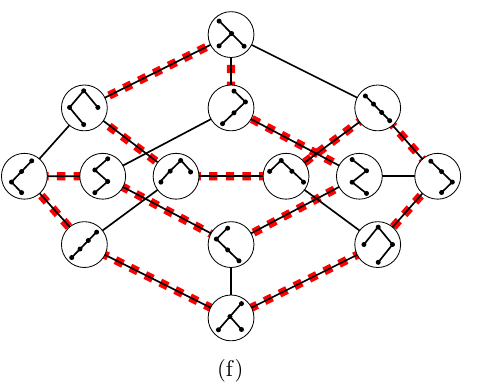}
\includegraphics{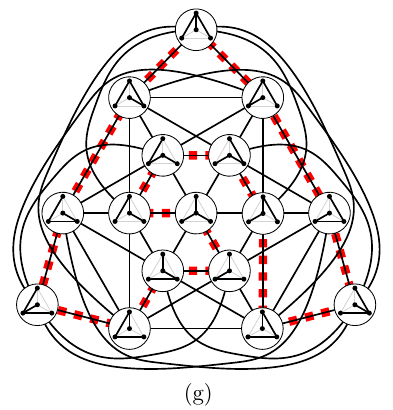} \\
\end{tabular}
}
\caption{Different flip graphs with a Hamilton path or cycle:
(a)~4-cube~$Q_4$; (b)~Johnson graph~$J(5,2)$; (c)~Kneser graph~$K(5,2)$ (=Petersen graph); (d)~permutahedron; (e)~sigma-tau graph (permutations under left shifts or swaps of the first two entries); (f)~associahedron; (g)~spanning trees of~$K_4$ under edge exchanges.
}
\label{fig:flip}
\end{figure}

Nearly every Gray code problem can be described by considering the corresponding \emph{flip graph}.
The vertices of this graph are the combinatorial objects to be listed in the Gray code, and the graph has an edge connecting any two objects that differ by the desired `small change' operation, sometimes generically referred to as a \emph{flip}.
If the flip operation is not an involution, then the edges of the flip graph are directed.
In this way, we obtain a huge family of interesting graphs; see Figure~\ref{fig:flip}.
For example, the flip graph of bitstrings of length~$n$ under flipping a single bit is the \emph{$n$-dimensional hypercube~$Q_n$}, or \emph{$n$-cube} for short (Figure~\ref{fig:flip}~(a)).
The flip graph of subsets of~$[n]:=\{1,2,\ldots,n\}$ of fixed size~$k$ under exchanging a single element is the \emph{Johnson graph~$J(n,k)$} (Figure~\ref{fig:flip}~(b)).
The flip graph of subsets of~$[n]$ of fixed size~$k$ under disjointness is the \emph{Kneser graph~$K(n,k)$} (Figure~\ref{fig:flip}~(c)).
The flip graph of permutations of~$[n]$ under adjacent transpositions is the \emph{permutahedron} (Figure~\ref{fig:flip}~(d)).
The flip graph of binary trees with $n$ nodes under tree rotations is the \emph{associahedron} (Figure~\ref{fig:flip}~(f)).
Clearly, finding a Gray code is equivalent to finding a Hamilton path or cycle in the corresponding flip graph.
The term `Gray code' can refer to both a Hamilton path or cycle, and we call a Gray code \emph{cyclic} to mean a Hamilton cycle, i.e., the last and first object in the Gray code listing also differ by a flip.

In some of these examples we see that a flip operation can change many, possibly all, entries of the string representation of an object, so the question of what constitutes a `small change' is up to debate (but can sometimes be justified via polytopes; see Sections~\ref{sec:popo} and~\ref{sec:01}).
For example, the directed edges in Figure~\ref{fig:flip}~(e) correspond to cyclic left shifts of the permutations, which result in a change of every entry of the permutation.
Similarly, the edges of Kneser graphs as shown in Figure~\ref{fig:flip}~(c) correspond to disjoint sets, i.e., a flip exchanges all elements.
It turns out that such flip operations often lead to very interesting Gray codes.
Insisting that a Gray code may only change a constant number of positions in the string, which is prevalent in many papers, is a misguided paradigm.
One reason is that the amount of change is entirely dependent on the representation of the object.
For example, we will encounter Gray codes for permutations that use cyclic shifts of substrings, and while these shifts may change many entries of the permutation, they change only a single entry in the corresponding inversion vector.
Another reason is that depending on the available hardware, these `not-so-small' changes may actually be implemented faster.

Flip graphs allow us to apply graph-theoretic tools to tackle the Gray code problem.
Moreover, they call for investigating all kinds of other interesting graph-theoretic parameters beyond Hamiltonicity, such as girth, diameter, connectivity, colorability, independence number etc.

In some Gray code problems, the admissible flip operations depend on the recent history of previous flips, and in this case there is no underlying `static' flip graph, which typically makes these problems hard.
Examples of this are the order-preserving paths of Felsner and Trotter and Beckett-Gray codes (see Sections~\ref{sec:weight} and~\ref{sec:beckett}).

\subsection{Necessary conditions}

There are a few simple graph-theoretic conditions that a graph~$G$ has to satisfy in order to have a Hamilton path or cycle.
Trivially, $G$ needs to be connected.
Furthermore, if $G$ has vertices of degree~1, then a Hamilton cycle cannot exist, and if $G$ has more than two vertices of degree~1, then a Hamilton path cannot exist.
Also, if $G$ has precisely two vertices of degree~1, then any Hamilton path in~$G$ must start and end at these two vertices.
Moreover, if $G$ is bipartite with partition classes of different sizes, then a Hamilton cycle cannot exist, and if the size difference between the partition classes is more than~1, than a Hamilton path cannot exist.
This condition generalizes straightforwardly to $k$-partite graphs~$G$ for $k\geq 2$, where one compares the size of one partition class against the sum of the sizes of all others.
We refer to this condition on~$G$ as \emph{balancedness constraint}.
As we shall see, for most interesting flip graphs~$G$, the aforementioned connectivity, minimum degree and balancedness constraints are the only obstacles for Hamiltonicity, i.e., once $G$ satisfies these three constraints, then it does have a Hamilton path or cycle.
This could be considered a `meta-conjecture' for most interesting flip graphs (a notable counterexample is mentioned in~\cite{MR1814541}).

\subsection{Stronger Hamiltonicity notions}
\label{sec:hamiltonicity}

For constructing Gray codes, it is often helpful to establish stronger Hamiltonicity properties of the corresponding flip graph.
Specifically, a graph is called \emph{Hamilton-connected}, if it admits a Hamilton path between every pair of vertices.
Moreover, a bipartite graph is called \emph{Hamilton-laceable}, if it admits a Hamilton path between every pair of vertices from the two partition classes.
For example, the hypercube and the permutahedron for $n\geq 4$ are Hamilton-laceable (see~\cite{MR527987} and~\cite{MR683982}, respectively), and the Johnson graph is Hamilton-connected~\cite{MR1298971,MR1301949}.
There are many more strengthened Hamiltonicity properties in addition to the aforementioned two, some of which we will encounter in Section~\ref{sec:graphs}; see also~\cite{MR2500822}.
Establishing these stronger properties yields more information about the flip graph, and maybe counter-intuitively, it often makes the proof of Hamiltonicity \emph{easier}, and \emph{not harder}.
This is because the inductive statements involving these stronger notions are more powerful and flexible.

\subsection{General Gray code results}
\label{sec:general}

As mentioned before, a Gray code corresponds to a Hamilton path or cycle in the corresponding flip graph~$G$.
More generally, a \emph{$k$-Gray code} is a sequence of vertices of~$G$ that contains every vertex exactly once, such that any two consecutive vertices in the sequence have distance at most~$k$ in~$G$.
Given a graph~$G$, the \emph{$k$th power of~$G$}, denoted $G^k$, is the graph with the same vertex set as~$G$, connecting any two vertices that have distance at most~$k$ in~$G$ with an edge.
With this definition, a $k$-Gray code corresponds to a Hamilton path in~$G^k$, and a cyclic $k$-Gray code corresponds to a Hamilton cycle in~$G^k$.

Sekanina~\cite{MR0140095} proved that the cube~$G^3$ of any connected graph~$G$ has a Hamilton cycle.
Equivalently, his result asserts that any connected graph~$G$ admits a cyclic 3-Gray code.
If a flip graph is not just connected, but 2-connected, even more can be said.
Specifically, Fleischner~\cite{MR332573} proved that the square~$G^2$ of any 2-connected graph~$G$ has a Hamilton cycle.
Equivalently, any 2-connected graph~$G$ admits a cyclic 2-Gray code.
Karaganis~\cite{MR230645} strengthened Sekanina's result, by showing that for any connected graph~$G$, the cube~$G^3$ is Hamilton-connected.
Similarly, Chartrand, Hobbs, Jung, Kapoor, and Nash-Williams~\cite{MR345865} strengthened Fleischner's theorem, by showing that for 2-connected~$G$, the square~$G^2$ is Hamilton-connected.
Alstrup, Georgakopoulos, Rotenberg and Thomassen~\cite{MR3775895} developed an algorithm for computing such Hamilton cycles and paths in~$G^2$, running in time and space~$O(|E|)$, where $|E|$ is the number of edges of~$G$.
Sekanina's aforementioned result about~$G^3$ is also algorithmic, with the same bounds, and a cycle can be obtained by traversing any spanning tree of~$G$ in prepostorder; see~\cite[Sec.~7.2.1.6]{MR3444818}.
To complement these results, it is known that if~$G$ is not 2-connected, then the Hamiltonicity problem for~$G^2$ is NP-complete~\cite{MR522906}.

These results, which deserve to be better known in the Gray code community, impose some criteria against which we should compare any Gray code results obtained for particular combinatorial objects.
Put differently, as soon as the flip graph of a Gray code problem is connected or 2-connected, we are immediately guaranteed a 3-Gray code or 2-Gray code, respectively, by the results from before, so any further work on the problem should either focus on obtaining 2- or 1-Gray codes, or on conceptual simplicity or algorithmic efficiency for the particular problem.

\subsection{Algorithmic questions}

Of course, the linear bound~$O(|E|)$ in the algorithms discussed in the previous section is not very satisfying in the context of Gray codes.
This is because the size of a flip graph for a particular Gray code problem is often exponential in the size of the parameter that specifies the objects.
For example, the $n$-cube has $2^n$ vertices and $2^{n-1}n$ edges, or the permutahedron for~$[n]$ has $n!$ vertices and $n!(n-1)/2$ edges.
In those situations, we will much prefer Gray code algorithms that work in time and space that are polynomial in~$n$.
The word `time' here refers to the time per generated object, or equivalently, to the delay between two consecutive objects.
The ultimate goal from an algorithmic point of view is to develop a Gray code algorithm with constant delay, i.e., each new object is produced in constant time.
Such algorithms are sometimes called \emph{loopless} or \emph{loop-free}, a term coined by Ehrlich~\cite{MR0366085}.
Another interesting class are \emph{constant amortized-delay} or \emph{constant amortized-time (CAT)} algorithms, which produce each new object in constant time on average (not necessarily in the worst case).
In addition to being fast, a Gray code algorithm ideally requires only linear space to store the current object and possibly some auxiliary data structures.
Information-theoretic trade-offs concerning the space- and memory-efficient computation of Gray codes are explored in the literature under the term \emph{quasi-Gray codes}, a direction initiated by Fredman~\cite{MR0483669}.

Apart from generating a Gray code listing, there are two closely related algorithmic tasks of interest, namely \emph{ranking} and \emph{unranking}.
Ranking means to compute, for a given combinatorial object, its position in the Gray code listing.
On the other hand, unranking means to compute, for a given position, the object at this position in the Gray code listing.
Clearly, efficient unranking can be used for efficient random generation of the objects.

In this survey, we comment on the aforementioned algorithmic issues for several Gray codes problems, but we do not delve into low-level implementation details.
Knuth covers a large number of these algorithms in his book~\cite{MR3444818}, which also contains numerous challenging research problems and interesting background information; see also~\cite{knuth_pre8A}.
Furthermore, Arndt's book~\cite{DBLP:books/Arndt2011} contains a vast number of recent C++ implementations and experimental comparisons of many of the Gray code algorithms discussed here.
Lastly, many of these algorithms are available on the Combinatorial Object Server~\cite{cos} for experimentation and download, an online project that was first launched by Frank Ruskey.

\subsection{Posets and polytopes}
\label{sec:popo}

Interestingly, many flip graphs can be equipped with a poset structure.
For example, the hypercube is the cover graph of the Boolean lattice, the permutahedron is the cover graph of the weak Bruhat order of the symmetric group~$S_n$, and the associahedron is the cover graph of the Tamari lattice (recall Figure~\ref{fig:flip}~(a)+(d)+(f)).
Moreover, many flip graphs are skeleta of high-dimensional polytopes, including the aforementioned hypercube, permutahedron and associahedron.
These observations link the Gray code problem to order theory and geometry, enabling us to apply tools and insights from those areas.
For example, Balinski's theorem yields information about the connectivity of the flip graph based on the dimension of the polytope.
Also, the geometric viewpoint can be used to justify the choice of `small change' condition used in a Gray code for a particular set of objects: The change is small if and only if it corresponds to walking along an edge of the polytope (see Section~\ref{sec:01} below).
Furthermore, this allows broad generalization of the classes of flip graphs under investigation, and unification of several known Gray code results through a common lens:

\begin{itemize}[leftmargin=5mm, noitemsep, topsep=3pt plus 3pt]
\item
The hypercube, permutahedron and associahedron arise as cover graphs of quotients of lattice congruences of the symmetric group~$S_n$~\cite{MR3221544}, and the corresponding polytopes are known as \emph{quotientopes}~\cite{MR3964495,MR4584712}.
This is a family of $2^{2^{n(1+o(1))}}$ many flip graphs, and there is a simple greedy algorithm for computing a Hamilton path in each of them~\cite{MR4344032}, which yields a unified description of several known Gray codes for bitstrings, permutations, binary trees, set partitions, rectangulations etc.~\cite{MR4391718}.

\item
\emph{Graph associahedra}, introduced by Carr and Devadoss~\cite{MR2239078,MR2479448} and Postnikov~\cite{MR2487491}, are another large class of polytopes that generalize hypercube, permutahedron and associahedron.
They yield a family of $2^{n^2/2(1+o(1))}$ many flip graphs, most of which have a Hamilton cycle~\cite{MR3383157,MR4415121}.

\item
The Johnson graph and the flip graph of spanning trees of a graph under edge exchanges (recall Figure~\ref{fig:flip}~(b)+(g)) arise as skeleta of \emph{matroid base polytopes}.
In fact, these polytopes are part of the much larger family of \emph{0/1-polytopes} discussed in the next section, which are all Hamilton-connected or Hamilton-laceable~\cite{MR762893}, and a Hamilton path on them can be computed by a simple greedy algorithm.

\item
Flip graphs on acyclic orientations of a graph induced by flipping the orientation of a directed cycle or a directed cut can be equipped with a distributive lattice structure~\cite{propp_2021}, and they are realized by distributive polytopes that are part of the large family of \emph{alcoved polytopes}~\cite{MR2352704,MR2727459}.

\item
Flip graphs on acyclic orientations of a graph induced by flipping the orientation of a single arc also have a poset structure, and they can be realized as polytopes~\cite{MR1609342}.
These posets sometimes form a lattice and admit lattice congruences that can also be realized as polytopes~\cite{MR4490908}, generalizing the aforementioned quotientopes.
\end{itemize}

\openproblem{1}{The question when Hamilton paths or cycles exist for the latter two classes of flip graphs is still very much open in general (see also Problem~P56).}

\subsection{0/1-polytopes}
\label{sec:01}

Given any set of bitstrings~$X\seq\{0,1\}^n$, the convex hull of~$X$ in~$\mathbb{R}^n$ is a \emph{0/1-polytope}.
If $X$ encodes a set of combinatorial objects, then the edges of this polytope `automatically' capture small change operations between the objects.
In the simplest case, if we take $X=\{0,1\}^n$, we obtain the hypercube, with edges corresponding to flipping a single bit.
Furthermore, if we take $X$ as characteristic vectors of bases of a matroid, we obtain the matroid base polytope, whose edges correspond to base exchanges.
If we take $X$ as characteristic vectors of matchings of a graph, we obtain the matching polytope~\cite{MR371732}, whose edges correspond to alternating paths/cycles.
In particular, perfect matchings of the complete bipartite graph~$K_{n,n}$ are in bijection with permutation matrices, and their convex hull yields the Birkhoff polytope, whose edges correspond to multiplication with a permutation cycle.
Furthermore, if we take $X$ as characteristic vectors of independent sets or vertex covers of a graph, we obtain the stable set polytope and vertex cover polytope~\cite{MR510371}, respectively, whose edges correspond to a symmetric difference exchange on a connected subgraph.
As a last example, if we take $X$ as characteristic vectors of antichains or ideals of a poset, we obtain the chain polytope or order polytope~\cite{MR824105}, respectively, whose edges correspond to a symmetric difference exchange on a connected subposet.

In fact, Naddef and Pulleyblank~\cite{MR638286,MR762893} proved that every 0/1-polytope is either Hamilton-connected or a hypercube, which is known to be Hamilton-laceable.
An efficient algorithm for computing a Hamilton path on the skeleton of any 0/1-polytope was proposed by Merino and M\"utze~\cite{MR4720318}.
Their algorithm greedily modifies the shortest possible prefix of the current vertex from~$X$ and moves to another vertex with minimum Hamming distance from the current one.
The algorithm can start at an arbitrary vertex of the polytope, and it can be implemented efficiently by solving instances of the corresponding linear optimization problems over~$X$.
This approach yields a large number of efficient Gray code algorithms for a variety of objects based on matroids, graphs and posets, such as matroid bases, spanning trees, forests, matchings, maximum matchings, and cost-optimal matchings in general graphs, vertex covers, minimum vertex covers and cost-optimal vertex covers in bipartite graphs, antichains, maximum antichains and cost-optimal antichains in posets, as well as ideals and cost-optimal ideals of posets.
In each case, the small change operations used are given by the edges of the corresponding 0/1-polytope.

\subsection{Lov\'asz' conjecture and Cayley graphs}
\label{sec:lovasz}

Many natural Gray code problems yield flip graphs that are highly symmetric.
Specifically, many are vertex-transitive, such as the hypercube, permutahedron, Johnson graphs or Kneser graphs, and others are regular, such as the associahedron, or more generally, all graph associahedra.
\openproblem{2}{This connects the Gray code problem to a well-known and long-standing conjecture due to Lov\'asz~{\cite{MR0263646}}, which asserts that every connected vertex-transitive graph has a Hamilton path.}
A stronger version of the conjecture asserts that such a graph admits a Hamilton cycle, unless it is one of following five exceptional graphs: a single edge, the Petersen graph, the Coxeter graph, and the graphs obtained from the latter two by replacing every vertex by a triangle.
Thomassen conjectured that there are only finitely many such counterexamples (see~\cite{MR543656}).
Note that by a theorem of Mader~\cite{MR289343}, every vertex-transitive graph with degree~$k$ is $k$-connected, provided that it does not contain a clique of size~4, which can be seen as evidence for Lov\'asz' conjecture.
In stark contrast, Babai~\cite{MR1373683} conjectured that for some $\varepsilon>0$, there are infinitely many connected vertex-transitive graphs with $n$ vertices (even Cayley graphs) without cycles of length at least $(1-\varepsilon)n$.

There are several results establishing Lov\'asz' conjecture subject to conditions on~$n$, the number of vertices of the graph.
For example, vertex-transitive graphs with prime~$n$ are circulant graphs~\cite{MR211908}, and therefore a Hamilton cycle can be traced by following vertices with a fixed step size (modulo $n$) along the circulant vertex ordering.
Similarly, the conjecture has been confirmed if $n$ equals one of the following expressions, where $p$ and $q$ denote prime numbers: $2p$~\cite{MR561039}, $3p$~\cite{MR978524}, $4p$~\cite{MR2388379}, $5p$~\cite{MR677056}, $6p$~\cite{MR2548563}, $10p$~\cite{MR2904975} (with some additional assumptions), $p^2$ and $p^3$~\cite{MR821510}, $p^4$~\cite{MR1604701}, $p^5$~\cite{MR3336590}, $2p^2$~\cite{MR900940}, and $pq$~\cite{MR4328721}, in some cases in the stronger Hamilton cycle version.
Christofides, Hladk\'{y} and M\'{a}th\'{e}~\cite{MR3269902} proved that for every $\varepsilon>0$ there is an $n_0$ such that every vertex-transitive graph with $n\geq n_0$ vertices and degree at least~$\varepsilon n$ has a Hamilton cycle.
For more results on Lov\'asz' conjecture, we refer to the survey~\cite{MR2548567} by Kutnar and Maru\v{s}i\v{c}.

\openproblem{3}{Cayley graphs are a major class of vertex-transitive graphs, and are thus also an important special case of Lov\'asz' conjecture.}
The question for Hamilton cycles in Cayley graphs was first raised by Rapaport-Strasser~\cite{MR111702}.
For a finite group~$G$ and a set of generators~$X$, the \emph{Cayley graph} $\Gamma(G,X)$ has vertex set~$G$ and edges $(g,gx),(g,gx^{-1})$ for all $g\in G$ and $x\in X$.
Every two vertices in this graph that are connected by an edge also have an edge in the opposite direction, so it can be considered as an undirected graph.
The \emph{Cayley digraph} $\vv{\Gamma}(G,X)$ has the vertex set~$G$ and directed edges $(g,gx)$ for all $g\in G$ and $x\in X$ (i.e., inverse operations must not be used).
Note that if $X^{-1}:=\{x^{-1}\mid x\in X\}=X$, then the distinction between these two graphs is irrelevant.
This happens for example if $X$ is a set of involutions, i.e., $x^2=1$ for all $x\in X$.
Also note that if $G=S_n$ is the symmetric group and $x=\sigma_k\in X$ is the operation that cyclically \emph{left-shifts} the first $k$ entries of the permutation, then for every directed edge~$(\pi,\pi x)=:(\pi,\pi')$, $\pi\in S_n$, of~$\vv{\Gamma}(G,X)$, we have $\pi'(i)=\pi(x(i))=\pi(\sigma_k(i))$ for all $i=1,\ldots,n$, i.e., the permutation~$\pi'$ is obtained from~$\pi$ by cyclically \emph{right-shifting} the first $k$~entries.
In this case, it is maybe more intuitive to consider the Cayley graph $\vv{\Gamma}(G,X^{-1})$, which is obtained from $\vv{\Gamma}(G,X)$ by reversing the orientation of all edges; see Figure~\ref{fig:flip}~(e).

There is a substantial number of cases for~$G$ and~$X$ for which $\Gamma(G,X)$ or $\vv{\Gamma}(G,X)$ are known to have a Hamilton cycle.
One instance is the Chen-Quimpo theorem~\cite{MR641233}, which asserts that if $G$ is an abelian group such that the Cayley graph $\Gamma(G,X)$ is connected and has degree at least~3, then it is Hamilton-connected if it is not bipartite and Hamilton-laceable if it is bipartite.
Conway, Sloane and Wilks~\cite{MR1032382} showed that if~$G$ is a direct product $G=G_1\times\cdots\times G_k$ of irreducible Coxeter groups, and if $X_i$ is the canonical basis for~$G_i$ for all $i=1,\ldots,k$, then the Cayley graph~$\Gamma(G,X)$ for the basis $X=X_1\times\cdots\times X_k$ has a Hamilton cycle.
For more results on the Hamiltonicity of Cayley graphs, see Sections~\ref{sec:rho-sigma} and~\ref{sec:perm-restr} and the surveys~\cite{MR641232,MR762322,MR1405010,MR2548568,MR4027130} (see also~\cite{MR2946084} for another relevant result on 3-regular Cayley graphs).

\subsection{Genlex ordering}

As mentioned before, this survey focuses on Gray codes and not on the broader topic of efficient combinatorial generation.
For this reason, we never discuss lexicographic or colexicographic orderings.
Put differently, for those orderings there is no corresponding flip operation and hence no flip graph.
Nevertheless, there is a generalization of lexicographic and colexicographic orderings that is very relevant in the context of Gray codes, namely genlex ordering, a term coined by Walsh~\cite{MR2062208}.
A listing of strings over some alphabet is in \emph{genlex order}, if all strings with the same prefix (or suffix) appear consecutively.
Such listings are also sometimes referred to as \emph{prefix-partitioned} (or \emph{suffix-partitioned})~\cite{MR2577686}.
This generalizes \emph{graylex order} introduced by Chase~\cite{MR995888}, in which the symbols after the common prefixes (or before the common suffixes) appear in monotonic order, i.e., increasingly or decreasingly.
Note that in lexicographic or colexicographic order, they always appear in increasing order.
We shall see that many different classes of objects admit genlex Gray codes.
This is significant, as any genlex listing can be constructed by a simple greedy algorithm; see Section~\ref{sec:greedy}.

\section{Binary strings}
\label{sec:bits}

In this section we discuss Gray codes for listing all $2^n$ bitstrings of length~$n$, subject to various additional conditions.
We also discuss Gray codes for subsets of bitstrings that arise from structural constraints in the hypercube~$Q_n$, from forbidding certain substrings, or from considering equivalence classes of strings.
Most Gray codes in this section flip a single bit in each step, i.e., they correspond to paths or cycles in~$Q_n$ or subgraphs of it.
In Sections~\ref{sec:factor-avoid}, \ref{sec:neck} and~\ref{sec:flip-swap} we also encounter 2- and 3-Gray codes, which flip at most 2 or 3 bits in each step, which corresponds to paths or cycles in~$Q_n^2$ or~$Q_n^3$.
Another exception is Section~\ref{sec:shift-transpose}, where we consider Gray codes that use cyclic shifts or transpositions of bits.
Note that bitstrings of length~$n$ can be interpreted as characteristic vectors of subsets of~$[n]$, so the Gray codes discussed in this section can be seen as listings of subsets of the ground set~$[n]$.
Flipping a single bit in a bitstring corresponds to adding or removing a single element in the corresponding subsets.

For any string~$x$ and any integer~$\ell\geq 0$, we write~$x^\ell$ for the $\ell$-fold repetition of~$x$.

\subsection{The binary reflected Gray code}

\begin{figure}
\centerline{
\setlength{\tabcolsep}{3pt}
\renewcommand{\arraystretch}{1.5}
\begin{tabular}{ccccccccccccccccc}
(a) & (b) & (c) & (c') & (d) & (e) & (e') & (f) & (g) & (h) & (i) & (i') & (j) & (j') & (k) & (k') & (l) \\[-5mm]
\raisebox{-\height}{\includegraphics{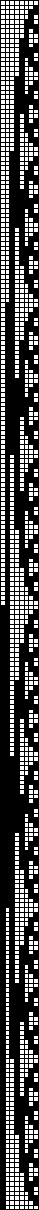}} &
\raisebox{-\height}{\includegraphics{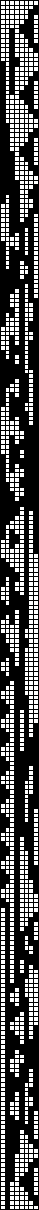}} &
\raisebox{-\height}{\includegraphics{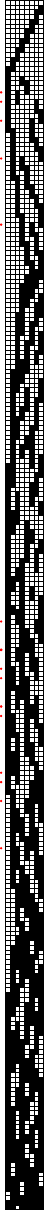}} &
\raisebox{-\height}{\includegraphics{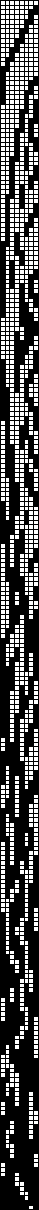}} &
\raisebox{-\height}{\includegraphics{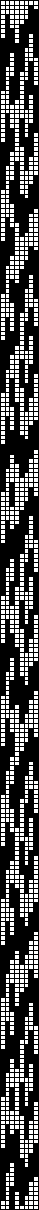}} &
\raisebox{-\height}{\includegraphics{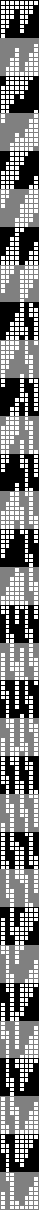}} &
\raisebox{-\height}{\includegraphics{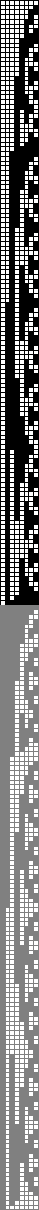}} &
\raisebox{-\height}{\includegraphics{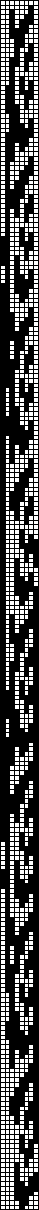}} &
\raisebox{-\height}{\includegraphics{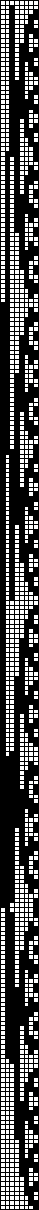}} &
\raisebox{-\height}{\includegraphics{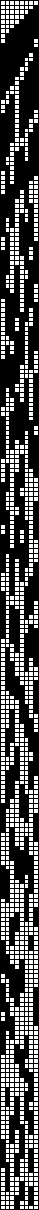}} &
\raisebox{-\height}{\includegraphics{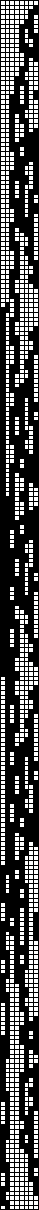}} &
\raisebox{-\height}{\includegraphics{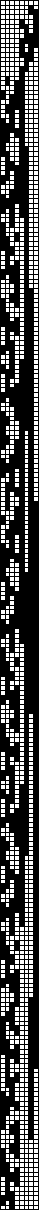}} &
\raisebox{-\height}{\includegraphics{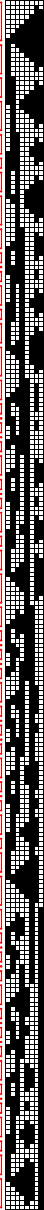}} &
\raisebox{-\height}{\includegraphics{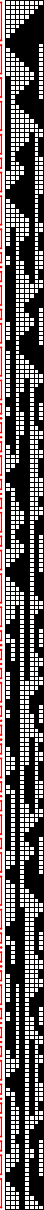}} &
\raisebox{-\height}{\includegraphics{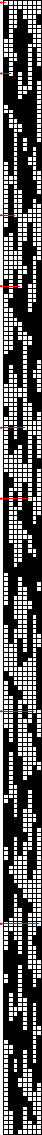}} &
\raisebox{-\height}{\includegraphics{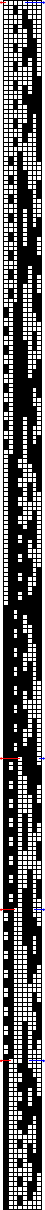}} &
\raisebox{-\height}{\includegraphics{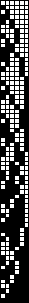}}
\end{tabular}
}
\caption{Different binary Gray codes (1-bits black, 0-bits white):
(a)~BRGC $\Gamma_8$;
(b)~balanced;
(c)~monotone;
(c')~order-preserving;
(d)~long-run;
(e)~8-antipodal;
(e')~128-antipodal;
(f)~non-local;
(g)~trend-free;
(h)~Beckett;
(i)~$K_8$-code;
(i')~$Q_3$-code;
(j)~symmetric chains;
(j')~Greene/Kleitman symmetric chains;
(k)~single-track (length~240);
(k')~2-track and 8-symmetric;
(l)~Savage-Shields-West (here $n=6$, unlike for the other codes where $n=8$).
}
\label{fig:bits1}
\end{figure}

The \emph{binary reflected Gray code (BRGC)}~\cite{gray_1953} lists all bitstrings of length~$n$ such that any two consecutive strings differ in a single bit, including the last and first bitstring.
The BRGC is defined inductively as $\Gamma_0:=\varepsilon$ and $\Gamma_{n+1}:=0\Gamma_n,1\rev(\Gamma_n)$ for $n\geq 0$, where $\varepsilon$ denotes the empty sequence and $\rev()$ denotes sequence reversal; see Figure~\ref{fig:bits1}~(a).
In words, the listing for strings of length~$n+1$ is obtained by prepending~0 to all strings in the listing for~$n$, followed by the reversed listing in which all strings are prepended~1.
The idea of this code was known long before Gray, and Knuth~\cite[Sec.~7.2.1.1]{MR3444818} tells more about its history.
In particular, the BRGC is useful for solving the Tower of Hanoi puzzle~\cite{MR3026271}, the Chinese rings puzzle, and to construct Brunnian links~\cite[Sec.~5.2.1]{ruskey_2003}.
The BRGC has very simple ranking/unranking functions, and many other fascinating properties.
For example, the BRGC ordering is genlex by definition, and genlex Gray codes for several other interesting combinatorial objects can be obtained by considering sublists of the BRGC; see Sections~\ref{sec:factor-avoid}--\ref{sec:flip-swap}.
Similarly, Tang and Liu~\cite{MR0349274} showed that in the sublist of the BRGC induced by bitstrings with a fixed number of 1s, any two consecutive bitstrings differ in two bits (a transposition of~0 and~1); see Section~\ref{sec:comb-trans}.
Also, it was shown by Yuen~\cite{MR0366527} and Cavior~\cite{MR0378994} that if any two $n$-bit strings differ in $t>0$ bits, then they are at least $\lceil 2^t/3\rceil$ steps and at most $2^n-\lceil 2^t/3\rceil$ steps apart in the BRGC.
More interesting properties of the BRGC are established in the papers by Conder~\cite{MR1663855}, and by Agrell, Lassing, Str\"om, and Ottosson~\cite{1362904}.
Bitner, Ehrlich and Reingold~\cite{MR0424386} describe a loopless algorithm for generating the BRGC.
\openproblem{4}{For many of the Gray codes with additional constraints discussed in the following sections, no such efficient algorithms are known, and it would be very interesting to devise them.}

\subsection{Transition counts}
\label{sec:transition}

Since the discovery of the BRGC, there has been continued interest in constructing Gray codes with additional desirable properties that the BRGC does not have.
A problem that has received a lot of attention is to construct a Gray code where the number of flips are distributed as equally as possible over the $n$~bit positions; see Figure~\ref{fig:bits1}~(b).
Specifically, let $c_i$ denote the number of times that bit~$i$ is flipped along a given Gray cycle, where we number bits from right to left with $i=1,\ldots,n$.
Clearly, each transition count~$c_i$ must be an even number, and they sum up to $2^n$, so we require that each~$c_i$ is approximately equal to~$2^n/n$.
Formally, in a \emph{balanced Gray code}, we require that $|c_i-2^n/n|<2$ for $i=1,\ldots,n$.
In particular, if $n$ is a power of 2, then $c_1=\cdots=c_n=2^n/n$.
Note that in the BRGC, we have $c_i=2^{n-i}$ for $i=1,\ldots,n-1$ and $c_n=2$, i.e., it is very unbalanced.
A balanced Gray code was first constructed by Tibor Bakos around~1955, and his proof is published in \'Ad\'am's book \cite[p.~28--37]{MR0242572}.
Apparently, this proof got forgotten, and the entire problem was rediscovered and resolved many years later.
In particular, several heuristics for constructing balanced Gray codes were suggested by Tootill~\cite{tootill_1956}, Vickers and Silverman~\cite{MR565503}, Ludman~\cite{1094886}, and by Robinson and Cohn~\cite{MR600765}.
Eventually, Wagner and West~\cite{MR1124876} proved the existence of balanced Gray codes when $n$ is a power of 2, attributing the problem to Stanley, and subsequently Bhat and Savage~\cite{MR1410880} solved the general case.
A very concise general solution is described in Knuth's book~\cite[Sec.~7.2.1.1]{MR3444818}.

A variant of balanced Gray codes are \emph{exponentially balanced Gray codes}, in which each transition count~$c_i$ is one of two consecutive powers of~2.
Note that if $n$ is a power of~2, then this is the same problem as before.
Van Zanten and Suparta~\cite{MR2113073} constructed exponentially balanced Gray codes for any $n\geq 1$, answering a conjecture posed in the aforementioned paper by Wagner and West, and a simplified construction was later given by Suparta~\cite{MR2176536}.

The most general version of this problem is to ask for an $n$-bit Gray code with arbitrary prescribed transition counts $c_1,\ldots,c_n$.
Clearly, we may assume without loss of generality that $c_1\leq c_2\leq \cdots\leq c_n$, and we have to require that each $c_i$ is even, that $\sum_{i=1}^n c_i=2^n$, and that $\sum_{i=1}^k c_i\geq 2^k$ for each $1\leq k\leq n$.
\openproblem{5}{The problem is whether every sequence $c_1,\ldots,c_n$ of integers satisfying those necessary conditions allows constructing a cyclic $n$-bit Gray code with exactly these transition counts.}
Observe that solving this problem for one particular value of~$n$ also solves it for all smaller values.
On the other hand, Potapov~\cite{MR2978613} positively resolved the problem for all $n>N:=\lceil 2^{2^{\mu/2}(4+\log_2 \mu)}\rceil$ where $\mu:=4(2^{35}-2^4)/31$, \emph{conditional} on the assumption that it has a positive solution for~$n=N$ (and therefore also for all~$n<N$).
Consequently, settling the base case~$n=N$ for his inductive argument is a very interesting open problem, but is infeasible to tackle solely with computer help.
Already Gilbert~\cite[Thm.~IV]{MR0094273}, using a lemma of Shannon~\cite{MR0029860}, showed how to construct a Gray \emph{path} with prescribed transition counts of the form $1,c_2,\ldots,c_n$ (i.e., the first bit flips exactly once), if they satisfy the corresponding obvious necessary conditions.

\subsection{Weight-monotonicity}
\label{sec:weight}

The \emph{weight} of a bitstring is defined as the number of~1s in it.
Savage and Winkler~\cite{MR1329390} considered Gray codes for which the weight of bitstrings increases (almost) monotonically along the code.
Specifically, a \emph{monotone} Gray path starts at~$0^n$ and once it has visited a bitstring of weight~$k$, it never returns to a bitstring of weight~$k-2$.
The last bitstring on the path will either have weight~$n$ if $n$~is odd or weight~$n-1$ if $n$~is even.
Savage and Winkler constructed such a monotone path for any $n\geq 1$, and their solution is shown in Figure~\ref{fig:bits1}~(c).
An even stronger requirement was proposed by Felsner and Trotter~\cite{MR1350586}; see also~\cite[Problem~483]{MR2548549}.
They ask for an \emph{order-preserving} Gray path $x_1,\ldots,x_N$ with $N=2^n$, which starts at~$0^n$ and if $x_i$ is obtained from~$x_j$ by a flip $1\rightarrow 0$, then $i\leq j+1$.
It is easy to check that any order-preserving path must also be monotone.
However, the Savage-Winkler paths are not order-preserving, and for $n=8$ the bitstrings violating the order-preserving property are marked in Figure~\ref{fig:bits1}~(c) with a small dot on the left.
Computer experiments show that for $n=1,\ldots,7$, the number of order-preserving Gray paths is $1,1,1,1,10,123,\num{1492723}$ [OEIS~A259839], and solutions have also been found for $n\in\{8,9,10\}$ (by Trotter and M\"utze, independently), one of which is shown in Figure~\ref{fig:bits1}~(c').
For comparison, the number of monotone paths for $n=1,\ldots,6$ is $1,1,1,1,14,\num{46935}$.
In all those counts, paths that differ only in renaming the coordinate directions are factored out and not counted separately.
\openproblem{6}{It is a challenging open problem whether order-preserving paths exist for $n\geq 11$.}
As a partial result, Felsner and Trotter~\cite{MR1350586} constructed an order-preserving Gray path of length $N=2^{n-1}(1+o(1))$.
Moreover, it was shown by Bir\'o and Howard~\cite{MR2525360} that for all $n\geq 11$, there is an order-preserving path that traverses all bitstrings with weights $0,1,2,3$ (and some bitstrings with weight~4, but no other weights), so it has length $N=\Theta(n^3)$.

\subsection{Antipodal Gray codes}

In a \emph{$t$-antipodal Gray code}, we require that the complement of any bitstring is visited either exactly $t$~steps before or later.
In principle, $t$~can be any integer in the interval~$[n,2^{n-1}]$, but~$t$ and~$n$ must have the same parity.
Killian and Savage~\cite{MR2047769} attribute this problem to Snevily, and they consider the case when~$t=n$.
Specifically, they proved that no solution exists for any odd~$n\geq 5$ or for~$n=6$, and presented a constructive solution when $n$~is a power of~2, which is shown in Figure~\ref{fig:bits1}~(e).
\openproblem{7}{They mention as an open problem the question whether an $n$-antipodal $n$-bit Gray code exists for even $n\geq 10$.}
Chang, Eu, and Yeh~\cite{MR2317471} showed that a $t$-antipodal $n$-bit code for odd~$n\geq 3$ exists if and only if~$t=2^{n-1}-1$.
Moreover, for even~$n\geq 6$ they constructed $t$-antipodal codes for all $t=2^{n-1}-2^k$ with $k$~odd and $1\leq k\leq n-3$, and for $t=2^{n-k}$ with $1\leq k\leq 3$.
Knuth~\cite[Ex.~49 in Sec.~7.2.1.1]{MR3444818} describes a simple solution for even $n\geq 2$ and $t=2^{n-1}$; see Figure~\ref{fig:bits1}~(e').
In Figures~\ref{fig:bits1}~(e) and (e'), any bitstring whose complement appears $t$ steps earlier is drawn shaded.

\subsection{Run lengths}

Another property of interest are \emph{run lengths} of a Gray code, i.e., the number of bitstrings visited between two flips of the same coordinate.
In Figure~\ref{fig:bits1}, the run lengths are the heights of maximal vertical blocks of either black or white squares.
We let $r(n)$ denote the maximum value $r$ such that an $n$-bit Gray code with all run lengths at least~$r$ exists.
\openproblem{8}{The exact value of $r(n)$ for $n=1,\ldots,7$ is $r(n)=1,2,2,2,4,4,5$, and it is unknown for $n\geq 8$~\mbox{\cite[Sec.~7.2.1.1]{MR3444818}}.}
An 8-bit Gray code with all run lengths at least~5 is shown in Figure~\ref{fig:bits1}~(d) (it is the same as in \cite[Fig.~34 in Sec.~7.2.1.1]{MR3444818}).
Clearly, we have $r(n)\leq n$, and Goddyn and Gvozdjak~\cite{MR2014514} proved the remarkable fact that $r(n)\geq \lceil n-2.001 \log_2(n)\rceil\geq n-3\log_2(n)=n-\cO(\log n)$ for $n\geq 2$, improving upon earlier results due to Evdokimov~\cite{MR625381} and Goddyn, Lawrence and Nemeth~\cite{MR982258}.
Note that in an $n$-antipodal $n$-bit Gray code, the run lengths within every pair of a block and its subsequent complemented block are equal to~$n$, but not necessarily across the boundaries of these pairs.
Codes with long run lengths have applications in image detectors~\cite{lawrence_mcclintock_1996} and for structured light~\cite{DBLP:conf/icra/KimRL08}.

\subsection{Non-locality}

In a \emph{non-local Gray code}, we require that among any $2^t$ consecutive bitstrings, strictly more than $t$ bit positions change, for any $2\leq t\leq n-1$; see Figure~\ref{fig:bits1}~(f).
This is impossible for $n=2,3,4$, but can be achieved for any $n\geq 5$, as was first shown by Mills~\cite{MR0151950}.
A simpler construction was later given by Ramras~\cite{MR1081842}.

\subsection{Trend-freeness}

In a \emph{trend-free Gray code}, the center of mass of 1-bits in each coordinate position is in the middle.
Equivalently, we require that $\sum_{i=1}^{2^n} i(-1)^{x_{i,k}}=0$ for all $k=1,\ldots,n$, where $x_{i,k}$ denotes the $k$th bit in the $i$th bitstring.
Such codes are useful for experiment design~\cite{doi:10.1080/00401706.1974.10489146, MR565503, cheng_1985}.
Note that the BRGC satisfies the trend-freeness conditions for all bit positions except the leftmost one ($k=n$).
Knuth \cite[Ex.~75 and~76 in Sec.~7.2.1.1]{MR3444818} gives a simple construction of trend-free Gray codes for all $n\geq 5$; see Figure~\ref{fig:bits1}~(g), improving upon an earlier result due to Cheng~\cite{cheng_1985}.

\subsection{Beckett-Gray codes}
\label{sec:beckett}

\emph{Beckett-Gray codes} are named after the Irish playwright Samuel Beckett, whose play `Quad' consists of a series of arrivals and departures of characters who appear on stage together in different combinations.
Whenever a character leaves, it is the one who has been on stage longest among all characters present.
Formally, the problem is to find an $n$-bit Gray code with the property that on every transition $1\rightarrow 0$, the 1-bit that disappears is the one that appeared first among all 1-bits.
There is no constraint with regards to the other transitions $0\rightarrow 1$.
We can think of the positions of 1s as being stored in a queue that operates on a first-in first-out basis.
Cooke, North, Dewar and Stevens~\cite{cooke_north_dewar_stevens_2016} found such codes with the help of computers for $n=2,5,6,7,8$, see Figure~\ref{fig:bits1}~(h), and showed that none exist for $n=3,4$.
Sawada and Wong~\cite{DBLP:journals/endm/SawadaW07} report on some computational advances and additional solutions found with computer help.
\openproblem{9}{It is an open problem whether Beckett-Gray codes exist for $n\geq 9$.}
This problem has a `dynamic' flavor, i.e., the way that the cycle can move through the $n$-cube depends on the recent history, similarly to the aforementioned order-preserving path problem by Felsner and Trotter, which makes both problems hard.

\subsection{Transition graphs}

We let $P_n$, $C_n$, $K_n$ and $K_{k,n-k}$ denote the path, cycle, complete graph and complete bipartite graph on the vertex set~$[n]$, where the partition classes of~$K_{k,n-k}$ have size~$k$ and~$n-k$.
We consider the \emph{transition sequence} of a Gray code, i.e., the cyclic sequence $\delta_1,\ldots,\delta_{2^n}$ of bit positions that change along the code.
For instance, starting at~$0000$ the BRGC $\Gamma_4$ shown in Figure~\ref{fig:flip}~(a) has the transition sequence $1213121412131214$ (the rightmost bit flips in every second step, and the leftmost bit flips only twice).
Such a transition sequence defines a graph~$G$ on the vertex set $[n]$, by adding all edges $\{\delta_i,\delta_{i+1}\}$ for $i=1,\ldots,2^n$, where indices are taken modulo~$2^n$.
We refer to a Gray code yielding a particular graph~$G$ as a \emph{$G$-code}.
Note that $G$ is always a connected graph, and the BRGC is a $K_{1,n-1}$-code.
The problem of which $n$-vertex graphs $G$ admit a $G$-code was first raised by Slater~\cite{slater_1979} for the case when~$G=P_n$.
\openproblem{10}{There is a $P_n$-code for $n\leq 6$, but none for $n=7$~{\cite{slater_1989}}, and it is open whether there is one for~$n\geq 8$.}
Darko Dimitrov, Petr Gregor and Borut Lu\v{z}ar (personal communication) ran a computer search that excludes the existence of a $P_n$-code also for the case~$n=8$.
It is open whether there is a $C_n$-code for $n\geq 8$ (there is one for $n=4$ and none for $n=3,6,7$, and none is known for $n=5$).
The papers~\cite{MR1377601} and~\cite{MR1935751} catalogue all trees~$G$ on up to $n=7$ vertices that do not admit a $G$-code.
Moreover, Bultena and Ruskey~\cite{MR1377601} proved that if $G$ is a tree with diameter~4, then there is a $G$-code, and if $G$ is a tree with diameter~3, then there is no $G$-code.
They also constructed $G$-codes such that $G$ is a subgraph of any multipartite graph.
Wilmer and Ernst~\cite{MR1935751} answered several open problems from the Bultena-Ruskey paper, by showing that there are $G$-codes such that $G$ is a tree with arbitrarily large diameter, and that there are $G$-codes such that $G$ is a subgraph of any two-dimensional grid.
Suparta and van~Zanten~\cite{MR2427745} constructed a $K_n$-code for any $n\geq 3$; see Figure~\ref{fig:bits1}~(i).
Moreover, Suparta~\cite{MR3716903} describes a $K_{k,n-k}$-code for any $1\leq k\leq n-1$.
Dimitrov, Dvo\v{r}\'ak, Gregor, and \v{S}krekovski~\cite{MR2974271} constructed an $n$-bit $G$-code such that $G$ is an induced subgraph of the $\lceil\log_2(n)\rceil$-cube.
In particular, if $n$ is a power of~2, then this is a $Q_{\log_2(n)}$-code; see Figure~\ref{fig:bits1}~(i').

\subsection{Avoiding and including other structures}

If we consider a Gray code as a Hamilton cycle in~$Q_n$ for~$n\geq 2$, we may ask for cycles that avoid or include certain other structures.

In this direction, Latifi, Zheng and Bagherzadeh~\cite{DBLP:conf/ftcs/LatifiZB92} showed that for any set~$F$ of at most~$n-2$ edges, there is a Hamilton cycle in~$Q_n$, $n\geq 2$, that avoids~$F$, and clearly the bound~$n-2$ is best possible.
Furthermore, Chan and Lee~\cite{MR1129389} showed that $Q_n$, $n\geq 3$, has a Hamilton cycle that avoids a given set~$F$ of at most $2n-5$~edges, provided that each vertex is incident to at least two edges not in~$F$, and the bound~$2n-5$ is best possible.
Moreover, Dimitrov, Dvo\v{r}\'ak, Gregor and \v{S}krekovski~\cite{MR2568029} showed that $Q_n$ has a Hamilton cycle that avoids a given perfect matching~$M$ if and only if~$Q_n\setminus M$ is connected.

Instead of forbidding edges, one may also consider the problem of forbidding vertices.
To this end, Fink and Gregor~\cite{MR2970496} proved that for any set~$F$ of at most~$\binom{n}{2}-2$ vertices of~$Q_n$, $n\ge 3$, there is a cycle of length at least~$2^n-2|F|$ that avoids all vertices in~$F$.
The bound~$\binom{n}{2}-2$ is best possible.
Note also that if $F$ contains only vertices of the same parity, then any extending cycle will also skip $|F|$ vertices of the opposite parity, so $2|F|$ vertices in total.

Dvo\v{r}\'ak~\cite{MR2178189} showed that any set of at most $2n-3$~edges that together form disjoint paths in~$Q_n$ can be extended to a Hamilton cycle, and this bound is best possible.
Note that any Hamilton cycle in~$Q_n$ can be partitioned into two perfect matchings.
Conversely, Kreweras~\cite{MR1374632} conjectured that every perfect matching in~$Q_n$ can be extended to a Hamilton cycle (via a second perfect matching).
More generally, Ruskey and Savage~\cite{MR1201997} asked whether any matching (not necessarily perfect) extends to a Hamilton cycle; see also~\cite[Problem~487]{MR2548549}.
Kreweras' conjecture was answered affirmatively by Fink~\cite{MR2354719}, and his beautiful and concise proof was later refined by Gregor~\cite{MR2510579}.
More recently, Fink~\cite{MR3936192} proved that every (not necessarily perfect) matching of~$Q_n$ extends to a cycle factor, i.e., to a collection of disjoint cycles that together visit all vertices.
\openproblem{11}{The original Ruskey-Savage conjecture is still open, and constitutes one of the most intriguing problems in this area.}

A \emph{symmetric chain} in~$Q_n$ is a path $x_k,\ldots,x_{n-k}$ where $x_i$ has weight~$i$ for all $i=k,\ldots,n-k$, and a \emph{symmetric chain decomposition} is a partition of all vertices of~$Q_n$ into symmetric chains.
This concept was pioneered by Greene and Kleitman~\cite{MR0389608} (see also \cite{MR0043115,MR205864,MR319772,MR450151}).
Streib and Trotter~\cite{MR3268651} constructed a Hamilton cycle in~$Q_n$ that contains a symmetric chain decomposition; see Figure~\ref{fig:bits1}~(j), where the chains are highlighted on the side.
Gregor, Mi\v{c}ka and M\"utze~\cite{MR4539109} showed that the aforementioned symmetric chain decomposition due to Greene and Kleitman can be extended to a Hamilton cycle of~$Q_n$; see Figure~\ref{fig:bits1}~(j').
A Hamilton cycle that extends a symmetric chain decomposition has the minimum number of \emph{peaks} and \emph{valleys}, which are triples of consecutive vertices $x_i,x_{i+1},x_{i+2}$ where $x_i$ and~$x_{i+2}$ have the same weight.
It is interesting to contrast this with the aforementioned monotone path due to Savage and Winkler, which has the maximum number of peaks and valleys.
The authors of~\cite{MR4539109} also conjectured that every symmetric chain decomposition of~$Q_n$ can be extended to a Hamilton cycle, but a counterexample in~$Q_6$ can be constructed easily:
We start with the Greene-Kleitman decomposition, and we consider the chains $C_1={-}{-}{-}{-}{-}{-}$, $C_2={-}{-}01{-}{-}$, $C_3={-}0011{-}$ and $C_4=000111$  of lengths 6, 4, 2, and 0, respectively, where 0s and 1s represent fixed bits, and $-$ represents a bit that changes from~0 to~1 when moving up the chain, with the flips occurring from left to right.
We remove the middle vertex of~$C_i$ and reconnect the remaining two pieces of the chain~$C_i$ via the middle vertex of~$C_{i+1}$ for $i=1,2,3$.
This yields a modified symmetric chain decomposition of~$Q_6$ in which the former middle vertex of~$C_1$, namely $111000$, is a chain of length~0, with the additional property that all of its neighbors are inner vertices of another chain, i.e., this length~0 chain cannot be connected to any other chain along a Hamilton cycle.

The results of Aubert and Schneider~\cite{MR676514} imply that all edges of~$Q_n$ can be partitioned into Hamilton cycles and possibly a perfect matching (see \cite{MR1096980}).
Saad and Schultz~\cite{2234} proved that $Q_n$ contains a cycle of length~$\ell$ for every even number~$\ell$ from~$4$ to~$2^n$.
Simmons~\cite{MR527987} first proved that $Q_n$ has a Hamilton path between any two prescribed vertices of opposite parity, i.e., it is Hamilton-laceable.
This was later generalized by Dvo\v{r}\'ak and Gregor~\cite{MR2326160}, who showed that any two prescribed vertices~$u$ and~$v$ of opposite parity can be joined by a Hamilton path containing at most $2n-4$ prescribed edges that together form disjoint paths, if none of the paths contains~$u$ or~$v$ internally or joins~$u$ and~$v$.
The bound~$2n-4$ is again best possible.
More results about cycles in the $n$-cube satisfying various additional constraints are surveyed by Xu and Ma~\cite{MR2500822}.
For some older results about the $n$-cube, we refer the reader to the survey~\cite{MR949280} by Harary, Hayes and Wu.

\subsection{The number of Gray codes}

Let $H(n)$ denote the number of Hamilton cycles in the $n$-cube, where we do not distinguish the starting vertex or direction of the cycle, but the names of the $n$~coordinate directions matter (factoring them out one needs to divide by~$n!/2$).
For $n=1,\ldots,6$ the exact values of~$H(n)$ are $1, 1, 6, \num{1344}, \num{906545760}, \num{35838213722570883870720}$ [OEIS~A066037], where the last value was computed by Haanp\"a\"a and \"Osterg{\aa}rd~\cite{MR3143701}.
The number of Hamilton cycles is closely related to the number of perfect matchings in the $n$-cube, which we denote by~$M(n)$ and which is known exactly for all~$n\leq 7$ [OEIS~A005271].
The results of Bregman~\cite{MR0327788} and Clarke, George, and Porter~\cite{MR1604993} imply that the quantity $m(n):=M(n)^{1/2^n}$ grows like $m(n)=(\frac{n}{e}(1+o(1))^{1/2}$.
Clearly, we have $H(n)\leq \binom{M(n)}{2}$, which shows that the quantity $h(n):=H(n)^{1/2^n}$ satisfies $h(n)\leq \frac{n}{e}(1+o(1))$.
The best known lower bound is $h(n)\geq \frac{n}{4e}+O(\log^2 n)=\frac{n}{e}(1/4+o(1))$, and it is due to Knuth~\cite[Ex.~47 in Sec.~7.2.1.1]{MR3444818}.
Earlier, slightly worse lower bounds were proved by Perezhogin and Potapov~\cite{MR1859856}, and by Feder and Subi~\cite{MR2483804}.
These results imply that $H(n)=M(n)^{2-o(1)}$, i.e., the number of Hamilton cycles essentially grows quadratically in the number of perfect matchings.
\openproblem{12}{The most interesting open question here is to tighten the bounds~$[\frac{1}{4},1]$ for the constant~$\lambda$ in $h(n)=\frac{n}{e}(\lambda+o(1))$.}

\subsection{Single-track Gray codes}
\label{sec:track}

\emph{Single-track Gray codes} were discovered by Hiltgen, Paterson, and Brandestini~\cite{DBLP:journals/tit/HiltgenPB96}.
In such a code, any two of the $n$~columns of bits differ only in a cyclic shift.
As before, any two consecutive rows/bitstrings differ in a single bit.
Clearly, the single-track property implies that the Gray code is balanced, i.e., each bit is flipped the same number of times.
Single-track Gray codes are useful for measuring angular positions, using only a single column of bits to encode the different positions, and $n$ appropriately rotated reading heads along this column; see Figure~\ref{fig:bits1}~(k), where the shifted columns are shaded and the shifts are indicated by horizontal bars.
In this example, 16~bitstrings are missing from the code.
Indeed, for~$n\geq 3$ the single-track property does not allow to visit all $2^n$~bitstrings, but only a subset of size~$N<2^n$.
In the aforementioned paper, Hiltgen, Paterson, and Brandestini show that $2n$ must divide~$N$ for such a code to exist.
Moreover, they provide constructions of single-track Gray codes of length~$N=kn$, for all $n\geq 4$ and every even $2\leq k\leq 2^{n-1-\lceil\sqrt{2(n-3)}\rceil}$.
Etzion and Paterson~\cite{MR1445874} constructed single-track Gray codes of length~$N=2^n-2n$ for the special case when $n$ is a power of~2.
Later, Schwartz and Etzion~\cite{MR1725126} proved that this is the maximum length of a cycle one can achieve under these conditions.
Their result shows that the single-track cycle of length~$2^8-2\cdot 8=240$ shown in Figure~\ref{fig:bits1}~(k) is in fact longest possible.
\openproblem{13}{Etzion~{\cite{MR2629553}} gives a short survey on the state-of-the art on single-track Gray codes, and mentions as the most interesting open problem the question whether for prime~$n\geq 3$ there is a single-track Gray code of length~$N=2^n-2$, which would be optimal.}
This has been tested and confirmed for all prime~$n\leq 19$.
A natural approach to construct such a Gray code is to find a Gray cycle through all $n$-bit necklaces except~$0^n$ and~$1^n$, which are aperiodic as $n$ is prime, so the number of such necklaces is $N/n$, selecting one representative from each necklace arbitrarily, in such a way that from the last representative a single flip leads to a shifted version of the first representative.
Applying the transition sequence of such a Gray cycle $n$ times, each time shifted appropriately, produces the desired single-track Gray code.
The open problem mentioned in Section~\ref{sec:neck} seeks a very similar Gray path through all $n$-bit necklaces (starting and ending in~$0^n$ and~$1^n$, respectively).

More generally, in a \emph{$t$-track Gray code}, the $n$~columns of bits can be partitioned into $t$ groups such that the columns in each group differ only in a cyclic shift.
Gregor, Merino and M\"utze~\cite{MR4747482} constructed a $t$-track Gray code of all $2^n$ bitstrings of length~$n$ (no bitstring is missed) for every $n$ that is a sum of~$t\geq 2$ powers of~2.
In particular, if $n$ is a sum of two powers of~2, then this is a 2-track Gray code; see Figure~\ref{fig:bits1}~(k').
Furthermore, these results also show that for every dimension~$n\geq 5$ there is a $t$-track Gray code with $t\leq \lfloor \log_2(n+1)\rfloor$, i.e., at most logarithmically many tracks.

\subsection{Symmetric Gray codes}
\label{sec:symm}

A Hamilton cycle~$(x_1,\ldots,x_N)$ of~$Q_n$, where $N=2^n$, is called \emph{$k$-symmetric}, if the mapping $x_i\rightarrow x_{i+N/k}$ for $i=1,\ldots,N$ is an automorphism of~$Q_n$ (indices are considered modulo~$N$).
In other words, when arranging the vertices in the order of the Hamilton cycle equidistantly on a circle and drawing the edges of~$Q_n$ as straight lines, then the resulting drawing has $k$-fold rotational symmetry.
Gregor, Merino and M\"utze~\cite{MR4747482} constructed a $k$-symmetric Hamilton cycle of~$Q_n$ for $k:=2^{\lceil \log_2 n\rceil}$, and they showed that this factor~$k$ is largest possible.
Clearly, if $n$ is a power of~2, then this optimum factor equals~$k=n$.
Such an 8-symmetric Hamilton cycle for $n=8$ is shown in Figure~\ref{fig:bits1}~(k'); the corresponding automorphism is $f(x_1,\ldots,x_8)=(x_2,x_3,x_4,\ol{x_1},x_6,x_7,x_8,\ol{x_5})$.

\subsection{Shifts and transpositions}
\label{sec:shift-transpose}

Feldmann and Mysliewietz~\cite{MR1397915} constructed a Gray code for all bitstrings of length~$n$ which in every step either flips the leftmost bit or shifts all bits cyclically one position to the left or right.

Savage, Shields and West~\cite{MR1951392} consider a Gray code for all bitstrings of length~$n$ which in every step either flips the leftmost bit or transposes two adjacent bits $01\leftrightarrow 10$.
Clearly, the listing must start and end with~$0^n$ and $1^n$, respectively, as these two only admit flipping the leftmost bit.
The authors constructed such a Gray code for all~$n>5$ for which $\binom{n+1}{2}$ is odd, and showed that it does not exist otherwise due to balancedness problems.
Figure~\ref{fig:bits1}~(l) shows their construction for~$n=6$.
Observe that such a Gray code corresponds to a Hamilton path in the cover graph of the poset of all lattice paths of length~$n$ with steps~$(+1,+1)$ or~$(+1,-1)$ (encoded by 1 and~0, respectively) ending at the fixed point~$(n,0)$ and ordered by the pointwise-below relation, which was introduced by Stanley in a slightly different context.
In this setting, the aforementioned results were considerably strengthened by Tasoulas, Manes and Sapounakis~\cite{MR4703269}, who showed that the cover graphs of various intervals in this poset all admit Hamilton cycles.

\subsection{Coils and snakes}

A \emph{coil-in-the-box} is an induced cycle in the $n$-cube, i.e., a cycle without chords.
Clearly, such a cycle cannot visit all $2^n$ vertices, and the problem is to determine the maximum length~$c(n)$ of a coil-in-the-box for each dimension~$n$.
For $n=1,\ldots,8$, the exact values of~$c(n)$ are known to be $0,4,6,8,14,26,48,96$ [OEIS~A000937], and an optimal solution for~$n=8$ is shown in Figure~\ref{fig:bits2}~(a).
In the literature one also encounters induced paths in the $n$-cube, which we call \emph{snakes-in-the-box} (in some papers, the term snake is used to refer to an induced cycle, not to an induced path).
The maximum length of a snake is denoted~$s(n)$ and for $n=1,\ldots,8$ the exact values of~$s(n)$ are $1,2,4,7,13,26,50,98$ [OEIS~A099155].
\openproblem{14}{Clearly, we have $s(n)\geq c(n)-2$, but this leaves open whether~$s(n)$ dominates~$c(n)$ or vice versa for larger values of~$n$.}
Kautz~\cite{DBLP:journals/tc/Kautz58} discovered coil-in-the-box and snake-in-the-box codes for single bit error detection in digital circuits, and he established simple lower and upper bounds for~$c(n)$.
He also determined the exact values of~$c(n)$ for~$n\leq 5$.
Subsequently, Davies~\cite{DBLP:journals/tc/Davies65}, Kochut~\cite{MR1376708}, and \"Osterg{\aa}rd and Pettersson~\cite{MR3281193} established the optimality of various coils for $n=6,7,8$, respectively, by exhaustive search.
The candidate of length~96 for~$n=8$ shown in Figure~\ref{fig:bits2}~(a) had already been found by Paterson and Tuliani~\cite{MR1616657}.
There was also much effort towards proving good asymptotic bounds for~$c(n)$ as $n$ grows.
The best known upper bound is due to Z\'emor~\cite{MR1479303}, who showed that
\begin{equation*}
c(n)\leq 1+2^{n-1}\frac{6n}{6n+\frac{1}{6\sqrt{6}}\sqrt{n}-7}\leq 2^{n-1}\Big(1-\frac{\sqrt{n}}{89}+O\big(\frac{1}{n}\big)\Big).
\end{equation*}
for all~$n\geq 1$.
This bound improves upon earlier results due to Douglas~\cite{MR0289334}, Deimer~\cite{MR815576}, and Snevily~\cite{MR1298987}.
As far as lower bounds are concerned, Evdokimov~\cite{evdokimov_1969} first established that $c(n)\geq \lambda 2^{n-1}$ for $\lambda:=c(8)/256=0.375$, i.e., there are coils that miss only a constant fraction of all vertices of the $n$-cube, a fact that is maybe not apparent from the first few values of~$c(n)$.
Wojciechowski~\cite{MR1010304} independently proved a similar bound, but with a slightly worse (i.e., smaller) constant.
The best known lower bound is due to Abbott and Katchalski~\cite{MR1143262}, who showed that one can take $\lambda:=77/128>0.601$ for~$n\geq 2$, following up on their earlier paper~\cite{MR953891}.
\openproblem{15}{Consequently, the most interesting open problem is whether $c(n)\leq \lambda 2^{n-1}$ for some constant~$\lambda<1$ (and large~$n$).}
Evdokimov summarizes many of these results in his survey~\cite{MR3474560} (in Russian).

\subsection{Factor-avoiding strings}
\label{sec:factor-avoid}

For any fixed bitstring~$\tau$, we consider all bitstrings of length~$n$ that \emph{avoid~$\tau$ as a factor}, i.e., that do not contain~$\tau$ as a contiguous substring.

\begin{figure}
\centerline{
\setlength{\tabcolsep}{2pt}
\renewcommand{\arraystretch}{1.5}
\begin{tabular}{cccccccccccccccccc}
(a) & (b) & (b') & (c) & (c') & (d1) & (d2) & (d3) & (e1) & (e1') & (e2) & (e2') & (e3) & (e4) & (e5) & (e6) & (e7) & (e8) \\[-5mm]
\raisebox{-\height}{\includegraphics{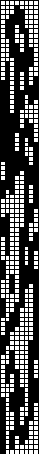}} &
\raisebox{-\height}{\includegraphics{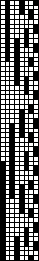}} &
\raisebox{-\height}{\includegraphics{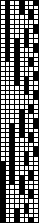}} &
\raisebox{-\height}{\includegraphics{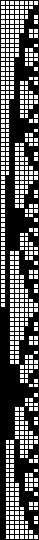}} &
\raisebox{-\height}{\includegraphics{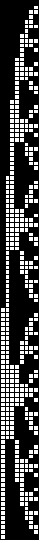}} &
\raisebox{-\height}{\includegraphics{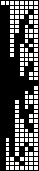}} &
\raisebox{-\height}{\includegraphics{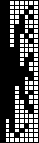}} &
\raisebox{-\height}{\includegraphics{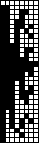}} &
\raisebox{-\height}{\includegraphics{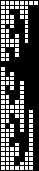}} &
\raisebox{-\height}{\includegraphics{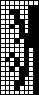}} &
\raisebox{-\height}{\includegraphics{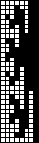}} &
\raisebox{-\height}{\includegraphics{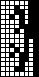}} &
\raisebox{-\height}{\includegraphics{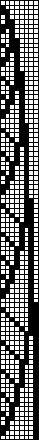}} &
\raisebox{-\height}{\includegraphics{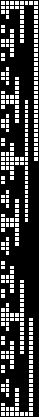}} &
\raisebox{-\height}{\includegraphics{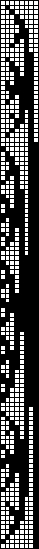}} &
\raisebox{-\height}{\includegraphics{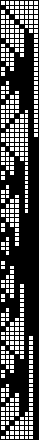}} &
\raisebox{-\height}{\includegraphics{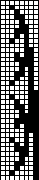}} &
\raisebox{-\height}{\includegraphics{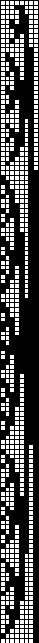}} \\
(e9') & (e10') & (e11') & (f) & (g) & & & & & & & & & & & & & \\[-5mm]  
\raisebox{-\height}{\includegraphics{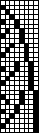}} &
\raisebox{-\height}{\includegraphics{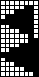}} &
\raisebox{-\height}{\includegraphics{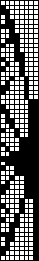}} &
\raisebox{-\height}{\includegraphics{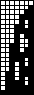}} &
\raisebox{-\height}{\includegraphics{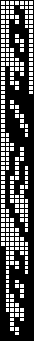}} & & & & & & & & & & & & &
\end{tabular}
}
\caption{Restricted binary Gray codes (1-bits black, 0-bits white):
(a)~coil-in-the-box (length~96);
(b)~11-avoiding (Fibonacci words);
(b')~11-avoiding circularly (Lucas words);
(c)~101-avoiding;
(c')~010-avoiding;
(d1)~necklaces (lex.\ largest);
(d2)~aperiodic necklaces (lex.\ largest);
(d3)~bracelets (lex.\ largest);
(e1)~necklaces (lex.\ smallest);
(e1')~unlabeled necklaces;
(e2)~Lyndon words (=lex.\ smallest aperiodic necklaces);
(e2')~unlabeled Lyndon words;
(e3)~weight $\leq 3$;
(e4)~100-avoiding;
(e5)~$\leq 6$ inversions w.r.t.~$0^*1^*$;
(e6)~$\leq 1$ transpositions w.r.t.~$0^*1^*$;
(e7)~$\leq 00100101$;
(e8)~$\leq$ their reversal (lex.\ smallest neckties);
(e9')~$\{11,101\}$-avoiding;
(e10')~$\{100,101\}$-avoiding;
(e11')~1-decreasing;
(f)~necklaces (repr.\ in BRGC order);
(g)~neckties (lex.\ smallest).
%
}
\label{fig:bits2}
\end{figure}

We say that a string~$\tau$ is \emph{bifix-free} if none of its proper prefixes equals one of its proper suffixes.
Squire~\cite{MR1385319} showed that if $\tau$ is bifix-free, or if $\tau=\sigma\sigma\cdots\sigma$ for some bifix-free string $\sigma$, then for all $n\geq 1$ there is a Gray code for all bitstrings of length~$n$ that avoid~$\tau$.
If $\tau$ does not satisfy these conditions, then for infinitely many $n\geq 1$ there are balancedness problems that prevent the existence of a Gray code.
Liu, Hsu and Chung~\cite{MR1301237} constructed Gray codes for \emph{Fibonacci words}, i.e., bitstrings that avoid the factor $\tau=1^\ell$ for some fixed integer~$\ell\geq 2$.
Subsequently, Vajnovszki~\cite{MR1934834} presented a simple genlex Gray code for Fibonacci words; see Figure~\ref{fig:bits2}~(b) and also~\cite{mikawa_semba_2005}.
Chinburg, Savage and Wilf~\cite{MR1487609} consider bitstrings in which any two 1s are at least~$d\geq 2$ positions apart.
For $d=2$ these are precisely Fibonacci words (with $\ell=2$), and the paper describes the same genlex Gray code as Vajnovszki for this case.
For $d\geq 3$, they showed that balancedness problems rule out the existence of a (1-)Gray code, and they construct a 2-Gray code that uses single bitflips or transpositions of a~0 and~1; see also \cite{baril_moreira-dos-santos_2009}.
Baril, Kirgizov and Vajnovszki~\cite{MR4463209} consider bitstrings in which every maximal run of~$b$ consecutive 1s that is not a suffix is followed by a run of strictly more than $b/q$ consecutive 0s, where $q\geq 1$ is some integer parameter.
They call these strings $q$-decreasing words, and show that they are in bijection to Fibonacci words with $\ell=q+1$.
Furthermore, they present a 3-Gray code for any~$q\geq 1$, shown in Figure~\ref{fig:bits2}~(e11') for $n=8$ and $q=1$, and a 1-Gray code for $q=1$.
\openproblem{16}{They also conjecture that a 1-Gray code exists for all $q\geq 1$ and $n\geq 1$.}


Baril and Vajnovszki~\cite{MR2187406} constructed a genlex Gray code for \emph{Lucas words}, i.e., bitstrings that avoid the factor $\tau=1^\ell$ for some $\ell\geq 2$, where factors that cyclically wrap around the boundary are also avoided; see Figure~\ref{fig:bits2}~(b').
It is a 1-Gray code if $\ell+1$ does not divide~$n$, and a 2-Gray code otherwise.
In the latter case, a 1-Gray code does not exist because of balancedness problems.

The subgraph of the $n$-cube induced by Fibonacci words or Lucas words for $\ell=2$ is known as \emph{Fibonacci cube}~\cite{DBLP:journals/tpds/Hsu93} or \emph{Lucas cube}~\cite{MR1812615}, respectively, and the number of those words are given by the Fibonacci numbers [OEIS~A000045] or Lucas numbers [OEIS~A000204].

Bernini et al.~\cite{MR3404701} showed that for `almost all' factors~$\tau$, the sublist of the BRGC induced by all $\tau$-avoiding words is a 3-Gray code for all $n\geq 1$; see Figure~\ref{fig:bits2}~(c).
Specifically, this holds for every $\tau$ that is not a suffix of~$\cdots(110^\ell)(110^\ell)0$ with at least two 1s for some $\ell\geq 0$.
This is a surprisingly general result, given that for every possible length~$\ell$ of the forbidden factor, all but $\ell-2$ factors $\tau\in\{0,1\}^\ell$ satisfy this condition.
These Gray codes are all cyclic, except if $\tau=0^\ell$ for $\ell\geq 2$ or $\tau=10^\ell$ for $\ell\geq 1$, in which case they are non-cyclic 1-Gray codes.
For $\tau=0^\ell$ these are complements of Fibonacci words, so we recover Vajnovszki's aforementioned result.
Moreover, for $\tau\in\{0^\ell 1\mid\ell\geq 1\}$ they are (cyclic) 2-Gray codes.
As the BRGC is genlex, all these sublist Gray codes are also genlex.

Note that any factor~$\tau$ that ends with~1 trivially satisfies the aforementioned condition of Bernini et al.
On the other hand, if $\tau$ ends with~0, then we may simply consider the complemented factor~$\ol{\tau}$, which does ends with~1.
Consequently, we get a 3-Gray code for~$\ol{\tau}$, and complementing all of its bitstrings we obtain a 3-Gray code for~$\tau$; see Figure~\ref{fig:bits2}~(c').
Combining these observations we see every factor~$\tau$ admits a 3-Gray code for all $\tau$-avoiding bitstrings.

Many of these results were later generalized substantially via flip languages discussed in Section~\ref{sec:flip-swap}.
In that section, we will consider the BRGC with each bitstring reversed, which corresponds to factors starting (and not ending) with a particular symbol.

\subsection{Necklaces and relatives}
\label{sec:neck}

\begin{figure}
\centerline{
\includegraphics{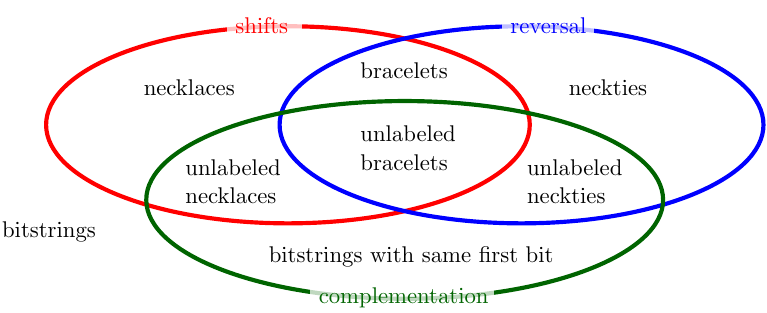}
}
\caption{Venn diagram of combinatorial objects obtained from bitstrings under shifts, reversal and/or complementation.}
\label{fig:equiv}
\end{figure}

In this section we consider equivalence classes of bitstrings under shifts, reversal and/or complementation; see Figure~\ref{fig:equiv}.
A \emph{necklace} is an equivalence class of bitstrings under cyclic shifts.
A necklace of $n$-bit strings is called \emph{aperiodic} if it has size~$n$, i.e., all $n$ cyclic shifts are distinct strings.
The lexicographically smallest representatives of aperiodic necklaces are called \emph{Lyndon words}.
A \emph{necktie} is an equivalence class of bitstrings under reversal, and a \emph{bracelet} is an equivalence class under cyclic shifts and reversal.
The goal is to find a Gray code that visits each equivalence class exactly once, where we may also impose restrictions on the choice of representatives, such as for example the lexicographically smallest one from each class.

Vajnovszki~\cite{MR2306525} showed that the sublists of the BRGC given by lexicographically largest representatives of necklaces, aperiodic necklaces or bracelets form a cyclic 3-Gray code; see Figure~\ref{fig:bits2}~(d1)--(d3).
If we are interested in the lexicographically smallest representatives, we can apply complementation to the strings in those Gray codes.
Weston and Vajnovszki~\cite{MR2380356} generalized these results for necklaces and aperiodic necklaces to larger alphabets.
It turns out we can do even better by reversing each string in the BRGC.
While this does not change the order of the (reversed) strings among each other, we will change the choice of representatives, which results in a different sublist.
Specifically, Vajnovszki~\cite{MR2399878} showed that the sublist approach applied to the BRGC with reversed strings gives a cyclic 2-Gray code for lexicographically smallest representatives of necklaces and aperiodic necklaces, in particular for Lyndon words; see Figure~\ref{fig:bits2}~(e1)+(e2).
He also obtains a cyclic 3-Gray code for subsets of these two classes, namely for strings that are lexicographically smallest under cyclic shifts and complementation, which he calls \emph{unlabeled necklaces} and \emph{unlabeled Lyndon words}; see Figure~\ref{fig:bits2}~(e1')+(e2').
As all of these Gray codes are sublists of the BRGC, they are clearly genlex.
Sawada, Williams and Wong~\cite{MR3719917} gave a constant amortized-time algorithm for Vajnovszki's Gray codes.
The aforementioned results for necklaces, Lyndon words and their unlabeled variants can be generalized substantially via flip languages and flip-swap languages discussed in the next section.

For even values of~$n$, a 1-Gray code for necklaces does not exist due balancedness problems (regardless of the choice of representatives); see~\cite{MR1413286}.
However, for odd values of~$n$ this problem is still open:
\openproblem{17}{Is there a 1-Gray code for necklaces of length~$n$ for odd~$n$?}
Note that~$0^n$ and~$1^n$ must be the first and last necklaces in such a Gray code, so the code cannot be cyclic.
The choice of representatives in this problem is arbitrary, ideally of course they are the lexicographically smallest ones.
Savage conjectured~\cite{MR1491049} that the answer to this question is `yes' (even in the more general setting of necklaces of $b$-ary strings for any~$b\geq 2$, when the flip operation is to change one entry by~$\pm 1$).
For prime values of~$n$, Degni and Drisko~\cite{MR2285811} conjectured a very simple choice of representatives and ordering, which they checked works for all primes~$n\leq 37$.
Specifically, they consider the sublist of the BRGC obtained by choosing for each necklace the first representative in this list.
Figure~\ref{fig:bits2}~(f) shows this 1-Gray code for~$n=7$.
Note that for prime~$n$, all necklaces except~$0^n$ and~$1^n$ are aperiodic.
As explained in Section~\ref{sec:track}, certain cyclic 1-Gray codes for all necklaces except~$0^n$ and~$1^n$ would yield optimal single-track Gray codes for prime~$n$.

Wang~\cite{wang_1993} proved that there is a Gray code for the lexicographically smallest representatives of neckties if and only if $n\geq 3$ is odd.
A simple inductive construction of a necktie Gray code~$N_n$ for odd~$n$ is given by~$N_1:=0,1$ and $N_n:=0N_{n-2}0,0\rev(M_{n-2})1,1N_{n-2}1$, where $M_n$ is any Gray code for all bitstrings of length~$n$ starting with~$0^n$ and ending with~$1^n$, for example the aforementioned monotone Gray code due to Savage and Winkler; see Figure~\ref{fig:bits2}~(g).

\openproblem{18}{No 2- or 1-Gray code is known for bracelets; see~{\cite{MR1491049}}.}
Balancedness problems rule out a 1-Gray code for even~$n$.

\openproblem{19}{Is there a Gray code for unlabeled neckties, i.e., strings under reversal and complementation, or for unlabeled bracelets, i.e., strings under cyclic shifts, reversal and complementation?}

\subsection{Flip-swap languages}
\label{sec:flip-swap}

A \emph{flip-swap language} is a subset~$S$ of bitstrings of length~$n$ such that~$S$ is closed under the following two operations (when applicable): (i)~flipping the leftmost~1; (ii)~transposing the leftmost~1 with the bit to its right.
If operation~(ii) is omitted, we call $S$ a \emph{flip language}.

Examples of flip-swap languages include all bitstrings, necklaces (lex.\ smallest representatives) and Lyndon words, strings with weight~$\leq \ell$, strings that avoid the factor~$10^\ell$ for some fixed integer~$\ell\geq 1$, strings with $\leq \ell$ inversions w.r.t.~$0^*1^*$, strings with $\leq \ell$ transpositions w.r.t.~$0^*1^*$, strings that are lexicographically~$\leq$ some fixed string, and strings that are $<$ or $\leq$ their reversal (in particular, lex.\ smallest neckties).
The number of \emph{inversions w.r.t.~$0^*1^*$} of a bitstring is the number of pairs of bits~$(x_i,x_j)$ where~$x_i$ is left of~$x_j$ and $(x_i,x_j)=(1,0)$, and the number of \emph{transpositions w.r.t.~$0^*1^*$} is the minimum number of exchanges of a~1 and~0 so that all 0s are left of all~1s.

Clearly, every flip-swap language is also a flip language.
Examples of flip languages that are not already mentioned before are unlabeled necklaces (lex.\ smallest representatives) and unlabeled Lyndon words, strings that avoid any factor~$\tau$ starting with~1, and $q$-decreasing words discussed in Section~\ref{sec:factor-avoid} for any $q\geq 1$.

Vajnovszki~\cite{MR2306525} proved that for the BRGC with reversed bitstrings, the sublist induced by the strings of any flip language forms a cyclic 3-Gray code; see Figure~\ref{fig:bits2}~(e1'), (e2') and~(e9')--(e11').
Sawada, Williams and Wong~\cite{MR4488724} strengthened this by showing that the sublist induced by the strings of any flip-swap language forms a cyclic 2-Gray code; see Figure~\ref{fig:bits2}~(e1)--(e8).

Note that complementing a flip-swap language 2-Gray code yields a 2-Gray code for the complemented flip-swap language, for example strings with weight~$\geq \ell$, necklaces and aperiodic necklaces with lexicographically largest representatives, strings that avoid~$01^\ell$, strings that are $\geq$ some fixed string, and strings that are $>$ or $\geq$ their reversal (in particular, lex.\ largest neckties).
Similarly, complementing a flip language 3-Gray code yields a 3-Gray code for the complemented flip language, for example unlabeled necklaces and aperiodic necklaces with lexicographically largest representatives, and strings that avoid any factor~$\tau$ starting with~0.

Observe furthermore that for any two flip(-swap) languages~$S$ and~$S'$, the intersection~$S\cap S'$ and the union~$S\cup S'$ are also flip(-swap) languages.
For example, Lyndon words that are $\leq$ their reversal are a flip-swap language.
Similarly, words avoiding any number of distinct factors all starting with~1 are a flip language.
As an application, recall that Chinburg, Savage and Wilf~\cite{MR1487609} considered bitstrings in which any two 1s are at least~$d\geq 2$ positions apart.
Those are characterized by avoiding the factors~$\{11,101,1001,\ldots,10^{d-2}1\}$, which is a flip language, so we immediately obtain a 3-Gray code for this problem, shown in Figure~\ref{fig:bits2}~(e9') for $d=2$ (as mentioned before, this can be improved by a 2-Gray code tailored to this problem).

Taking intersections of flip languages enables us to consider avoidance of factors that contain wildcard symbols~$*$, which can match both~0 or~1, also known \emph{partial words}.
For example, the factor $10*$ stands for `10 followed by another arbitrary symbol', and it can be modelled by avoiding both~$\{100,101\}$; see Figure~\ref{fig:bits2}~(e10').

\subsection{Larger alphabets}
\label{sec:larger}

The definition of BRGC can be generalized straightforwardly to larger alphabets, which yields a listings of all $b$-ary strings of length~$n$ where only one of the entries changes by~$\pm 1$ at a time: $\Gamma_0:=\varepsilon$ and $\Gamma_{n+1}:=0\Gamma_n,1\rev(\Gamma_n),2\Gamma_n,3\rev(\Gamma_n),\ldots,(b-1)\Gamma'$ for $n\geq 0$, where $\Gamma':=\rev(\Gamma_n)$ if $b$ is even and $\Gamma':=\Gamma_n$ otherwise; see \cite{DBLP:journals/tc/Flores56,MR0156723,MR538673,5009360,535751,MR1349941}.
Recently, Herter and Rote~\cite{MR3906879} devised loopless algorithms for generating $b$-ary Gray codes based on generalizations of the Tower of Hanoi puzzle.
Flahive~\cite{MR2448878} describes a construction of balanced $b$-ary Gray codes with transitions counts that differ by at most~2 from the average~$b^n/n$, improving upon the earlier paper by Flahive and Bose~\cite{MR2302538}.
\openproblem{20}{Many other of the properties discussed before for the binary case (weight-monotonicity, runs lengths, trend-freeness, single-track etc.) are still unexplored for the general $b$-ary case.}


A \emph{signed binary representation} (SBR) of an integer~$N$ is a string $(a_{n-1},\ldots,a_1,a_0)$ over the alphabet $\{-1,0,1\}$ such that $N=\sum_{i=0}^{n-1} a_i2^i$.
An SBR is said to be \emph{minimal} if the number of nonzero digits is minimum.
Some values of~$N$ have only a single minimal SBR, whereas other values admit an exponential (in $n$) number of minimal SBRs.
Sawada~\cite{MR2299690} developed a genlex 3-Gray code for listing all minimal SBRs of an integer~$N$.

\section{Combinations}
\label{sec:comb}

Combinations are the subsets of the ground set~$[n]$ of fixed size~$k$, where $0\leq k\leq n$.
We sometimes refer to them as \emph{$(n,k)$-combinations}, and their number is~$\binom{n}{k}$.
In a computer an $(n,k)$-combination can be conveniently represented by its characteristic vector~$x=(x_1,\ldots,x_n)$, i.e., $x_i=1$ if $i$ is in the set and $x_i=0$ otherwise.
Clearly, $x$ contains precisely $k$ many 1s, i.e., $x$ has weight~$k$.

In this section, we discuss Gray codes for combinations that use different flip operations, as well as Gray codes for restricted classes of combinations, including necklaces and Catalan objects (the latter can be seen as special $(2k,k)$-combinations).
Lastly, we consider more general Gray codes that list bitstrings for certain admissible weight values and weight ranges, and also some generalizations and variants of combinations.

\begin{figure}
\centerline{
\setlength{\tabcolsep}{2pt}
\renewcommand{\arraystretch}{1.5}
\begin{tabular}{ccccccccccccccccc}
(a) & (b) & (c) & (d) & (d') & (e) & (e') & (f) & (g1) & (g2) & (g3) & (g4) & (g5) & (g6) & (g7) & (g8) & (h) \\[-5mm]
\raisebox{-\height}{\includegraphics{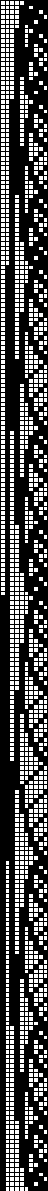}} &
\raisebox{-\height}{\includegraphics{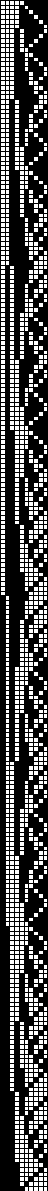}} &
\raisebox{-\height}{\includegraphics{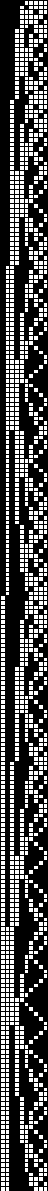}} &
\raisebox{-\height}{\includegraphics{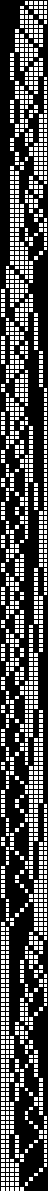}} &
\raisebox{-\height}{\includegraphics{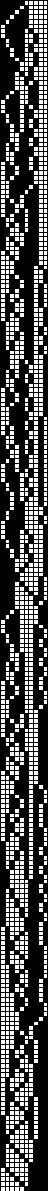}} &
\raisebox{-\height}{\includegraphics{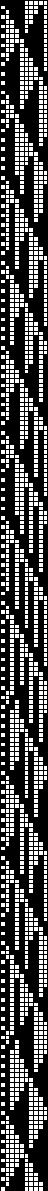}} &
\raisebox{-\height}{\includegraphics{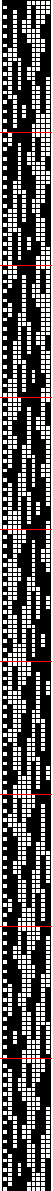}} &
\raisebox{-\height}{\includegraphics{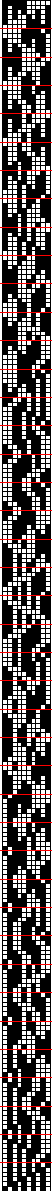}} &
\raisebox{-\height}{\includegraphics{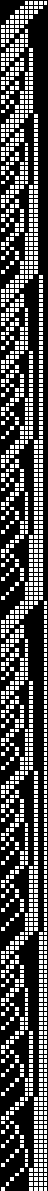}} &
\raisebox{-\height}{\includegraphics{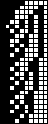}} &
\raisebox{-\height}{\includegraphics{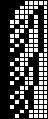}} &
\raisebox{-\height}{\includegraphics{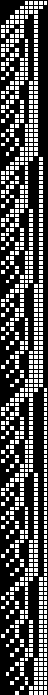}} &
\raisebox{-\height}{\includegraphics{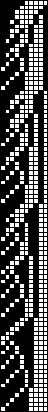}} &
\raisebox{-\height}{\includegraphics{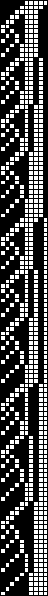}} &
\raisebox{-\height}{\includegraphics{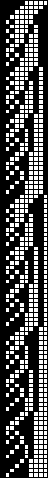}} &
\raisebox{-\height}{\includegraphics{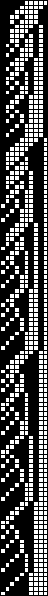}} &
\raisebox{-\height}{\includegraphics{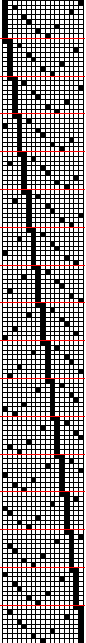}}
\end{tabular}
}
\caption{Combination Gray codes (1-bits black, 0-bits white):
(a)~revolving door from BRGC;
(b)~homogeneous transpositions $00\cdots 01\leftrightarrow 10\cdots 00$;
(c)~transpositions $01\leftrightarrow 10$ or $001\leftrightarrow 100$;
(d)~adjacent transpositions (Eades-Hickey-Read);
(d')~adjacent transpositions (Ruskey);
(e)~star transpositions;
(e')~star transpositions + cyclic blocks;
(f)~alternating flaws;
(g)~prefix/substring shifts (cool-lex):
(g1) all $(10,5)$-combinations,
(g2)~necklaces (lex.\ largest),
(g3)~aperiodic necklaces (lex.\ largest),
(g4)~0111-avoiding,
(g5)~$\leq 10$ inversions w.r.t.~$1^*0^*$,
(g6)~$\leq 2$ transpositions w.r.t.~$1^*0^*$,
(g7)~$\geq 1001100010$,
(g8)~$\geq$ their reversal;
(h)~balanced $(17,2)$-combinations by transpositions with cyclic blocks.
}
\label{fig:comb1}
\end{figure}

\subsection{Transpositions}
\label{sec:comb-trans}

We first focus on transposition Gray codes or \emph{revolving door} Gray codes for combinations, i.e., any two consecutive combinations differ by swapping a~0 and~1 in the bitstring representation.
Tang and Liu~\cite{MR0349274} showed that restricting the binary reflected Gray code (BRGC) for bitstrings of length~$n$ to any fixed weight~$k$ yields a cyclic transposition Gray code for $(n,k)$-combinations; see Figure~\ref{fig:comb1}~(a).
The resulting listing is genlex, and its ordered set representation has the interesting property that for any sublist with the same fixed prefix the next larger element is always non-decreasing or non-increasing, alternatingly depending on the parity of the prefix length, which leads to an efficient generation algorithm; see~\cite{DBLP:journals/toms/PayneI79} and~\cite[Sec.~7.2.1.3]{MR3444818}.
Bitner, Ehrlich and Reingold~\cite{MR0424386} devised a loopless algorithm to generate this ordering.
The flip graphs of $(n,k)$-combinations under arbitrary transpositions are known as \emph{Johnson graphs~$J(n,k)$}, which were shown to be Hamilton-connected by Jiang and Ruskey~\cite{MR1298971} and by Knor~\cite{MR1301949}; see Figure~\ref{fig:flip}~(b) and also~\cite{MR2928441}.

Chase~\cite{DBLP:journals/cacm/Chase70} and Eades and McKay~\cite{MR782221} constructed so-called \emph{homogeneous} Gray codes for combinations, which use only transpositions of the form $00\cdots 01\leftrightarrow 10\cdots 00$, i.e., the bits between the swapped~0 and~1 are all~0s; see Figure~\ref{fig:comb1}~(b).
We can think of this as an ordering that plays all possible combinations of~$k$ keys out of~$n$ available keys on a piano, without ever crossing any fingers.
Their homogeneous Gray code also has the genlex property.
Chase~\cite{MR995888} showed that we can restrict the allowed swaps further and only allow transpositions of the form $01\leftrightarrow 10$ or $001\leftrightarrow 100$, and this listing is also genlex and can be implemented efficiently; see Figure~\ref{fig:comb1}~(c) and also~\cite{MR1288775,MR1352777}.

Buck and Wiedemann~\cite{MR737262} and Eades, Hickey and Read~\cite{MR821383} independently proved that all $(n,k)$-combinations can be generated by using only adjacent transpositions $01\leftrightarrow 10$ if and only if $k\in\{0,1\}$ or $n-k\in\{0,1\}$ or $k$ and $n-k$ are both odd; see Figure~\ref{fig:comb1}~(d).
Note that such a Gray code must start and end at the combinations~$1^k0^{n-k}$ and~$0^{n-k}1^k$, as these admit only a single transposition (they are degree-1 vertices in the corresponding flip graph).
An efficient algorithm to compute the Eades-Hickey-Read ordering was provided by Hough and Ruskey~\cite{MR978522}, and an efficient algorithm for computing a distinct adjacent transposition listing was given by Ruskey~\cite{MR936104}; see Figure~\ref{fig:comb1}~(d').

Extending earlier partial results of Joichi and White~\cite{MR578058}, Enns~\cite{MR1246674} showed that $(n,k)$-combinations can be generated by circularly adjacent transpositions, i.e., including transpositions that swap the first and last bit, if and only if at least one of $k$ or $n-k$ is odd.

Felsner, Kleist, M\"utze and Sering~\cite{MR4046775} showed that for odd~$n$, all $(n,2)$-combinations can be generated by transpositions so that each of the $\binom{n}{2}$ possible transpositions is used exactly once.
\openproblem{21}{They raised the question where such balanced Gray codes are possible for $(n,k)$-combinations in general, whenever $n$ is odd and $\binom{n}{2}$ divides $\binom{n}{k}$ (it is easy to argue that for even~$n$ such Gray codes cannot exist).}
Their balanced Gray code for $(n,2)$-combinations can be partitioned into~$n$ blocks, where each block is a cyclic right-shift of the previous one; see Figure~\ref{fig:comb1}~(h).
Consequently, it is also a single-track Gray code in which column~$i$ is a cyclic shift of the first column by $\frac{1}{n}\binom{n}{2}\cdot(i-1)$ for all $i=2,\ldots,n$.

\subsection{Cool-lex order and bubble languages}
\label{sec:bubble}

Another elegant and efficient method for generating combinations based on prefix shifts was described by Ruskey and Williams~\cite{MR2548545}, who named their ordering \emph{cool-lex}, because of its similarities to colex and its numerous striking features.
A \emph{prefix shift} cyclically left-shifts a prefix of the bitstring by one position to generate the next combination.
The following simple successor rule generates the next combination in the cool-lex ordering: Shift the shortest prefix of length at least~3 that ends with~01, or if no such prefix exists shift the entire string; see Figure~\ref{fig:comb1}~(g1) and cf.~\cite{MR3193758}.
The cool-lex ordering is cyclic, any two consecutive combinations differ by one or two transpositions, and it admits simple loopless algorithms.
Moreover, the ordering is genlex in the ordered set representation (but not in the bitstring representation).

Ruskey, Sawada and Williams~\cite{MR2844089} present another major application of the cool-lex approach.
They consider subsets~$B$ of $(n,k)$-combinations that are \emph{bubble languages}, which means that if a string~$x$ is in~$B$ and~$x$ contains~01, then the string obtained from~$x$ by flipping the first such occurrence to~10 must also be in~$B$.
Examples of bubble languages include combinations, $b$-ary Dyck words, necklaces and aperiodic necklaces with lexicographically largest representatives, strings that avoid the factor~$01^\ell$ for some fixed integer~$\ell\geq 1$, strings with $\leq \ell$ inversions w.r.t.~$1^*0^*$, strings with $\leq \ell$ transpositions w.r.t.~$1^*0^*$, strings that are lexicographically~$\geq$ some fixed string, and strings that are $>$ or $\geq$ their reversal (all of these families are considered for a fixed number~$k$ of 1s).
The number of \emph{inversions w.r.t.~$1^*0^*$} of a combination $x=(x_1,\ldots,x_n)$ is the number of pairs of bits $(x_i,x_j)$ with~$i<j$ and $(x_i,x_j)=(0,1)$, and the number of \emph{transpositions w.r.t.~$1^*0^*$} is the minimum number of exchanges of a~0 and~1 so that all 1s are left of all 0s.
The main result of~\cite{MR2844089} is that the words of any bubble language, when arranged in the order induced by the cool-lex order, form a cyclic Gray code in which any two consecutive strings differ by a left-shift of a substring by one position, which results in one or two transpositions; see Figure~\ref{fig:comb1}~(g2)--(g8) and Figure~\ref{fig:comb2}~(f2)+(f3).
Note that complementing such a listing yields a shift Gray code for the complemented bubble language, for example necklaces and aperiodic necklaces with lexicographically smallest representatives (i.e., fixed-density Lyndon words), strings that avoid~$10^\ell$, strings that are~$\leq$ some fixed string, and strings that are $<$ or $\leq$ their reversal.
Observe furthermore that for any two bubble languages~$B$ and~$B'$, the intersection~$B\cap B'$ and the union~$B\cup B'$ are also bubble languages.
For example, $b$-ary Dyck words that avoid~$01^\ell$ are a bubble language.
Following up on these results, Sawada and Williams~\cite{MR2900417,MR3101690} devised constant amortized-time algorithms for computing many of the aforementioned bubble language Gray codes.

\openproblem{22}{Is there a unifying principle behind flip-swap languages and bubble languages (recall Section~{\ref{sec:flip-swap}})?}

\openproblem{23}{Another interesting open question in this context is whether all $(n,k)$-combinations can be generated by \emph{prefix reversals}, i.e., in each step, a prefix of the bitstring representation is reversed to obtain the next combination.}
Such orderings can be constructed easily for the cases $k\in\{1,2\}$ or $n-k\in\{1,2\}$, but no general construction is known.

\subsection{Fixed-density necklaces}

We now consider equivalence classes of $(n,k)$-combinations under cyclic shifts, so-called \emph{necklaces with fixed density}.
Wang and Savage~\cite{MR1413286} devised a cyclic transposition Gray code, valid for all~$n$ and~$k$.
The representatives from each equivalence class are not the lexicographically smallest ones, but those rotated left until the first bit equals~1.
Independently, Ueda~\cite{MR1761724} described another transposition Gray code for necklaces with fixed density, and also a Gray code for aperiodic necklaces with fixed density.
However, the chosen representatives in his two Gray codes are non-standard, in particular, they are neither the lexicographically smallest nor largest ones.
He also established stronger Hamiltonicity properties about the underlying flip graphs.

Shift Gray codes for fixed density necklaces and aperiodic necklaces, either for lexicographically smallest or largest representatives, i.e., in particular for fixed density Lyndon words, can be derived from the bubble-language framework discussed in the previous section; recall Figure~\ref{fig:comb1}~(g2)+(g3).
These Gray codes operate by substring shifts that result in one or two transpositions per step.

\subsection{Catalan objects}

\emph{Dyck words} are $(2k,k)$-combinations with the additional property that every prefix contains at most as many~0s as~1s.
By interpreting the 1s as opening brackets and the 0s as closing brackets, we see that Dyck words are in bijection with well-formed parenthesis expressions with $k$ pairs of parenthesis.
They are also in bijection to preorder sequences of binary trees with $k$ nodes, when recording a 1 for every internal node and a 0 for every leaf (the 0 corresponding to the last leaf being implicit).
It is well known that they are counted by the $k$th Catalan number $C_k=\frac{1}{k+1}\binom{2k}{k}$ [OEIS~A000108].

\begin{figure}
\centerline{
\setlength{\tabcolsep}{3pt}
\renewcommand{\arraystretch}{1.5}
\begin{tabular}{ccccccccccccc}
(a) & (b) & (c) & (d2) & (d3) & (e2) & (e3) & (f2) & (f3) & (g5) & (g6) & (g7) & (h) \\[-5mm]
\raisebox{-\height}{\includegraphics{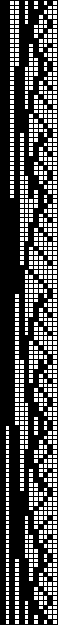}} &
\raisebox{-\height}{\includegraphics{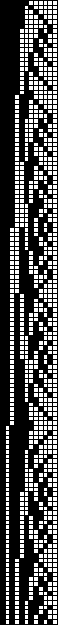}} &
\raisebox{-\height}{\includegraphics{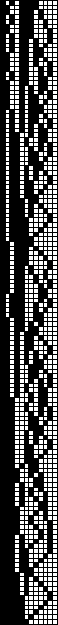}} &
\raisebox{-\height}{\includegraphics{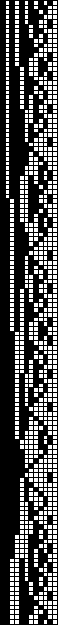}} &
\raisebox{-\height}{\includegraphics{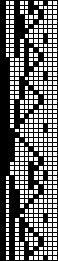}} &
\raisebox{-\height}{\includegraphics{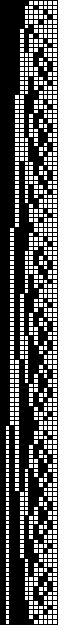}} &
\raisebox{-\height}{\includegraphics{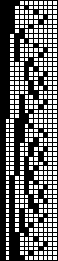}} &
\raisebox{-\height}{\includegraphics{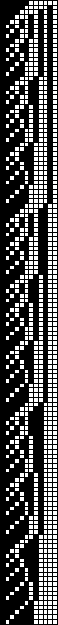}} &
\raisebox{-\height}{\includegraphics{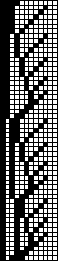}} &
\raisebox{-\height}{\includegraphics{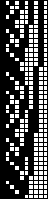}} &
\raisebox{-\height}{\includegraphics{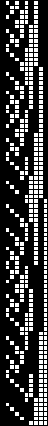}} &
\raisebox{-\height}{\includegraphics{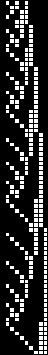}} &
\raisebox{-\height}{\includegraphics{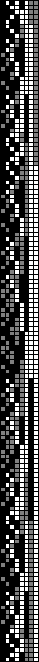}}
\end{tabular}
}
\caption{Gray codes for Dyck words and relatives (1-bits black, 0-bits white; $-$ in Motzkin words gray):
(a)~transpositions;
(b)~homogeneous transpositions;
(c)~adjacent transpositions;
(d2)+(d3)~2/3-ary, homogeneous transpositions $00\cdots 01\leftrightarrow 10\cdots 00$;
(e2)+(e3)~2/3-ary, transpositions $01\leftrightarrow 10$ or $001\leftrightarrow 100$;
(f2)+(f3)~2/3-ary, prefix shifts (cool-lex);
(g5)--(g7)~Dyck prefixes of length~10 with weight 5 (=Dyck words), 6 or 7;
(h)~Motzkin words with fixed content $\{1,1,1,0,0,0,-,-\}$.
}
\label{fig:comb2}
\end{figure}

Proskurowski and Ruskey~\cite{MR789905} devised a transposition Gray code for Dyck words, and in a subsequent paper~\cite{MR1041167}, they present a simpler genlex construction for this problem; see Figure~\ref{fig:comb2}~(a) and also~\cite{MR1645246}.
They also showed that an adjacent transposition Gray code exists if and only if $k$ is even or $k<5$; see Figure~\ref{fig:comb2}~(c).
A homogeneous transposition Gray code for Dyck words, which is also genlex, was given by Bultena and Ruskey~\cite{MR1664740}; see Figure~\ref{fig:comb2}~(b).
More generally, Roelants van Baronaigien~\cite{MR1747718} developed a homogeneous transposition Gray code \emph{$b$-ary Dyck words}, valid for any integer $b\geq 2$; see Figure~\ref{fig:comb2}~(d2)+(d3); see also~\cite{MR1807675}.
These are are strings with $k$ many 1s and $k(b-1)$ many 0s such that in every prefix, the number of 0s is at most $b-1$ times the number of 1s, which are in bijection to preorder sequences of $b$-ary trees.
Vajnovszki and Walsh~\cite{MR2577686} present a genlex Gray code for this problem that is even more restrictive and uses only transpositions~$01\leftrightarrow 10$ or $001\leftrightarrow 100$; see Figure~\ref{fig:comb2}~(e2)+(e3).

Ruskey and Williams~\cite{DBLP:conf/cats/RuskeyW08} present a Gray code for Dyck words that uses prefix shifts, which cyclically left-shifts a substring that starts at the second position (the first bit clearly has to always be~1); see Figure~\ref{fig:comb2}~(f2).
Any two consecutive words in this ordering differ by one or two transpositions, and the ordering has a very simple successor rule and loopless algorithm.
Durocher, Li, Mondal, Ruskey and Williams~\cite{MR2960361} generalized this construction to higher arity~$b\geq 2$; see Figure~\ref{fig:comb2}~(f3).

Dyck words can be generalized to \emph{words with $e$ flaws}, where $0\leq e\leq k$ is an integer parameter, which are $(2k,k)$-combinations with the property that precisely $e$ prefixes ending with~0 have strictly more~0s than~1s.
We denote them by $F_{k,e}$.
Clearly, $F_{k,0}$ are Dyck words, and $F_{k,k}$ are complemented Dyck words.
The well-known Chung-Feller theorem asserts that $|F_{k,0}|=|F_{k,1}|=\cdots=|F_{k,k}|=\frac{1}{k+1}\binom{2k}{k}=C_k$, i.e., the number of words with $e$ flaws is the same for each possible value of~$e$.
M\"utze, Standke and Wiechert~\cite{MR3738156} constructed a bijection between $F_{k,e}$ and $F_{k,e+1}$, for any $e=0,\ldots,k-1$, that uses only adjacent transpositions, and that maps any Dyck word in~$F_{k,0}$ to its complement in~$F_{k,k}$ when applied $k$ times.
By combining this Gray code with any of the aforementioned transposition Gray codes for Dyck words, we obtain a transposition Gray code for $(2k,k)$-combinations in which the number of flaws is alternatingly increasing and decreasing; see Figure~\ref{fig:comb1}~(f).

Alternatively, Dyck words can be generalized to \emph{Dyck prefixes with fixed weight}, which are $(n,k)$-combinations for some $k\geq n/2$ such that every prefix contains at most as many~0s as~1s.
A transposition Gray code for those objects, which is furthermore genlex, was derived by Vajnovszki and Wong~\cite{DBLP:conf/cocoon/VajnovszkiW23}, and it has a simple greedy description; see Figure~\ref{fig:comb2}~(g5)--(g7).
The authors also interleave these listings into a Gray code for \emph{ballot sequences}, which is the set of all Dyck prefixes of length~$n$ with weights $n/2\leq k\leq n$, and this Gray code only uses transpositions and flips of single bits.
In fact, the transpositions used in~\cite{DBLP:conf/cocoon/VajnovszkiW23} are all of the form~$11\cdots 110\leftrightarrow 01\cdots 11$, i.e., they are homogeneous in the complemented bitstrings.

A \emph{Motzkin word with fixed content} is a word over the alphabet~$\{0,1,-\}$ with $k$ many 1s, $k$ many 0s, and a fixed number~$m$ of occurrences of~$-$, with the property that every prefix contains at most as many~0s as~1s.
Consequently, Motzkin words for $m=0$ are Dyck words.
Vajnovszki~\cite{MR1943198} presents a genlex Gray code for Motzkin words with fixed content that either transposes two symbols or permutes three symbols in each step; see Figure~\ref{fig:comb2}~(h).
These can be combined into a Gray code for Motzkin words of fixed length with an arbitrary number of~$-$s.
Lapey and Williams~\cite{DBLP:conf/iwoca/LapeyW22} devised a shift Gray code for {\L}ukasiewicz words with fixed content, which generalize both Dyck words and Motzkin words with fixed content.

Further Gray codes for Catalan objects will be discussed in Sections~\ref{sec:avoid}, \ref{sec:posets}, \ref{sec:setp-geom}, \ref{sec:triang} and~\ref{sec:otrees}.

\subsection{Weight ranges}
\label{sec:range}

As discussed before, combinations can be represented by bitstrings of length~$n$ with weight~$k$.
On the other hand, in Section~\ref{sec:bits} we discussed Gray codes for all bitstrings of length~$n$, irrespective of their weight.
Generalizing both types of objects, we may consider bitstrings of length~$n$ whose weight is in some interval~$[k,\ell]$, where $0\leq k\leq \ell\leq n$.
The corresponding flip graph, denoted by~$Q_{n,[k,\ell]}$, is the subgraph of the $n$-cube induced by all levels from~$k$ to~$\ell$, where the $k$th \emph{level of~$Q_n$} contains all bitstrings with weight~$k$.
The graph~$Q_{n,[k,\ell]}$ is clearly bipartite, and it is easy to show that the partition classes have the same size if and only if $n$ is even and $(k,\ell)=(0,n)$, or $n$ is odd and $k+\ell=n$, i.e., the boundary levels~$k$ and~$\ell$ are symmetric around the middle.
Consequently, a Hamilton cycle in~$Q_{n,[k,\ell]}$ can only exist in one of those cases.

One of most intriguing instances of this problem is the \emph{middle levels conjecture}, which asks for a Hamilton cycle in~$Q_{2k+1,[k,k+1]}$ for all~$k\geq 1$, i.e., a Hamilton cycle in the middle two levels of the $(2k+1)$-cube.
This question was first raised by Havel~\cite{MR737021} and Buck and Wiedemann~\cite{MR737262} independently, and it received considerable attention in the literature.
The conjecture was answered affirmatively by M\"utze~\cite{MR3483129}, who showed that there are in fact double-exponentially (in $k$) many distinct Hamilton cycles.
Subsequently, a shorter and more accessible proof of the middle levels conjecture was presented by Gregor, M\"utze and Nummenpalo~\cite{MR3819051}, which reduces the problem to finding a spanning tree in an auxiliary graph on plane trees.
This construction was turned into a loopless algorithm by M\"utze and Nummenpalo~\cite{MR4075363}.
An even shorter 2-page `book proof' of the middle levels conjecture was presented by M\"utze~\cite{MR4707569}.
Consider a Hamilton cycle in the middle levels graph~$Q_{2k+1,[k,k+1]}$, and prefix every bitstring of weight~$k$ with~1 and every bitstring of weight~$k+1$ with~0.
Note that the resulting listing is a star transposition ordering of $(2k+2,k+1)$-combinations, i.e., in each step, the first bit is swapped with some later bit; see Figure~\ref{fig:comb1}~(e).
Knuth~\cite{MR3444818} raised a strengthened version of the middle levels conjecture, asking for a Hamilton cycle in~$Q_{2k+1,[k,k+1]}$ that can be partitioned into $2k+1$ blocks such that each block is a cyclic shift of the previous one (without the first column).
Such an ordering has the property that any two of the columns~$2,\ldots,2k+2$ of bits differ by a cyclic shift, i.e., it is a single-track Gray code, with the exception of the first column which is simply~$101010\cdots$.
Furthermore, such a Hamilton cycle corresponds to a Gray code for necklaces of length~$2k+1$ with weights~$k$ or~$k+1$.
A general algorithmic solution for this problem was presented by Merino, Mi\v{c}ka and M\"utze~\cite{MR4417008}; see Figure~\ref{fig:comb1}~(e').
Moreover, the middle levels graph~$Q_{2k+1,[k,k+1]}$ was shown to be Hamilton-laceable by Gregor, Merino and M\"utze~\cite{MR4580020}.
While the algorithms described in~\cite{MR4075363} and~\cite{MR4417008} solve the middle levels problem and its strengthened version efficiently, these algorithms are admittedly rather complicated.
\openproblem{24}{It remains open to find a simple algorithm for solving either of these problems, ideally in the form of a successor rule.}
The authors of~\cite{MR4417008} offer a reward of $2^5$~EUR for an algorithm that can be implemented with at most 4000~ASCII characters of C++ code.

\begin{figure}
\centerline{
\setlength{\tabcolsep}{3pt}
\renewcommand{\arraystretch}{1.5}
\begin{tabular}{cccccccccc}
(a0) & (a1) & (a2) & (a3) & (b0) & (b1) & (b2) & (b3) \\[-5mm]
\raisebox{-\height}{\includegraphics{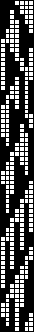}} &
\raisebox{-\height}{\includegraphics{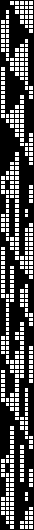}} &
\raisebox{-\height}{\includegraphics{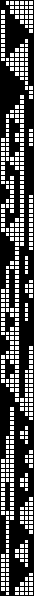}} &
\raisebox{-\height}{\includegraphics{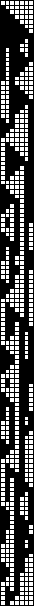}} &
\raisebox{-\height}{\includegraphics{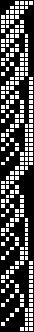}} &
\raisebox{-\height}{\includegraphics{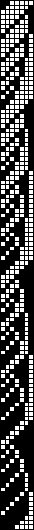}} &
\raisebox{-\height}{\includegraphics{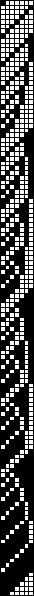}} &
\raisebox{-\height}{\includegraphics{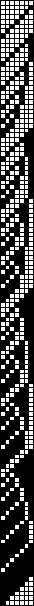}}
\end{tabular}
}
\caption{Weight-range combination Gray codes (1-bits black, 0-bits white):
(a0)--(a3)/(b0)--(b3)~central levels graph $Q_{7,[3-c,4+c]}$ for $c=0,1,2,3$ by single bitflips or prefix shifts, respectively.
}
\label{fig:comb3}
\end{figure}

The more general \emph{central levels problem} asks for a Hamilton cycle in~$Q_{2k+1,[k-c,k+1+c]}$ for all~$k\geq 1$ and $0\leq c\leq k$.
This problem was raised independently by Buck and Wiedemann~\cite{MR737262}, Savage~\cite{MR1275228}, Gregor and \v{S}krekovski~\cite{MR2609124}, and by Shen and Williams~\cite{shen_williams_2019}.
A solution for the case~$c=k$ (all levels) is given by the BRGC.
The case $c=k-1$ (all levels except the two extreme vertices) was solved by Buck and Wiedemann~\cite{MR737262} and in a more general setting by Locke and Stong~\cite{locke_stong_2003}; see also~\cite{MR1887372}.
Gregor and \v{S}krekovski~\cite{MR2609124} solved the case $c=k-2$ (all levels except the top two and bottom two).
The case $c=0$ (middle two levels) follows from the aforementioned solution of the middle levels conjecture.
Later, Gregor, J\"ager, M\"utze, Sawada and Wille~\cite{MR4342529} presented a solution for the case $c=1$ (middle four levels), and their construction is based on symmetric chain decompositions.
Extending these ideas further, a full solution of the central levels problem, valid for all $k\geq 1$ and $0\leq c\leq k$, was presented by Gregor, Mi\v{c}ka and M\"utze~\cite{MR4539109}; see Figure~\ref{fig:comb3}~(a0)--(a3).
\openproblem{25}{It is open to find an efficient algorithm for computing a Hamilton cycle in~$Q_{2k+1,[k-c,k+1+c]}$ for general values of~$c$.}
For $c=0$ and $c=k$ such algorithms are known by the aforementioned results.

Coming back to the most general version of the problem about~$Q_{n,[k,\ell]}$, Gregor and M\"utze~\cite{MR3758308} showed that for all~$n\geq 3$ and~$0\leq k\leq \ell\leq n$, in those cases where the partition classes of~$Q_{n,[k,\ell]}$ have different sizes, the graph admits a saturating cycle and a tight enumeration, where a \emph{saturating cycle} is a cycle in~$Q_{n,[k,\ell]}$ that visits all vertices in the smaller partition class, and a \emph{tight enumeration} is a cyclic listing of the bitstrings in both partition classes, flipping a single bit in most steps, and exchanging a~0 and~1 in $d$~steps, where $d$ is the difference in size between the partition classes (the steps flipping two bits are within the larger partition class).
Saturating cycles and tight enumerations are natural generalizations of Hamilton cycles for bipartite graphs whose partition classes are not necessarily of the same size.
Several of the weight range Gray codes devised in~\cite{MR3758308} are obtained by considering sublists of the BRGC within an interval~$[k,\ell]$.

A \emph{Boolean layer cake}~\cite{MR1688615} is obtained from~$Q_n$ by selecting any subset of complete levels (not necessarily consecutive), connecting any two vertices in subsequent levels by an edge if they are connected by a weight-monotone path in~$Q_n$ (equivalently, if the vertices are from subsequent levels~$k$ and~$\ell$, the path has length~$|k-\ell|$).
Obviously, the graph $Q_{n,[k,\ell]}$ is the Boolean layer cake obtained by selecting all levels $k,k+1,\ldots,\ell$.
\openproblem{26}{It is open which Boolean layer cakes admit Hamilton cycles, or more generally, saturating cycles and tight enumerations.}

Stevens and Williams~\cite{MR3193758} modified the aforementioned cool-lex successor rule for combinations to be able to generate all bitstrings of length~$n$ as follows: Shift the shortest prefix of length at least~3 that ends with~01, or if no such prefix exists complement the first bit and shift the entire string; see Figure~\ref{fig:comb3}~(b3).
Moreover, they generalized the rule even further so as to generate all bitstrings with weight in the interval~$[k,\ell]$ by prefix shifts:
Shift the shortest prefix of length at least~3 that ends with~01, or if no such prefix exists, complement the first bit if this keeps the weight in the interval $[k,\ell]$, and then shift the entire string; see Figure~\ref{fig:comb3}~(b0)--(b3).
In the resulting cyclic Gray codes, any two consecutive bitstrings differ in at most 4~positions.

\subsection{Bipartite Kneser graphs and Kneser graphs}

\begin{figure}
\centerline{
\setlength{\tabcolsep}{3pt}
\renewcommand{\arraystretch}{1.5}
\begin{tabular}{ccccccccccccc}
(a1) & (a2) & (a3) & (a4) & (b1) & (b2) & (b3) & (b4) & (c1) & (c2) & (d) \\[-5mm]
\raisebox{-\height}{\includegraphics{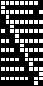}} &
\raisebox{-\height}{\includegraphics{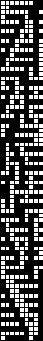}} &
\raisebox{-\height}{\includegraphics{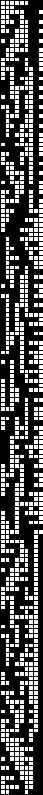}} &
\raisebox{-\height}{\includegraphics{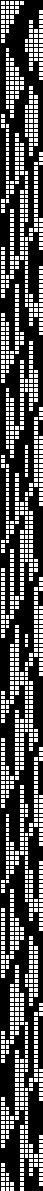}} &
\raisebox{-\height}{\includegraphics{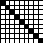}} &
\raisebox{-\height}{\includegraphics{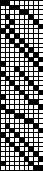}} &
\raisebox{-\height}{\includegraphics{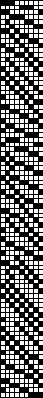}} &
\raisebox{-\height}{\includegraphics{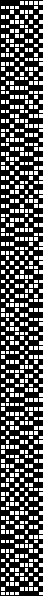}} &
\raisebox{-\height}{\includegraphics{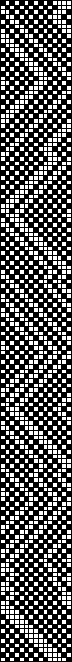}} &
\raisebox{-\height}{\includegraphics{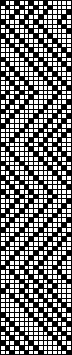}} &
\raisebox{-\height}{\includegraphics{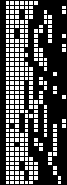}}
\end{tabular}
}
\caption{Hamilton cycles in (bipartite) Kneser graphs and relatives (1-bits black, 0-bits white):
(a1)--(a4)~bipartite Kneser graphs $H(9,k)$, $k=1,2,3,4$;
(b1)--(b4)~Kneser graphs $K(9,k)$, $k=1,2,3,4$;
(c1)~Schrijver graph $S(15,6)=S(15,6,2)$;
(c2)~$s$-stable Kneser graph $S(15,4,3)$;
(d)~numerical semigroups with genus~7.
}
\label{fig:comb4}
\end{figure}

For integers~$k\geq 1$ and~$n\geq 2k+1$, the \emph{bipartite Kneser graph~$H(n,k)$} has as vertices all $k$-element and $(n-k)$-element subsets of~$[n]$, and an edge between a set~$A$ of size~$k$ and a set~$B$ of size~$n-k$ if $A\seq B$.
Considering the characteristic vectors of those subsets, we see that $H(n,k)$ is isomorphic to the Boolean layer cake given by levels~$k$ and~$n-k$ in the $n$-cube.
In particular, $H(2k+1,k)$ is isomorphic to the middle levels graph~$Q_{2k+1,[k,k+1]}$.
Hurlbert~\cite{MR1271867} and Roth (see~\cite{MR1106528}) conjectured that $H(n,k)$ has a Hamilton cycle for all~$k\geq 1$ and~$n\geq 2k+1$.
This conjecture was proved by M\"utze and Su~\cite{MR3759914}, using an inductive construction in~$Q_n$ that uses the middle levels theorem as a base case.
Rephrased in terms of a Gray code, this result asserts that all $k$-element and $(n-k)$-element subsets of~$[n]$ (i.e., $(n,k)$- and $(n,n-k)$-combinations) can be listed cyclically, so that we alternatingly move to a subset or a superset of the current set.
The construction from~\cite{MR3759914} has the additional property that any two consecutive $k$-sets, which occur after every two steps, differ only in exchanging one element; see Figure~\ref{fig:comb4}~(a1)--(a4).

For integers~$k\geq 1$ and~$n\geq 2k+1$, the \emph{Kneser graph~$K(n,k)$} has as vertices all $k$-element subsets of~$[n]$, and an edge between any two disjoint sets.
The Kneser graphs~$K(2k+1,k)$ are also known as \emph{odd graphs~$O_k$}.
A Hamilton cycle in the Kneser graph corresponds to a cyclic listing of all $k$-element subsets of~$[n]$ (i.e., $(n,k)$-combinations), such that any two consecutive sets are disjoint.
\solvedproblem{27}{It has long been conjectured that all Kneser graphs~$K(n,k)$ have a Hamilton cycle, with one notable exception, the Petersen graph~$K(5,2)=O_2$; see Figure~{\ref{fig:flip}}~(c).}
For the odd graphs this conjecture was raised by Meredith and Lloyd~\cite{MR0457282,MR321782} and Biggs~\cite{MR556008}.
M\"utze, Nummenpalo and Walczak~\cite{MR4273468} proved that the odd graphs~$O_k=K(2k+1,k)$ have a Hamilton cycle for all $k\geq 3$ (even double-exponentially many distinct ones), using as a starting point the cycle factor in~$O_k$ from~\cite{MR3738156} whose cycles all have the same length~$k+1$; see Figure~\ref{fig:comb4}~(b4).
Combined with an inductive construction of Johnson~\cite{MR2836824}, their result implies that $K(2k+2^a,k)$ has a Hamilton cycle for all $k\geq 3$ and $a\geq 0$, in particular for $n=2k+2$.
A complete solution for all remaining cases $n\geq 2k+3$ was given by Merino, M\"utze and Namrata~\cite{MR4617441}; see Figure~\ref{fig:comb4}~(b1)--(b3).

For integers $k\geq 1$, $0\leq t\leq k-1$, and $n\geq 2k+1$ if $t=0$ or $n\geq 2k-t$ if $t>0$, the \emph{uniform subset graph~$U(n,k,t)$} has as vertices all $k$-element subsets of~$[n]$, and an edge between any two sets whose intersection has size~$t$.
By taking complements, we see that $U(n,k,t)$ is isomorphic to $U(n,n-k,n-2k+t)$.
For $t=0$ this precisely gives the Kneser graphs.
For $t=k-1$ we get the flip graphs of $(n,k)$-combinations under transpositions, the Johnson graphs~$J(n,k)$.
For this reason uniform subset graphs are also known as generalized Johnson graphs in the literature.
\solvedproblem{28}{The graphs~$U(n,k,t)$ were introduced by Chen and Lih~{\cite{MR888679}}, who conjectured that they all have a Hamilton cycle, except the Petersen graph, which arises as~$U(5,2,0)$ or~$U(5,3,1)$ (cf.~{\cite{MR1106528}}).}
A solution of this conjecture for the case $t=k-1$ follows from the transposition Gray codes for $(n,k)$-combinations presented earlier.
In addition to this, Chen and Lih settled the cases $t=k-2$ and $t=k-3$ in their paper.
A complete solution of the problem was presented by Merino, M\"utze and Namrata~\cite{MR4617441}.
These results settle Lov\'asz' conjecture mentioned in Section~\ref{sec:lovasz} for all known families of vertex-transitive graphs defined by intersecting set systems (see~\cite{MR4733111}).

For integers $k\geq 1$ and~$n\geq 2k+1$, the \emph{Schrijver graph} $S(n,k)$ is the subgraph of $K(n,k)$ induced by all subsets of~$[n]$ that contain no two circularly adjacent elements.
More generally, for integers $k\geq 1$, $s\geq 2$, and $n\geq sk+1$, the \emph{$s$-stable Kneser graph $S(n,k,s)$} is the subgraph of~$K(n,k)$ induced by all subsets of~$[n]$ in which any two elements are in circular distance at least~$s$.
Clearly, we have $S(n,k,2)=S(n,k)$.
Schrijver graphs are named after Schrijver~\cite{MR512648}, who showed that they are vertex-critical subgraphs of~$K(n,k)$, i.e. they have the same chromatic number~$n-2k+2$ as~$K(n,k)$, but removing any of their vertices yields a graph with smaller chromatic number.
In particular, $S(2k+1,k)$ is an odd cycle with chromatic number~3.
Schrijver graphs and more generally $s$-stable Kneser graphs were shown to be Hamiltonian by 
Ledezma and Pastine~\cite{ledezma_pastine_2024} and M\"{u}tze and Namrata~\cite{mutze_namrata_2024}, independently; see Figure~\ref{fig:comb4}~(c1)+(c2).

\openproblem{61}{Biggs~\cite{MR556008} conjectured that for $k\geq 3$ the edges of the odd graph~$O_k$ can be partitioned into $\lfloor (k+1)/2\rfloor$ edge-disjoint Hamilton cycles.
We do not even know two edge-disjoint Hamilton cycles in the odd graphs, the middle levels graphs~$Q_{2k+1,[k,k+1]}$, or more generally, in Kneser graphs~$K(n,k)$ or bipartite Kneser graphs~$H(n,k)$.}

\subsection{Other generalizations and variants}

Combinations are generalized by bounded integer compositions, and Gray codes for those are discussed in Section~\ref{sec:comp}.

Joichi, White and Williamson~\cite{MR557834} consider classes of combinatorial objects that satisfy the two-parameter recursion~$C_{n,k}=a_{n,k}\cdot C_{n-1,k-1}+b_{n,k}\cdot C_{n-1,k}$, which includes $(n,k)$-combinations (for $a_{n,k}=b_{n,k}=1$), set partitions of~$[n]$ into $k$~parts, permutations of~$[n]$ with $k$ cycles, and permutations of~$[n]$ with $k$ runs.
These objects are in bijective correspondence with monotone paths in the integer grid from~$(0,0)$ to~$(n-k,k)$, where the steps have edge-multiplicities as specified by the coefficients~$a_{n,k}$ and~$b_{n,k}$, i.e., they are described by an $(n,k)$-combination for the path, and an additional word of length~$n$ that describes which of the multi-edges are taken.
The paper describes a general Gray code scheme for these pairs of strings, obtained from fusing a combination Gray code with a straightforward product space Gray code.

A \emph{numerical semigroup} is a subset~$S$ of the non-negative integers~$\mathbb{N}$ such that $0\in S$, $\mathbb{N}\setminus S$ is finite, and $a,b\in S$ implies that $a+b\in S$.
The numbers in~$G:=\mathbb{N}\setminus S$ are referred to as \emph{gaps} of the semigroup, and the quantity $g:=|G|=|\mathbb{N}\setminus S|$ is referred to as its \emph{genus}.
For any semigroup with genus~$g$, we have $\max G<2g$, so any such semigroup can be represented as a $(2g,g)$-combination by considering the first $2g$ entries of the characteristic vector of~$G$.
\openproblem{29}{Is there a Gray code for numerical semigroups of genus~$g$, represented as $(2g,g)$-combinations, for every~$g\geq 1$?}
Figure~\ref{fig:comb4}~(d) shows a solution for $g=7$.
Jack Barron~\cite{barron_2023} made some progress by showing that the corresponding flip graph is 2-connected, which implies a 2-Gray code by Fleischner's~\cite{MR332573} theorem discussed in Section~\ref{sec:general}.
This problem is interesting, as for the counting sequence of numerical semigroups with genus~$g$, even basic facts such as monotonicity in~$g$ are not proved; see~\cite{MR2377597} and~[OEIS~A007323].

\section{Permutations}
\label{sec:perm}

In this section we first discuss Gray codes for listing all $n!$ permutations of~$[n]=\{1,2,\ldots,n\}$ using different flip operations.
We then discuss Gray code for interesting subsets of permutations, obtained by imposing additional constraints.
Lastly, we discuss Gray codes for generalizations and variants of permutations.

\subsection{Transpositions}

Possibly the most basic operation on a permutation is a \emph{transposition}, i.e., the exchange of two entries of the permutation.
Transposition Gray codes for permutations can be constructed inductively by appending an additional symbol in all possible ways to a Gray code with one symbol less, using transpositions with the newly appended symbol at the transitions between blocks, which yields a genlex ordering.
Wells~\cite{MR127507} devised an explicit rule for producing such a transposition ordering of permutations; see Figure~\ref{fig:perm1}~(a).
Another simple rule was found by Heap~\cite{DBLP:journals/cj/Heap63}, which is to perform the $i$th transposition involving the entry at position~$n$ in a block with the first entry if~$n$ is odd, and with the $i$th entry if $n$ is even; see Figure~\ref{fig:perm1}~(b).
Lipski~\cite{MR619774} considerably generalized the methods proposed by Wells and Heap, and came up with an entire family of transposition Gray codes for permutations, two of which are shown in Figures~\ref{fig:perm1}~(c)+(c'); see~\cite{DBLP:books/Arndt2011}.

The famous \emph{Steinhaus-Johnson-Trotter (SJT)} algorithm generates permutations by \emph{adjacent transpositions}, i.e., transpositions of neighboring entries of the permutation; see Figure~\ref{fig:flip}~(d) and Figure~\ref{fig:perm1}~(d).
This ordering can also be constructed inductively, but instead of appending an additional symbol in all possible ways, we insert the new symbol in all possible ways from right to left alternatingly.
This algorithm was discovered independently by Trotter~\cite{DBLP:journals/cacm/Trotter62} and Johnson~\cite{MR159764}, and is also implicit in a book by Steinhaus~\cite[p.~49--50]{MR0157881}, where he describes a puzzle about generating all permutations of beads on an abacus.
The SJT algorithm is also known as the method of \emph{plain changes}; see Knuth's book~\cite[Sec.~7.2.1.2]{MR3444818} for the historical origins of this algorithm in the art of change ringing of church bells.
Ehrlich gives a loopless implementation of the SJT algorithm~\cite{DBLP:journals/cacm/Ehrlich73,MR0366085}.
Sedgewick~\cite{MR464682} provides an insightful analysis and fair comparison of the computational efficiency of several permutation generation algorithms; see also Arndt's book~\cite{DBLP:books/Arndt2011}.

\begin{figure}
\centerline{
\setlength{\tabcolsep}{3pt}
\renewcommand{\arraystretch}{1.5}
\begin{tabular}{cccccccccccccc}
\raisebox{-\height}{\includegraphics{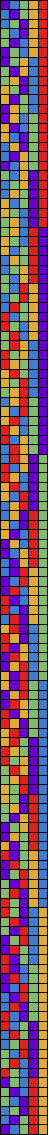}} &
\raisebox{-\height}{\includegraphics{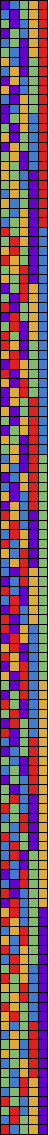}} &
\raisebox{-\height}{\includegraphics{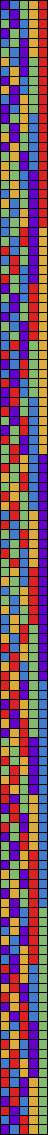}} &
\raisebox{-\height}{\includegraphics{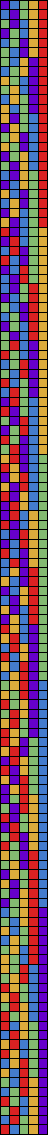}} &
\raisebox{-\height}{\includegraphics{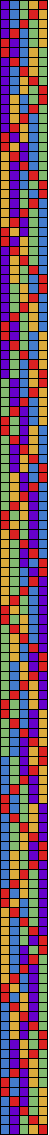}} &
\raisebox{-\height}{\includegraphics{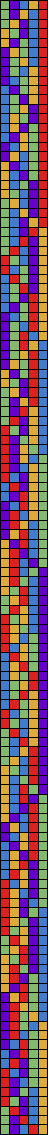}} &
\raisebox{-\height}{\includegraphics{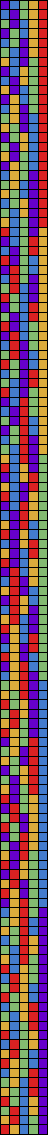}} &
\raisebox{-\height}{\includegraphics{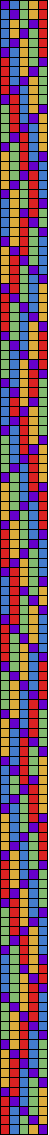}} &
\raisebox{-\height}{\includegraphics{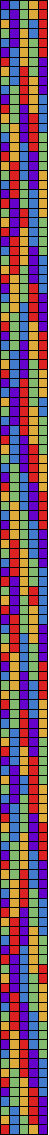}} &
\raisebox{-\height}{\includegraphics{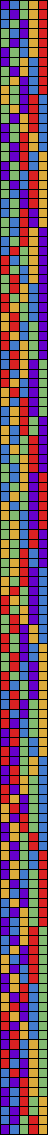}} &
\raisebox{-\height}{\includegraphics{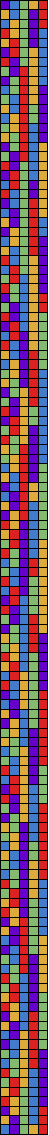}} &
\raisebox{-\height}{\includegraphics{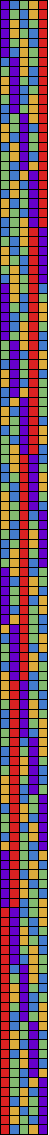}} &
\raisebox{-\height}{\includegraphics{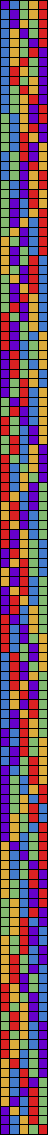}} &
\raisebox{-\height}{\includegraphics{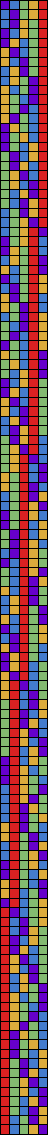}}
\\
(a) & (b) & (c) & (c') & (d) & (e) & (f) & (g) & (g$^{-1}$) & (h) & (h') & (i) & (j) & (j') \\
\end{tabular}
}
\caption{Gray codes for permutations of length~5 (1=purple, 2=blue, 3=green, 4=yellow, 5=red):
(a)~transpositions (Wells);
(b)~transpositions (Heap);
(c)~transpositions (Lipski 10);
(c')~transpositions (Lipski 16);
(d)~adjacent transpositions (SJT);
(e)~doubly-adjacent;
(f)~star transpositions;
(g)~hot potato;
(g$^{-1}$)~inverse of hot potato;
(h)~$(1,2)$-alternating adjacent transpositions;
(h')~$(1,2)$-alternating star transpositions;
(i)~balanced transpositions;
(j)~balanced adjacent transpositions;
(j')~balanced circularly adjacent transpositions;
(figure continues on next page)
}
\label{fig:perm1}
\end{figure}

\begin{figure}
\ContinuedFloat
\centerline{
\setlength{\tabcolsep}{3pt}
\renewcommand{\arraystretch}{1.5}
\begin{tabular}{ccccccc}
\raisebox{-\height}{\includegraphics{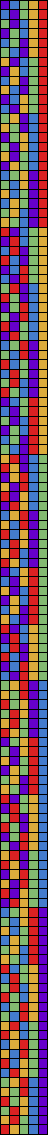}} &
\raisebox{-\height}{\includegraphics{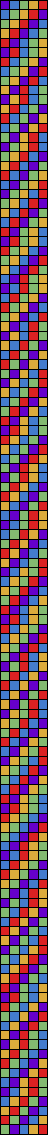}} &
\raisebox{-\height}{\includegraphics{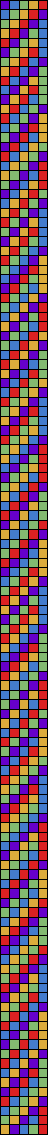}} &
\raisebox{-\height}{\includegraphics{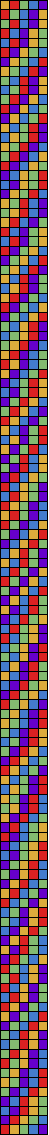}} &
\raisebox{-\height}{\includegraphics{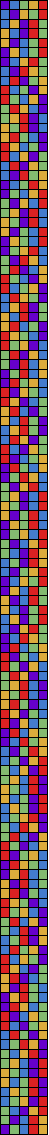}} &
\raisebox{-\height}{\includegraphics{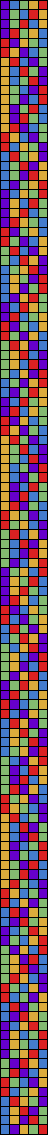}} &
\raisebox{-\height}{\includegraphics{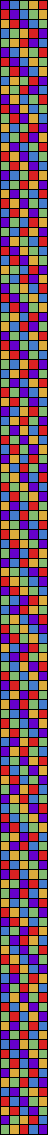}}
\\
(k) & (l) & (m) & (n) & (o) & (o') & (p) \\
\end{tabular}
}
\caption{(continued)
(k)~prefix reversals;
(l)~prefix shifts;
(m)~$\{\sigma_5,\sigma_4\}$;
(n)~$\{\sigma_5,\sigma_2\}$;
(o)~$\{(1\,2)(3\,4),(2\,3)(4\,5),(1\,2)\}$;
(o')~$\{(1\,2)(3\,4),(2\,3)(4\,5),(3\,4)\}$;
(p) derangement order.
}
\end{figure}

Compton and Williamson's~\cite{MR1308693} \emph{doubly-adjacent} Gray code is an ordering of permutations by adjacent transpositions with the additional property that for any two consecutive transpositions~$(i,i+1)$ and~$(j,j+1)$, we have $|i-j|=1$; see Figure~\ref{fig:perm1}~(e).
The inductive construction underlying their proof is considerably more involved than for the Steinhaus-Johnson-Trotter ordering.

Another special type of transpositions are \emph{star transpositions}, i.e., an exchange of the first entry of the permutation with some later entry.
Ehrlich proposed a general method for generating permutations by star transpositions, which he communicated to Martin Gardner in~1987~\cite[Sec.~7.2.1.2]{MR3444818}; see Figure~\ref{fig:perm1}~(f) and also~\cite[Sec.~10]{MR2500822}.
Curiously, the Ehrlich ordering is cyclic for $n=2,3,4,6,12,20,40$, but for no other $n\leq \num{250000}$ \cite[Ex.~58 in Sec.~7.2.1.2]{MR3444818}.

Shen and Williams \cite{DBLP:journals/endm/ShenW13} constructed a \emph{`hot potato'} Gray code for permutations, where the value~1 moves in each step, and it moves by circularly adjacent or semi-adjacent transpositions, i.e., transpositions of the form~$(i,i+1)$ or $(i,i+2)$ modulo~$n$; see Figure~\ref{fig:perm1}~(g).
We can think of the value~1 as being the hot potato that is passed around in a circle of $n$ children, and it is always passed to a neighbor in distance~1 or~2 along the circle.
By considering the inverse permutations of this ordering, the meaning of values and positions is exchanged, so we obtain a star transposition ordering of permutations (the entry at position~1 changes in each step), with the additional property that the value at position~1 always changes by $\pm 1$ or $\pm 2$ modulo~$n$; see Figure~\ref{fig:perm1}~(g$^{-1}$) and cf.~Figure~\ref{fig:perm1}~(f).
Shen and Williams' construction uses the aforementioned doubly-adjacent Gray code of Compton and Williamson as a building block.

The \emph{parity} of a permutation~$\pi$ is the parity of the number of its inversions, where an \emph{inversion} is a pair~$(\pi(i),\pi(j))$ with~$i<j$ and $\pi(i)>\pi(j)$.
Equivalently, the parity of~$\pi$ is given by the parity of the number of transpositions whose product is~$\pi$.
Alternatively, the parity of~$\pi$ is the parity of $n$ minus the number of cycles of~$\pi$.
Note that every transposition changes the parity of a permutation.

For a particular set of admissible transpositions, we consider the \emph{transposition graph} on the vertex~$[n]$, which contains as edges all admissible transpositions~$(i,j)$.
For adjacent transpositions, the transposition graph is a path, whereas for star transpositions, the transposition graph is a star.
Clearly, to be able to generate all permutations by a set of transpositions, the transposition graph must be connected.
Kompel'maher and Liskovec~\cite{MR0498276} and Slater~\cite{MR504868} independently generalized the aforementioned results considerably by proving that all permutations of~$[n]$ can be generated for any transposition tree on~$[n]$.
Even more generally, Tchuente~\cite{MR683982} showed that for any transposition tree on~$[n]$, $n\geq 4$, and for any two permutations of opposite parity, there is a transposition Gray code starting and ending with the two chosen permutations.
As every transposition changes the parity, any transposition Gray code lists permutations with even and odd parity alternatingly.
In graph-theoretic language, the underlying flip graph for any tree of transpositions is bipartite, and Tchuente's result asserts that it is Hamilton-laceable; see also~\cite[Sec.~12]{MR2500822}.

Ruskey and Savage~\cite{MR1201997} showed that for any transposition tree on~$[n]$, $n\geq 5$, and any fixed transposition from that tree, there is a transposition Gray code that uses this fixed transposition in every second step.
In other words, the fixed transposition induces a perfect matching in the corresponding flip graph, and this perfect matching extends to a Hamilton cycle, i.e., the cycle alternatingly uses matching and non-matching edges.
Their Gray code for a path and star on $n=5$ vertices and for the fixed transposition~$(1,2)$ is shown in Figures~\ref{fig:perm1}~(h) and~(h'), respectively.

\solvedproblem{30}{Felsner, Kleist, M\"utze and Sering~{\cite{MR4046775}} raised the question whether there is a transposition Gray code for permutations that uses each of the $\binom{n}{2}$ transpositions equally often.}
This is the analogue of balancing the transition counts in Gray codes for binary strings; recall Section~\ref{sec:transition}.
This problem was solved by Lucca Tiemens~\cite{tiemens_2024a} with a beautiful inductive construction; see Figure~\ref{fig:perm1}~(i).
\openproblem{62}{Similarly, can all $n!$ permutations be listed by adjacent transpositions, using each of the $n-1$ adjacent transpositions equally often?}
While there is no solution for $n=4$, there is one for $n=5$; see Figure~\ref{fig:perm1}~(j).
As a partial result, Gregor, Merino and M\"utze~\cite{MR4747482} constructed a permutation Gray code for odd~$n$ that uses each of the $n$ circularly adjacent transpositions equally often; see Figure~\ref{fig:perm1}~(j').
Tiemens~\cite{tiemens_2024b} complemented this result by providing such a Gray code also for the even values of~$n$.

\openproblem{63}{Furthermore, can all $n!$ permutations be listed by star transpositions, using each of the $n-1$ star transpositions equally often?}
While there is no solution for $n=4$, Gregor, Merino and M\"utze~\cite{MR4747482} found one for~$n=6$.

\openproblem{64}{Lastly, is there a single-track Gray code for permutations by adjacent transpositions, i.e., can all $n!$ permutations be generated by adjacent transpositions such that any two columns in the row-by-row listing of permutations only differ in a cyclic shift?}
For arbitrary transpositions and circularly adjacent transpositions the aforementioned balanced Gray codes also have the single-track property.

\subsection{Reversals, shifts and other generators}
\label{sec:rho-sigma}

Ord-Smith~\cite{DBLP:journals/cacm/Ord-Smith67} devised a method to generate all permutations by \emph{prefix reversals}; see Figure~\ref{fig:perm1}~(k).
His listing also has the genlex property, and the average length of reversed prefixes approaches~$e=2.71828\ldots$ (a constant, which allows for an efficient algorithm).
This ordering was rediscovered by Zaks~\cite{MR753548}, who was motivated by sorting problems involving stacks of pancakes; see also~\cite[Sec.~11]{MR2500822}.
A prefix reversal corresponds to flipping over some of the pancakes at the top of the stack.
In the following, we write $\rho_k$ for the operation that reverses the first $k$ entries of a permutation, i.e., in cycle notation we have $\rho_k:=(1\;k)(2\;k-1)\cdots$.

\emph{Prefix shifts} have also been investigated heavily for permutation generation.
Unlike transpositions and prefix reversals discussed before, prefix shifts are in general not involutions.
Formally, for any $2\leq k\leq n$, we write $\sigma_k$ for the operation that cyclically left-shifts the first $k$ entries of a permutation by one position, i.e., $\sigma_k=(k\;k-1\cdots 2\,1)$.
Note that $\sigma_2=\rho_2$ simply transposes the first two entries.
Corbett~\cite{DBLP:journals/tpds/Corbett92} showed that all permutations can be generated cyclically by prefix shifts of arbitrarily length, i.e., using the generators $\{\sigma_n,\sigma_{n-1},\ldots,\sigma_2\}$; see Figure~\ref{fig:perm1}~(l).
Compton and Williamson~\cite{MR1308693}, as an application of their doubly-adjacent Gray code, showed that we can also use the generators~$\{\sigma_n,\sigma_n^{-1},\sigma_2\}$.
Stevens and Williams~\cite{MR2991865} showed the same result for the generators $\{\sigma_n,\sigma_3,\sigma_2\}$, and Ruskey and Williams~\cite{MR2682614} proved it for~$\{\sigma_n,\sigma_{n-1}\}$; see Figure~\ref{fig:perm1}~(m).
In a recent breakthrough, Sawada and Williams~\cite{MR4060409} give a simple and explicit rule to generate all permutations of length~$n$ using only the two generators~$\{\sigma_n,\sigma_2\}$, and this listing is cyclic if and only if $n$ is odd; see Figure~\ref{fig:flip}~(e) and Figure~\ref{fig:perm1}~(n).
Subsequently, a simple loopless algorithm for this problem was devised by Lipt{\'{a}}k, Masillo, Navarro and Williams~\cite{DBLP:conf/spire/LiptakMNW23}.
Instead of cyclically shifting the array that holds the permutation when applying~$\sigma_n$, they increment an index to what is considered the `first' entry.
This problem was known as \emph{sigma-tau problem}, and was first mentioned in Nijenhuis and Wilf's book~\cite[Ex.~6]{MR0396274}.
All of the preceding results mentioned in this paragraph assert the existence of Hamilton paths or cycles in the Cayley digraph~$\vv{\Gamma}(S_n,X)$ for various sets~$X$ of shifts as generators.

We may also consider other generators of the symmetric group~$S_n$.
In this direction, Rapaport-Strasser~\cite{MR111702} proved that the Cayley graph~$\Gamma(S_n,X)$ has a Hamilton cycle for the three involutions $X=\{(1\,2)(3\,4)(5\,6)\cdots,(2\,3)(4\,5)(6\,7)\cdots,(1\,2)\}$.
Knuth~\cite[Ex.~69 in Sec.~7.2.1.2]{MR3444818} devised an efficient algorithm for computing such a cycle for the three involutions $X=\{(1\,2)(3\,4)(5\,6)\cdots,(2\,3)(4\,5)(6\,7)\cdots,(3\,4)(5\,6)(7\,8)\cdots\}$, which can be adapted to also work for Rapaport-Strasser's set of generators.
The corresponding Gray codes for $n=5$ are shown in Figures~\ref{fig:perm1}~(o) and~(o'), respectively.
In Knuth's Gray code, the largest value~$n$ moves back and forth between the first and last position, in exactly the same way as in the SJT listing (but the other values move differently).
The fact that $\Gamma(S_n,X)$ has a Hamilton cycle for the two aforementioned sets~$X$ of generators follows from the general fact that the Cayley graph of any group generated by three involutions~$\rho,\sigma,\tau$ has a Hamilton cycle if $\rho\tau=\tau\rho$; see \cite[Lemma~1]{MR2548568}.
Note that by reordering positions, Rapaport-Strasser's and Knuth's generators are equivalent to~$X=\{\rho_n,\rho_{n-1},(1\,n)\}$ and $X=\{\rho_n,\rho_{n-1},\rho'\}$, respectively, where $\rho'$ reverses the middle $n-2$ entries of a permutation.

Savage~\cite{MR1069115} showed that for any $n\geq k\geq 2$ and $k\neq 3$, all permutations of length~$n$ can be listed cyclically so that consecutive permutations differ in exactly $k$ positions.
For $k=n$ such a listing is sometimes called a \emph{derangement order}.
Figure~\ref{fig:perm1}~(p) shows an example for this, using a construction due to Yarbrough described in~\cite{MR1491049}.
Rasmussen and Savage~\cite{MR1298976} later generalized this result considerably, by showing that the corresponding flip graph is Hamilton-connected, even if we impose the additional constraint that the $k$ modified positions are consecutive.

\openproblem{31}{Despite all of the aforementioned results, in general it is wide open for which sets of generators~$X$ the Cayley graph~$\Gamma(S_n,X)$ or the Cayley digraph~$\vv{\Gamma}(S_n,X)$ have a Hamilton cycle or path.}
Generalizing the sigma-tau problem, we may consider $X=\{\sigma_n,(1\;k)\}$ or $X=\{\sigma_n,\sigma_k\}$ for suitable values of~$n$ and $2\leq k\leq n$.
Specifically, the case $X=\{\sigma_n,\sigma_3\}$ for even~$n$ was suggested by Stevens and Williams~\cite{MR2991865}.
Furthermore, Savage~\cite{MR1491049} mentions $X=\{\rho_n,\rho_{n-1},(n-1\;n)\}$ as an interesting open problem.
Lastly, Bass and Sudborough~\cite{DBLP:journals/jpdc/BassS03} and Sawada and Williams~\cite{MR3426940} suggest the instances $X=\{\rho_n,\rho_{n-1},\rho_{n-2}\}$ and~$X=\{\rho_n,\rho_{n-1},\rho_2\}$.

\subsection{Restricted permutations}
\label{sec:perm-restr}

\begin{figure}
\centerline{
\setlength{\tabcolsep}{3pt}
\renewcommand{\arraystretch}{1.5}
\begin{tabular}{ccccccccccccccc}
(a) & (a') & (a'') & (b) & (c) & (d1') & (d1) & (d2) & (e1) & (e2) & (f) & (g) & (h) & (h') & (i) \\[-5mm]
\raisebox{-\height}{\includegraphics{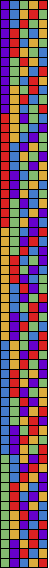}} &
\raisebox{-\height}{\includegraphics{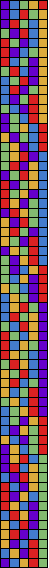}} &
\raisebox{-\height}{\includegraphics{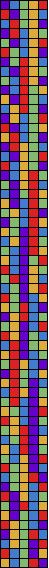}} &
\raisebox{-\height}{\includegraphics{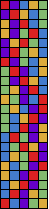}} &
\raisebox{-\height}{\includegraphics{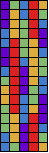}} &
\raisebox{-\height}{\includegraphics{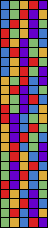}} &
\raisebox{-\height}{\includegraphics{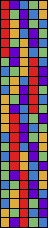}} &
\raisebox{-\height}{\includegraphics{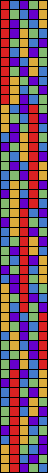}} &
\raisebox{-\height}{\includegraphics{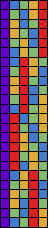}} &
\raisebox{-\height}{\includegraphics{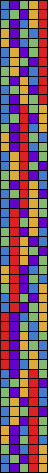}} &
\raisebox{-\height}{\includegraphics{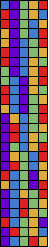}} &
\raisebox{-\height}{\includegraphics{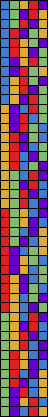}} &
\raisebox{-\height}{\includegraphics{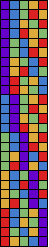}} &
\raisebox{-\height}{\includegraphics{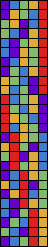}} &
\raisebox{-\height}{\includegraphics{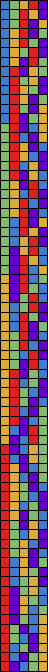}} \\
(j) & (k) & (l) & & & & & & & & & & & & \\[-5mm]
\raisebox{-\height}{\includegraphics{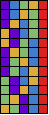}} &
\raisebox{-\height}{\includegraphics{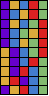}} &
\raisebox{-\height}{\includegraphics{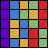}} & & & & & & & & & & & &
\end{tabular}
}
\caption{Gray codes for restricted permutations of length~5:
(a)--(a'')~even permutations;
(b)~five inversions;
(c) alternating;
(d1')+(d1) one cycle (=cyclic permutations); (d2) two cycles;
(e1)+(e2) one/two left-to-right minima;
(f) one exceedance;
(g)~derangements;
(h)+(h')~involutions;
(i) indecomposable;
(j)~rosary permutations;
(k)~linear extensions of~\includegraphics[page=1,scale=0.5]{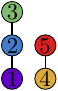};
(l)~linear extensions of~\includegraphics[page=2,scale=0.5]{poset}.
}
\label{fig:perm2}
\end{figure}

\emph{Even permutations} are permutations with even parity, i.e., with an even number of inversions.
As mentioned before, any transposition Gray code visits permutations with even and odd parity alternatingly, so from any transposition Gray code for all permutations we can obtain a 2-transposition Gray code for even permutations (or odd permutations, or permutations with an even number of cycles, or an odd number of cycles); cf.~\cite{MR2817083}.
Savage~\cite{MR1069115} showed that they can be listed by applying a \emph{3-rotation}, i.e., a product of two transpositions that share one element, in each step.
Recall that each transposition changes the parity, so single transpositions cannot be used.
Her listing is shown in Figure~\ref{fig:perm2}~(a), and it is genlex.
She also came up with a listing that uses only 3-rotations of consecutive entries, to the left or right by one position if the middle entry is at an odd or even position, respectively (these are only $n-2$ different operations in total).
This listing is described in~\cite[Ex.~103 in Sec.~7.2.1.2]{MR3444818}, and it is shown in Figure~\ref{fig:perm2}~(a').
Gould and Roth~\cite{MR900932} proved that for $n\geq 5$ all even permutations can be listed cyclically by 3-rotations of the form~$(1\,i\,n)$ for some $1<i<n$, i.e., in every step, the first, $i$th and last entry of the permutation rotate cyclically; see Figure~\ref{fig:perm2}~(a'').
Equivalently, their result asserts that for $n\geq 5$ the Cayley digraph~$\vv{\Gamma}(A_n,X)$ of the alternating group~$A_n$ with generators $X=\{(1\,i\,n)\mid 1<i<n\}$ has a Hamilton cycle.
By reordering positions, Gould and Roth's generators are equivalent to~$X=\{(1\,2\,i)\mid 3\leq i\leq n\}$, which can be seen as `star-like' 3-rotations, analogous to Ehrlich's star transpositions for~$S_n$.
Holroyd~\cite{MR3599935} proved that for all odd~$n\neq 5$, $\vv{\Gamma}(A_n,X)$ with generators $X=\{\sigma_3,\sigma_5,\ldots,\sigma_n\}$, i.e., odd prefix shifts, has a Hamilton cycle.

Walsh~\cite{MR1772772} constructed a Gray code for permutations of length~$n$ with a fixed number~$k$ of inversions, $0\leq k\leq \binom{n}{2}$, which uses a product of two adjacent transpositions in each step; see Figure~\ref{fig:perm2}~(b).
He derived this from a genlex Gray code for integer compositions $b_1+\cdots+b_n=k$ with bounded parts~$b_i\leq i-1$ for all $i\in [n]$, by interpreting the composition~$(b_1,\ldots,b_n)$ as the inversion vector of a permutation.
The \emph{inversion vector} of a permutation~$\pi$ of length~$n$ is a vector $(b_1,\ldots,b_n)$, where $b_i$ counts the number of values in~$\pi$ that are to the right of the value~$i$ and smaller than~$i$.
Consequently, any permutation with exactly $k$ inversions has an inversion vector $(b_1,\ldots,b_n)$ with $b_1+\cdots+b_n=k$.
In Walsh's composition Gray code, consecutive compositions differ by changing one entry~$b_i$ by~$+1$ and another entry~$b_j$ by~$-1$, and this translates into a product of two adjacent transpositions in the permutation: Specifically, $i$ is transposed with the first entry to its left smaller than~$i$, and $j$ is transposed with the first entry to its right larger than~$j$.
Vajnovszki and Vernay~\cite{MR2817083} generalized this Gray code to permutations for which the number of inversions is in any fixed interval.

An \emph{alternating permutation} is a permutation~$\pi$ such that $\pi(1)<\pi(2)>\pi(3)<\pi(4)>\cdots$.
Ruskey~\cite{MR1048093} showed that for odd~$n$, all alternating permutations can generated by transpositions (which must clearly be non-adjacent); see Figure~\ref{fig:perm2}~(c).
\openproblem{32}{He raised the question whether such a Gray code exists also for even values of~$n$.}
A more general notion are \emph{permutations with a signature}, which is sequence~$s=(s_1,\ldots,s_{n-1})$ with~$s_i\in\{+,-\}$, and the requirement is that $\pi(i)<\pi(i+1)$ if $s_i=+$ and $\pi(i)>\pi(i+1)$ if $s_i=-$ for all $1\leq i<n$.
Note that alternating permutations arise from the signature $s=(+,-,+,-,\ldots)$.
Efficient generation algorithms for permutations with any fixed signature have been devised by Roelants van Baronaigien and Ruskey~\cite{MR1158770}, and later by Korsh~\cite{MR1855048}.
\openproblem{33}{However, none of these methods is based on a Gray code, so the question for the existence of a (transposition?) Gray code for permutations with a given signature remains open in general.}

A \emph{cycle} in a permutation~$\pi$ is sequence $i_1,\ldots,i_k$ such that $\pi(i_j)=i_{j+1}$ for $j=1,\ldots,k-1$ and $\pi(i_k)=i_1$.
Baril~\cite{MR2316665} devised a Gray code to generate permutations of length~$n$ and any fixed number~$k$ of cycles, $1\leq k\leq n$, with 3-rotations as transitions.
Figure~\ref{fig:perm2}~(d1)+(d2) shows these orderings for $n=5$ and $k=1,2$.
Another listing for \emph{cyclic permutations}~($k=1$) based on the SJT algorithm was proposed by Knuth~\cite[Ex.~96 in Sec.~7.2.1.2]{MR3444818}; see Figure~\ref{fig:perm2}~(d1').

A \emph{left-to-right minimum} in a permutation~$\pi$ is an entry $\pi(i)$, $i\in[n]$, for which $\pi(j)>\pi(i)$ for all $1\leq j<i$.
In particular, the first entry is always a left-to-right minimum.
Baril~\cite{MR3100174} developed a Gray code for permutations of length~$n$ with a fixed number~$k$ of left-to-right minima, $1\leq k\leq n$, which uses transpositions or 3-rotations as transitions; see Figure~\ref{fig:perm2}~(e1)+(e2).
It is well known that for any~$n$ and~$k$, permutations of~$[n]$ with $k$ left-to-right minima and permutations of~$[n]$ with $k$ cycles are equinumerous, and both are counted by the Stirling numbers of the first kind [OEIS~A132393].

An \emph{exceedance} in a permutation is an entry~$\pi(i)$, $i\in[n]$, with $\pi(i)>i$.
Baril~\cite{MR2566173} describes a Gray code for permutations~$\pi$ that have exactly one exceedance, where consecutive permutations differ by a transposition or a 3-rotation; see Figure~\ref{fig:perm2}~(f).
\openproblem{34}{The question whether there is a Gray code for permutations with a fixed number~$k$ of exceedances is open for~$k\geq 2$.}

A \emph{derangement} is a permutation $\pi:[n]\rightarrow [n]$ that has no fixpoints, i.e., $\pi(i)\neq i$ for all $i=1,\ldots,n$ [OEIS~A000166].
Korsh and LaFollette~\cite{MR2050887} present a Gray code for derangements, using either one transposition or a 3-rotation in each step.
Their listing is genlex; see Figure~\ref{fig:perm2}~(g).
They also generalized their Gray code to visit all permutations with a given lower and upper bound on the number of fixpoints; see also~\cite{MR2064117}.

An \emph{involution} is a permutation $\pi:[n]\rightarrow [n]$ in which all cycles have length~1 or~2, i.e., for all $i=1,\ldots,n$ we have $\pi(i)=i$ or $\pi(i)=j$ and $\pi(j)=i$ for some $j\in[n]\setminus\{i\}$ [OEIS~A000085].
Walsh~\cite{MR1821628} presents a Gray code for involutions, which uses either one or two transpositions or a 3-rotation in each step, and this listing is genlex; see Figure~\ref{fig:perm2}~(h).
He also extended this result to involutions with a given number of fixpoints, in particular for fixpoint-free involutions (in the latter case no 3-rotations are needed).
Gutierres, Mamede and Santos~\cite{MR3947655} developed a Gray code for involutions that uses only one transposition or a 3-rotation in each step, i.e., consecutive permutations differ in at most three entries; see Figure~\ref{fig:perm2}~(h').

A permutation~$\pi$ of~$[n]$ is \emph{indecomposable} if no proper prefix of~$\pi$ is a permutation, i.e., $\pi([i])\neq [i]$ for all $1\leq i<n$.
King~\cite{MR2212519} constructed a transposition Gray code for indecomposable permutations; see Figure~\ref{fig:perm2}~(i).
\openproblem{35}{He also raised the question whether they admit an adjacent transposition Gray code.}

A \emph{rosary permutation} is bracelet of permutations, i.e., an equivalence class of permutations under cyclic shifts and reversal.
We can pick as representatives the permutations that have $n$ at their last position and the value~1 before the value~2, and then these representatives can be listed in adjacent transposition order by considering the first $(n-1)!/2$ permutations in the SJT listing for~$n$; see Figure~\ref{fig:perm2}~(j) and~\cite{DBLP:journals/cacm/Roy73}.

\subsection{Pattern-avoiding permutations}
\label{sec:avoid}

Given a permutation~$\pi$ of length~$n$ and a permutation~$\tau$ of length~$k$, where $k\leq n$, we say that \emph{$\pi$ contains the pattern~$\tau$}, if there is a subsequence of indices $i_1<\cdots<i_k$ such that the elements of $\pi(i_1),\ldots,\pi(i_k)$ appear in the same relative order as in~$\tau$.
Otherwise, we say that \emph{$\pi$ avoids~$\tau$}.
In a \emph{vincular} pattern~$\tau$, two consecutive entries are underlined, and for $\pi$ to contain~$\tau$, these entries must match consecutive entries of~$\pi$.
For example, $\pi=\colorbox{black!30}{\!53\!}1\colorbox{black!30}{\!4\!}6\colorbox{black!30}{\!2\!}$ contains $\tau=4231$, as evidenced by the shaded entries of~$\pi$, whereas $\pi$ avoids~$\tau'=4\underline{23}1$ and also $\tau''=2413$.
The study of pattern-avoidance in permutations is a major branch of combinatorics, and the reason is that a vast number of combinatorial objects are bijectively equivalent to pattern-avoiding permutations.
For example, 321-avoiding permutations and 231-avoiding permutations are both counted by the Catalan numbers, and are thus bijectively equivalent to many other Catalan objects~\cite{MR3467982}.
For any pattern~$\tau$, we write $S_n(\tau)$ for the set of permutations of~$[n]$ that avoid~$\tau$.
Moreover, for any set~$X$ of permutation patterns (possibly of different lengths), we write $S_n(X)$ for the permutations of~$[n]$ that avoid each of the patterns from~$X$, i.e., $S_n(X):=\bigcap_{\tau\in X}S_n(\tau)$.

\begin{figure}[t!]
\centerline{
\setlength{\tabcolsep}{3pt}
\renewcommand{\arraystretch}{1.5}
\begin{tabular}{ccccccccccccccc}
(a) & (b) & (c) & (d) & (e) & (f) & (g) & (h) & (i) & (j) & (k) & (l) & (m) & (n) & (n') \\[-5mm]
\raisebox{-\height}{\includegraphics{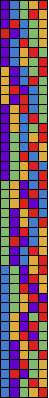}} &
\raisebox{-\height}{\includegraphics{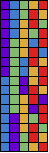}} &
\raisebox{-\height}{\includegraphics{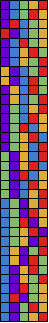}} &
\raisebox{-\height}{\includegraphics{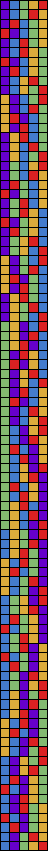}} &
\raisebox{-\height}{\includegraphics{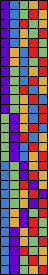}} &
\raisebox{-\height}{\includegraphics{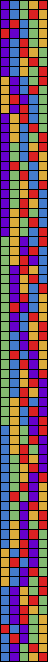}} &
\raisebox{-\height}{\includegraphics{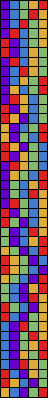}} &
\raisebox{-\height}{\includegraphics{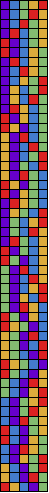}} &
\raisebox{-\height}{\includegraphics{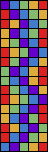}} &
\raisebox{-\height}{\includegraphics{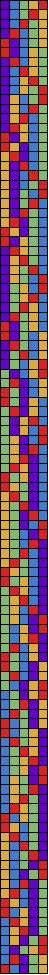}} &
\raisebox{-\height}{\includegraphics{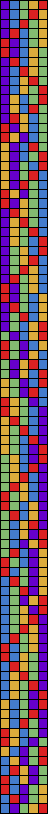}} &
\raisebox{-\height}{\includegraphics{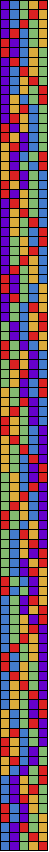}} &
\raisebox{-\height}{\includegraphics{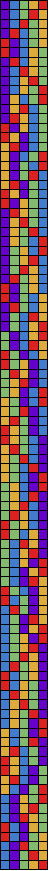}} &
\raisebox{-\height}{\includegraphics{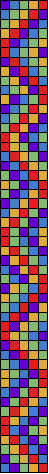}} &
\raisebox{-\height}{\includegraphics{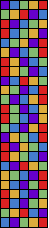}}
\end{tabular}
}
\caption{Gray codes for pattern-avoiding permutations of length~5 ((a)--(f) are for regular patterns~$X$, and (g)--(m) are for tame patterns~$X$):
(a)~$X=\{321\}$;
(b)~$X=\{321,312\}$;
(c)~$X=\{321,3412\}$;
(d)~$X=\{4321,4312\}$;
(e)~$X=\{321,3412,4123\}$;
(f)~$X=\{4321,4312,4132,4231\}$;
(g)~$X=\{231\}$;
(h)~$X=\underline{23}1$;
(i)~$X=\{231,132\}$ (permutations without peaks);
(j)~$X=\{2143\}$ (vexillary permutations);
(k)~$X=\{2143,3412\}$ (skew-merged permutations);
(l)~$X=\{2413,3142\}$ (separable permutations);
(m)~$X=\{2\underline{41}3,3\underline{14}2\}$ (Baxter permutations);
(n)~stamp foldings;
(n')~semi-meanders.
}
\label{fig:perm3}
\end{figure}

Dukes, Flanagan, Mansour and Vajnovszki~\cite{MR2412243} developed Gray codes for pattern-avoiding permutations.
To state their result, we call a set~$X$ of permutation patterns \emph{regular}, if every $\tau\in X$ of length~$k$ satisfies the following two conditions:
(i)~two entries smaller than~$k$ are to the right of~$k$ in~$\tau$; (ii)~if the entry~$k$ is not at the leftmost position of~$\tau$, then~$X$ contains at least one subpattern of~$\tau'$, where $\tau'$ is obtained from~$\tau$ by transposing~$k$ with the entry to its left.
For example, if $\tau=436512\in X$, then by~(ii) one subpattern of~$\tau':=463512$ must be in~$X$, for example $\sigma:=4632=3421\in X$.
Consequently, if $\sigma=3421\in X$, then by~(ii) one subpattern of~$\sigma':=4321$ must be in~$X$, for example $\sigma'\in X$ itself.
Overall, $X=\{\tau,\sigma,\sigma'\}=\{436512,3421,4321\}$ is a regular set of patterns.
Note that if all patterns in~$X$ have the largest entry at the leftmost position, then condition~(ii) is void, and (i) requires that all patterns have length at least~3.
The paper~\cite{MR2412243} shows how to list, for any regular set of patterns~$X$, all pattern-avoiding permutations~$S_n(X)$ such that at most five entries change in each step.
This construction was refined by Baril~\cite{MR2530633} to most three changes in each step (which must be a 3-rotation).
Examples of regular patterns to which this result applies are $X=\{321\}$ or $X=\{312\}$ (counted by the Catalan numbers), $X=\{321,312\}$ (counted by $2^{n-1}$), $X=\{321,3412\}$ or $X=\{321,4123\}$ or $X=\{321,4321\}$  (counted by the even-index Fibonacci numbers), $X=\{4321,4312\}$ or $X=\{4213,4123\}$ or $X=\{4132,4231\}$ (counted by the Schr\"oder numbers [OEIS~A006318]), $X=\{321,3412,4123\}$ or $X=\{312,3421,4321\}$ (counted by the Pell numbers [OEIS~A000129]), and $X=\{4321,4312,4132,4231\}$ or $X=\{4213,4123,4132,4231\}$ (counted by the central binomial coefficients $\binom{2n-2}{n-1}$); see Figure~\ref{fig:perm3}~(a)--(f).

\begin{figure}
\includegraphics{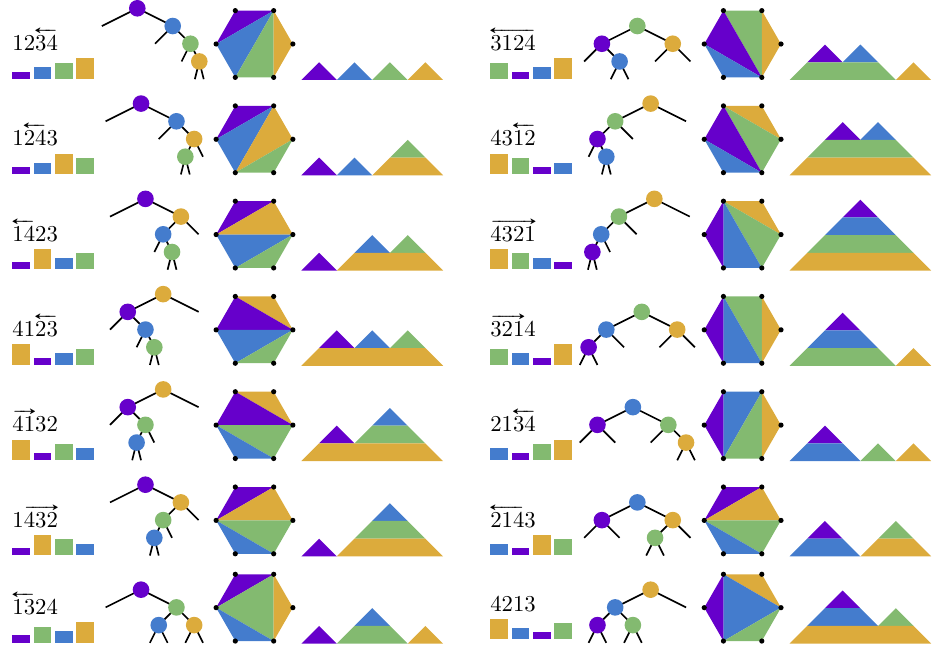}
\caption{231-avoiding permutations of length $n=4$ generated by minimal jumps and the resulting Gray codes for Catalan objects (binary trees, triangulations, Dyck paths), which coincide with the ordering due to Lucas, Roelants van Baronaigien, and Ruskey~\cite{MR1239499}.
Minimal jumps in the permutations correspond to rotations in the binary trees, to flipping a diagonal in the triangulations, and to flipping a hill in the Dyck paths.}
\label{fig:231cat}
\end{figure}

Hartung, Hoang, M\"utze and Williams~\cite{MR4391718} propose a different Gray code for generating pattern-avoiding permutations, using cyclic shifts of substrings.
Specifically, a \emph{right jump of the value~$a_i$} in a permutation $\pi=a_1,\ldots,a_n$ cyclically shifts a substring $a_i,\ldots,a_j$ with $a_i>a_{i+1},\ldots,a_j$ one position to the left, yielding $a_{i+1},\ldots,a_j,a_i$, whereas a \emph{left jump of the value~$a_j$} shifts a substring $a_i,\ldots,a_j$ with $a_j>a_i,\ldots,a_{j-1}$ one position to the right, yielding $a_j,a_i,\ldots,a_{j-1}$, and the \emph{distance} of the jump is the number~$j-i$.
A permutation pattern~$\tau$ of length~$k$ is called \emph{tame}, if it does not have the largest value~$k$ at the leftmost or rightmost position.
A vincular pattern is tame if in addition to this condition the largest value~$k$ is part of the vincular pair.
For example, $\tau=231$ and $\tau=2\underline{41}3$ are tame, whereas $\tau=321$ and $\tau=24\underline{13}$ are not.
The paper shows how to list, for any set of tame patterns~$X$, all pattern-avoiding permutations~$S_n(X)$ by jumps.
Examples of pattern-avoiding permutations to which this result applies are permutations without peaks ($X=\{231,132\}$), vexillary permutations ($X=\{2143\}$), skew-merged permutations ($X=\{2143,3412\}$), separable permutations ($X=\{2413,3142\}$), and Baxter permutations ($X=\{2\underline{41}3,3\underline{14}2\}$); see Figure~\ref{fig:perm3}~(g)--(m).
The listings are obtained from a simple greedy algorithm, which repeatedly performs a minimal distance jump with the largest possible value, so that a new permutation in the list is generated.
The main feature of this approach is that applying standard bijections to these jump orderings of permutations yields Gray codes for many different combinatorial objects; see Figure~\ref{fig:231cat}.
This generation approach therefore generalizes several known Gray codes, such as the BRGC (for $X=\{231,132\}$), the SJT algorithm (for $X=\emptyset$ minimal jumps are adjacent transpositions), the Gray code for binary trees by rotations due to Lucas, Roelants van Baronaigien and Ruskey~\cite{MR1239499} (for $X=\{231\}$), and the Gray code for set partitions by element exchanges described by Kaye~\cite{MR443423} (for $X=\{\underline{23}1\}$).
It also yields Gray codes for many different classes of rectangulations~\cite{MR4598046}, and it produces Hamilton paths and cycles on large classes of polytopes~\cite{MR4344032,MR4415121,MR4614413}.
The reason why minimal jumps are natural becomes clear when considering the inversion vector of the permutations.
Observe that a right jump of~$a_i$ as indicated in the definition above decreases the $a_i$th entry of the inversion vector of~$\pi$ by~$j-i$, whereas a left jump of~$a_j$ as indicated before increases the $a_j$th entry by~$j-i$ (all values that are jumped across are smaller than the value that jumps).
Consequently, minimal jumps change the inversion vector in a single entry by a minimal amount, and a jump Gray code of permutations corresponds to a genlex ordering of the inversion vectors.

\subsection{Linear extensions and ideals of posets}
\label{sec:posets}

A \emph{poset~$P$} is a reflexive, antisymmetric and transitive binary relation.
A \emph{linear extension} of~$P$ is a total order (=chain) that extends~$P$.
For further terminology regarding posets used below, see for example~\cite{MR1169299}.
Linear extensions are sometimes also referred to as \emph{topological sorts}, and they generalize many interesting combinatorial objects.

For example, if $P$ is an antichain on~$[n]$, then the linear extensions are all $n!$ permutations of~$[n]$.
More generally, if~$P$ is a set of $k$ disjoint chains of sizes~$c_1,\ldots,c_k$, then the linear extensions are multiset permutations with frequencies~$c_1,\ldots,c_k$ (see Section~\ref{sec:multiperm}), obtained by identifying the elements of each chain in the linear extension; see Figure~\ref{fig:perm2}~(k).
Furthermore, linear extensions of a $2\times k$ grid tilted by 45~degrees correspond to Catalan objects.
Indeed, if we regard the $k$ elements in the lower chain as opening brackets and the $k$ elements of the upper chain as closing brackets, then every linear extension corresponds to a well-formed parenthesis expression with $k$ pairs of parentheses.
By subdividing the lower chain and appending additional minima, this construction can be generalized to encode $b$-ary Catalan objects for any $b\geq 2$.
Moreover, note that if the poset~$P$ is a fence, i.e., $1<2>3<4>\cdots$, then its linear extensions are the inverses of alternating permutations.
Similarly, all Young tableaux of a given shape arise as the inverse permutations of linear extensions of the tableau turned so that its corner is the minimum; see Figure~\ref{fig:perm2}~(l).
Lastly, partitions of the ground set~$[k\cdot c]$ into $k$ sets of size~$c$ are obtained as inverse permutations of the linear extensions of the comb poset, obtained from $k$ chains of size~$c$ by combining their minima into a chain of size~$k$.
This idea can be generalized to partitions of~$[n]$ into subsets of arbitrary prescribed sizes~$c_1,\ldots,c_k$, by forming one comb for every group of identical integers~$c_i$ as before, and taking the disjoint union of those combs.

Clearly, the flip graph of linear extensions under transpositions is bipartite, with the partition classes given by the parity of the permutations.
We say that a poset is \emph{balanced} if the two partition classes have the same size.
\openproblem{36}{Ruskey~{\cite{ruskey_1988}} conjectured that all balanced posets admit a transposition Gray code, and he also asked whether an adjacent transposition Gray code exists; see also~\mbox{\cite[Problem~485]{MR2548549}}.
This conjecture still stands as one of the most challenging open problems in this area, and has been confirmed for several special cases.}
In~\cite{MR1142265}, Ruskey also raised a seemingly easier version of his conjecture, asking for a perfect matching in the flip graph of balanced posets.
Ruskey's original conjecture is known to be true for the $2\times k$ grid~\cite{MR1041167}, which is balanced if and only if $k$ is even.
It is also known for $k$ disjoint chains of sizes~$c_1,\ldots,c_k$, $k\geq 2$, which form a balanced poset if and only if at least two of the~$c_i$ are odd~\cite{MR1142265}, and this positive result also holds for the more restrictive requirement of adjacent transpositions~\cite{MR1157583}.
Ruskey~\cite{MR1142265} also showed that if~$P$ is a forest, i.e., the Hasse diagram of~$P$ is a collection of rooted trees, then $P$ is balanced if and only if the forest does not admit a matching in which all but at most one of the trees' roots are matched.
To state the next result, we say that a poset is \emph{ranked}, if it admits a rank function that assigns rank~0 to all maximal elements and the rank changes by~$-1$ when moving down any cover relation.
Pruesse and Ruskey~\cite{MR1105947} confirmed the conjecture for ranked posets~$P$ in which every non-minimal element has at least two lower covers or every non-maximal element has at least two upper covers.
Using the results from~\cite{MR1201997}, this statement even holds for adjacent transpositions.
This confirmed case of the conjecture generalizes Ruskey's~\cite{MR1048093} earlier transposition Gray code for alternating permutations for odd~$n$, obtained when $P$ is an odd fence.
It also applies to forest posets in which no vertex has exactly one child, i.e., leaves have none, and all other vertices have at least two.
West~\cite{MR1214892} gives another proof for such forest posets, by showing that if the linear extensions of~$P$ admit an adjacent transposition Gray code, then the poset $P\wedge Q$ also has this property, which is obtained by making a single element of~$P$ cover an antichain~$Q$ of size at least~2.
In a similar vein, Stachowiak~\cite{MR1157583} showed that if a poset~$P$ has an even number of linear extensions that can be listed in adjacent transposition order, then for any poset~$Q$, the linear extensions of the disjoint union $P\cup Q$ also admit an adjacent transposition Gray code.
Another relevant result due to Pruesse and Ruskey~\cite{MR1267216} is that the linear extensions of any poset~$P$ can be generated by one or a product of two adjacent transpositions in each step, and this scheme can be implemented efficiently (cf.~\cite{MR1336537,MR1922914}); see Figure~\ref{fig:perm2}~(k)+(l).
They proved this by showing that the Cartesian product of the linear extension flip graph~$G(P)$ with a single edge, i.e., the \emph{prism} of~$G(P)$, has a Hamilton cycle.
As~$G(P)$ is bipartite, this easily implies that the square~$G^2(P)$ has a Hamilton cycle.

\begin{figure}[h!]
\includegraphics{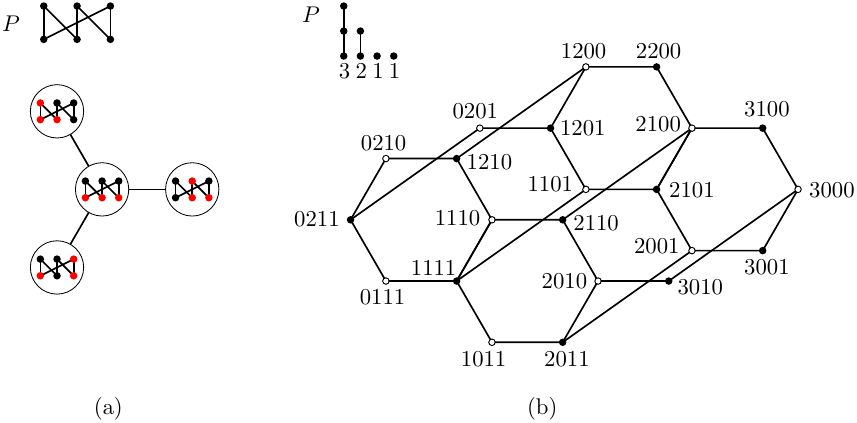}
\caption{Counterexamples to two conjectures of Pruesse and Ruskey~\cite{MR1267190}.
In (b), the ideals of the poset~$P$ are encoded by specifiying the number of elements taken from each of the four chains of~$P$.}
\label{fig:counter}
\end{figure}

An \emph{ideal} of a poset~$P$ is a downward-closed subset of~$P$.
Koda and Ruskey~\cite{MR1231447} constructed a Gray code for the ideals of a forest poset, where consecutive ideals differ by exactly one element.
Pruesse and Ruskey~\cite{MR1267190} showed that the ideals of any poset can be listed by changing (i.e., adding or removing) at most two elements in each step, by proving that the prism of the ideal flip graph admits a Hamilton cycle.
They also proved that such a listing exists for the ideals of fixed size of a series-parallel poset.
Habib, Nourine and Steiner~\cite{MR1474594} generalized this result about ideals of fixed size to interval orders, and they showed that all ideals of a poset can be listed in non-decreasing order of cardinality by changing at most two elements in each step.
Pruesse and Ruskey~\cite{MR1267190} also conjectured that for any poset, its ideals of any fixed size can be listed by element exchanges.
Furthermore, they conjectured that the ideals of size~$k$ or $k+1$ admit a Gray code by element addition/removal, provided that there is the same number of ideals of each of the two sizes.
This generalizes the middle levels conjecture discussed in Section~\ref{sec:range}, obtained when the poset is an antichain of size~$2k+1$.
It turns out that there are easy counterexamples for both of these conjectures; see Figure~\ref{fig:counter}~(a) and~(b), respectively.
For the first one, take the 6-element crown poset and its 3-element ideals. The flip graph formed by element exchanges is a star with 3 rays, so it does not have a Hamilton path.
For the second one, consider the 7-element poset formed by four chains of sizes $3,2,1,1$, which has 11 ideals of size~3 and~4 each, but the bipartite graph formed by them through element addition/removal does not admit a Hamilton path.
\openproblem{65}{Nevertheless, we may still ask, for which posets can the ideals of fixed size be listed by element exchanges, or the ideals of sizes~$k$ or~$k+1$ be listed by element addition/removal?}

\openproblem{37}{Also, it is a long-standing open problem whether the ideals of any poset can be generated in constant amortized-time per ideal, either in a Gray code fashion or not.}

\subsection{Multiset permutations}
\label{sec:multiperm}

A useful generalization of permutations are multiset permutation, which allow repetition of symbols.
Formally, a multiset permutation over the alphabet~$[k]$ with frequencies $c_1,\ldots,c_k\geq 1$ is a string that contains exactly $c_i$ entries~$i$ for all $i\in[k]$.
Permutations are the special case where $c_1=\cdots=c_k=1$.
Multiset permutations also generalize combinations, which are obtained for~$k=2$.

Ko and Ruskey~\cite{MR1161985} describe a genlex Gray code for multiset permutations that uses transpositions, and that is constructed from a generation tree.
Takaoka~\cite{MR1793533} used a similar idea to construct another transposition-based genlex Gray code, which moreover yields a loopless algorithm.
Ruskey~\cite{ruskey_2003} also devised a Gray code that operates by homogeneous transpositions, i.e., whenever entries~$a_i$ and~$a_j$ are transposed, all entries in between satisfy $a_{i+1},\ldots,a_{j-1}\leq \min\{a_i,a_j\}$.
This construction is a generalization of the Eades-McKay homogeneous transposition Gray code for combinations; see also~\cite{MR2022587,DBLP:journals/cj/KorshL04}.
Already much earlier, Chase~\cite{DBLP:journals/cacm/Chase70a} generalized his Gray code for combinations to multiset permutations, such that it only uses transpositions of the form $a_i,a_{i+1}\leftrightarrow a_{i+1},a_i$ or $a_{i-1},a_i,a_{i+1}\leftrightarrow a_{i+1},a_i,a_{i-1}$ with $a_i=\min\{a_{i-1},a_{i+1}\}$.

Lehmer~\cite{MR172816} raised the question under which conditions multiset permutations can be generated by adjacent transpositions.
He phrased the problem in terms of a motel, which hosts $k$ families with $c_i$ members for each~$i\in[k]$, and the task is to create all possible assignments of these guests to $n:=\sum_{i\in[k]}c_i$ rooms, where changes take place only between adjacent rooms.
Stachowiak~\cite{MR1157583} answered this question completely, proving that an adjacent transposition Gray code exists if and only if at least two of the~$c_i$ are odd, or in the trivial cases $k=1$, or $k=2$ and $\min\{c_1,c_2\}=1$; see also~\cite{MR1142265}.

Williams~\cite{MR2807540} showed that multiset permutations can be generated by prefix shifts, and the average length of the shifts is less than~$3$ when $c_1=\cdots=c_k=1$.
He obtained this result by generalizing his cool-lex order to multiset permutations.
In his thesis~\cite{williams_2009} (results also published in~\cite{MR4594331}), he also shows that prefix shifts of length~$n-1$ or~$n$ are sufficient to generate all multiset permutations, which generalizes the aforementioned result for permutations published in~\cite{MR2682614}.

\openproblem{38}{Shen and Williams~{\cite{shen_williams_2021}} conjectured that multiset permutations with equal frequencies $c_1=\cdots=c_k=:c$ can be generated by star transpositions.}
This generalizes Ehrlich's star transposition Gray code for permutations ($c=1$), and the middle levels conjecture ($k=2$).
Such star transposition Gray codes were constructed by Gregor, Merino and M\"utze~\cite{MR4580020} for all frequencies satisfying $c_i\leq 4$ for $i\in[k]$ and $2\cdot\max_{i\in[k]}c_i\leq \sum_{i\in[k]}c_i$, settling in particular the cases $c=2,3,4$ of Shen and Williams' conjecture.
All other cases are still open.

\subsection{Other generalizations and variants}

A \emph{signed permutation} is a permutation in which every entry has a sign~$\{+,-\}$.
Korsh, LaFollette and Lipschutz~\cite{MR2841338} combined the BRGC and the SJT algorithm to a Gray code for signed permutations that operates by adjacent transpositions or inverting the sign of the last entry.

In a complemented prefix reversal we reverse a prefix of a signed permutation and complement all of its signs.
In terms of the pancake flipping analogy, we can imagine pancakes with one side burnt and the other side unburnt, so flipping them around brings up the other side.
Suzuki, N.~Sawada and Kaneko~\cite{DBLP:conf/ISCApdcs/SuzukiSK05} showed that the flip graph on signed permutations under those operations is Hamilton-connected.
Moreover, J.~Sawada and Williams~\cite{MR3513761} describe two simple algorithms for computing a Hamilton cycle in this graph.
The first algorithm greedily applies complemented prefix reversals to the \emph{shortest} possible prefix, so that a new signed permutation is generated.
The average length of the reversed prefixes approaches~$\sqrt{e}=1.64872\ldots$.
Interestingly, the Ord-Smith/Zaks ordering of (unsigned) permutations is also obtained by such a greedy rule.
Moreover, it is shown that greedily reversing (and complementing) the \emph{longest} possible prefix also yields a Gray code listing of permutations (and signed permutations, respectively).
\openproblem{39}{In the follow-up paper~{\cite{MR3426940}}, Sawada and Williams conjectured that signed permutations can be generated using only complemented prefix reversals of lengths~$n$, $n-1$ or $n-2$.}

Cameron, Sawada and Williams~\cite{MR4346446} generalized these results to \emph{colored permutations}, where each entry of the permutation has one of $k$ available colors.
For $k=1$ these are (uncolored) permutations and for $k=2$ these are signed permutations.
This paper establishes that the greedy strategy that modifies the shortest possible prefix by reversing and cyclically shifting colors works for any number~$k$ of colors, and the average length of the reversed prefixes approaches~$e^{1/k}$.

Gutierres, Mamede and Santos~\cite{MR3828769} present a Gray code for signed involutions that uses a transposition, a 3-rotation and/or a sign change in each step.

A \emph{weak order} is a ranking of $n$ objects where ties are allowed.
For example, all thirteen weak orders on three elements~$\{1,2,3\}$ are $1\!\equiv\!2\!\equiv\!3$, $1\!\equiv\!2\!\prec\!3$, $1\!\prec\!2\!\equiv\!3$, $1\!\prec\!2\!\prec\!3$, $1\!\equiv\!3\!\prec\!2$, $1\!\prec\!3\!\prec\!2$, $2\!\prec\!1\!\equiv\!3$, $2\!\prec\!1\!\prec\!3$, $2\!\equiv\!3\!\prec\!1$, $2\!\prec\!3\!\prec\!1$, $3\!\prec\!1\!\equiv\!2$, $3\!\prec\!1\!\prec\!2$, $3\!\prec\!2\!\prec\!1$, where $\prec$ denotes a strict precedence and $\equiv$ denotes a tie.
The number of weak orders is given by the Fubini numbers~[OEIS~000670].
Any weak order can be described by a string $a_1,\ldots,a_n$, where $a_i$ denotes the number of~$\prec$ symbols that precede~$i$ in the weak order (these strings are called \emph{Cayley permutations} in~\cite{MR732206}).
The corresponding strings for the weak orders for $n=3$ listed before are $000,001,011,012,010,021,101,102,100,201,110,120,210$.
Baril~\cite{MR2027523} and Jacques, Wong and Woo~\cite{MR4115514} describe 2-Gray codes for weak orders in this string representation, i.e., any two consecutive strings differ in at most two positions.
The following two challenging open problems are mentioned in~\cite[Problem~482]{MR2548549}:
\openproblem{40}{Savage raised the question whether there is a 1-Gray code for weak orders.}
For $n=3$ such a (non-cyclic) Gray code is $000, 001, 201, 101, 102, 100, 120, 110, 210, 010, 012, 011, 021$, and computer experiments show that they exist also for $n=4,5,6$ (and these are even cyclic).
\openproblem{41}{Moreover, Knuth asked for a 2-Gray code that uses only adjacent transpositions $a_i\leftrightarrow a_{i+1}$ and adjacent left-assignments~$a_i:=a_{i+1}$; see also~\mbox{\cite[Ex.~106 in Sec.~7.2.1.2]{MR3444818}}.}

Cardinal, Merino and M\"utze~\cite{MR4415121} consider \emph{partial permutations}, which are all linearly ordered subsets of~$[n]$ of arbitrary size [OEIS~A000522].
Their Gray code uses either an adjacent transposition, deletion or insertion of a trailing element in each step.

\emph{Stamp foldings}, \emph{semi-meanders} and \emph{open meanders} are permutations that are characterized by avoiding certain geometric patterns (self intersections), which can be captured by permutation patterns with additional parity constraints.
\openproblem{42}{Sawada and Li~{\cite{MR2946101}} raised the question which of these families of objects admit Gray codes.}
Liu and Wong~\cite{MR4607055} made progress towards this problem, by providing Gray codes for stamp foldings and semi-meanders; see Figure~\ref{fig:perm3}~(n)+(n').
Their code operates by cyclic substring shifts (by more than one position) in the permutation representation.

\section{Set partitions}

\begin{figure}
\centerline{
\setlength{\tabcolsep}{3pt}
\renewcommand{\arraystretch}{1.5}
\begin{tabular}{ccccccccccccccc}
(a) & (b) & (c) & (d) & (e) & (f) & (g2) & (g3) & (g4) & (h2) & (h3) & (h4) & (i3) & (i4) & (i5) \\[-5mm]
\raisebox{-\height}{\includegraphics{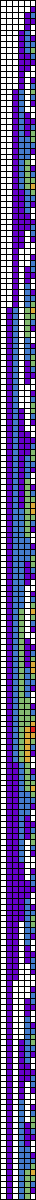}} &
\raisebox{-\height}{\includegraphics{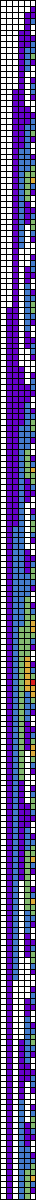}} &
\raisebox{-\height}{\includegraphics{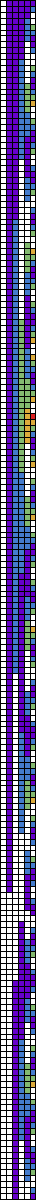}} &
\raisebox{-\height}{\includegraphics{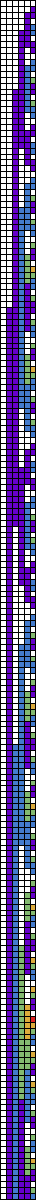}} &
\raisebox{-\height}{\includegraphics{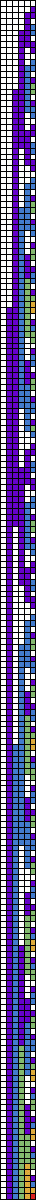}} &
\raisebox{-\height}{\includegraphics{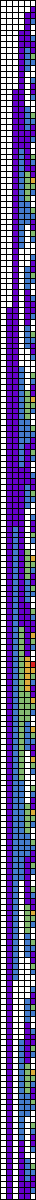}} &
\raisebox{-\height}{\includegraphics{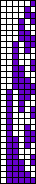}} &
\raisebox{-\height}{\includegraphics{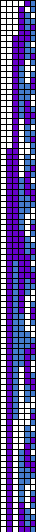}} &
\raisebox{-\height}{\includegraphics{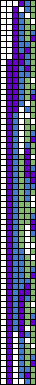}} &
\raisebox{-\height}{\includegraphics{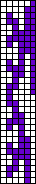}} &
\raisebox{-\height}{\includegraphics{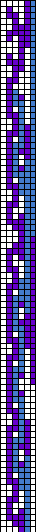}} &
\raisebox{-\height}{\includegraphics{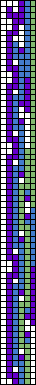}} &
\raisebox{-\height}{\includegraphics{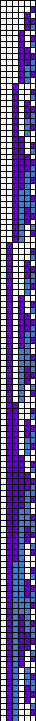}} &
\raisebox{-\height}{\includegraphics{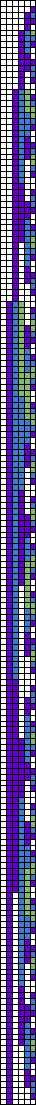}} &
\raisebox{-\height}{\includegraphics{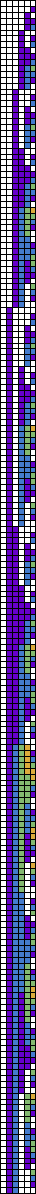}}
\end{tabular}
}
\caption{Gray codes for set partitions of~$[6]$ or~$[5]$ represented by restricted growth strings (0=white, 1=purple, 2=blue, 3=green, 4=yellow, 5=red, 6=brown):
(a)~Kaye;
(b)~Ruskey;
(c)~Ehrlich;
(d) Mansour-Nassar;
(e) greedy (Williams);
(f) greedy (Kerr-Sawada);
%
%
(g2)--(g4)~with 2, 3 or 4 parts (Ehrlich);
(h2)--(h4)~with 2, 3 or 4 parts (Ruskey);
(i3)--(i5)~with $\leq$ 3, 4 or 5 parts (Sabri-Vajnovszki);
(figure continues on next page)
}
\label{fig:setp}
\end{figure}

\begin{figure}
\ContinuedFloat
\centerline{
\setlength{\tabcolsep}{3pt}
\renewcommand{\arraystretch}{1.5}
\begin{tabular}{cccccccc}
(j) & (k) & (l) & (m) & (m') & (n) & (n') & (n'') \\[-5mm]
\raisebox{-\height}{\includegraphics{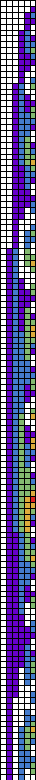}} &
\raisebox{-\height}{\includegraphics{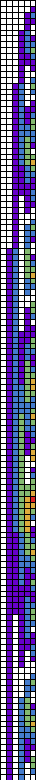}} &
\raisebox{-\height}{\includegraphics{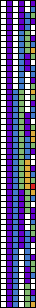}} &
\raisebox{-\height}{\includegraphics{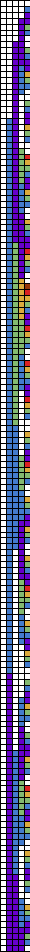}} &
\raisebox{-\height}{\includegraphics{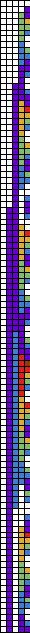}} &
\raisebox{-\height}{\includegraphics{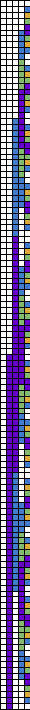}} &
\raisebox{-\height}{\includegraphics{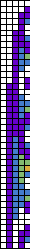}} &
\raisebox{-\height}{\includegraphics{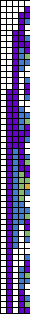}}
\end{tabular}
}
\caption{(continued)
(j)~non-crossing;
(k)~non-nesting;
(l)~sparse;
(m)+(m')~restricted growth functions: $b=(b_2,b_3,b_4,b_5)=(2,1,1,2)$ and $b=(b_2,b_3,b_4,b_5)=(1,1,3,1)$;
(n) subexcedant sequences: $f(a_1,\ldots,a_{i-1})=i-2$;
(n') staircase strings: $f(a_1,\ldots,a_{i-1})=a_{i-1}$;
(n'') ascent sequences: $f(a_1,\ldots,a_{i-1})=|\{1\leq j<i-1\mid a_j<a_{j+1}\}|$.
}
\end{figure}

A \emph{set partition} is a partition of~$[n]$ into a disjoint union of nonempty sets.
We consider every set partition in its canonic representation, where every subset is sorted increasingly, and the subsets are sorted increasingly by their smallest elements.
A \emph{restricted growth string} (RGS) is a string $a_1,\ldots,a_n$ such that $a_1=0$ and $a_i\leq 1+\max\{a_1,\ldots,a_{i-1}\}$ for all $i=2,\ldots,n$.
It is well known that set partitions are in bijection with RGSs.
Indeed, given a set partition in canonic representation, we give labels $j=0,1,\ldots$ to each subset according to the canonic ordering, and we then define $a_i:=j$ if the element~$i$ is in the subset labeled~$j$.
For example, the set partition $\{1,7,8\},\{2,3,5,6\},\{4\},\{9\}$ of $[9]$ corresponds to the RGS~$011211003$.
In this section we consider Gray codes for all set partitions (unrestricted) and some restricted variants obtained by fixing the number of parts or forbidding certain patterns, as well generalizations of RGSs.
Unrestricted set partitions are counted by the Bell numbers [OEIS~A000110], whereas set partitions with a fixed number of parts are counted by the Stirling numbers of the second kind [OEIS~A008277].
Furthermore, non-crossing and non-nesting partitions are counted by the Catalan numbers, and non-crossing and non-nesting partitions with a fixed number of parts are counted by the Narayana numbers [OEIS~A001263].

\subsection{Unrestricted set partitions}
\label{sec:setp-all}

Kaye~\cite{MR443423} describes a Gray code for set partitions where in each step, one element moves from one subset to an adjacent subset (in the canonic representation); see Figure~\ref{fig:setp}~(a).
This Gray code for the partitions of~$[n]$ can be constructed inductively by considering the listing for partitions of~$[n-1]$ and adding~$n$ to every subset and as a singleton, alternatingly from left to right or vice versa.
In Kaye's paper and in~\cite{MR993775}, this construction idea is attributed to Knuth, who solved the corresponding problem stated in~\cite[p.~292, Problem~35]{MR0396274}.
The corresponding RGS listing is genlex, but consecutive RGSs may differ in multiple positions.
Ruskey~\cite{ruskey_2003} proposed a simple modification of this ordering which maintains the genlex property and ensures that consecutive RGSs differ in exactly one position, and by $\pm 1$ or $\pm 2$ in that position; see Figure~\ref{fig:setp}~(b).
Already earlier, Ehrlich~\cite{MR0366085} presented a Gray code for set partitions which is genlex and consecutive RGSs differ in exactly one position, either by $\pm 1$ or the value in that position changes from $0$ to the maximum value or vice versa; see Figure~\ref{fig:setp}~(c).
In the canonic set representation, this means that an element moves to a circularly adjacent subset, i.e., the first and last subset are also considered adjacent.
Ehrlich also proved that under the stricter requirement that two consecutive RGSs differ only by $\pm 1$, Gray codes do not exist for infinitely many values of~$n$ because of balancedness problems.
Mansour and Nassar~\cite{MR2443226} present another genlex Gray code for set partitions in which consecutive RGSs differ in exactly one position.

Williams~\cite{MR3126386} describes a genlex Gray code for set partitions that is produced by starting at the partition $\{1,2,\ldots,n\}$, and then repeatedly and greedily moving the largest possible element to the leftmost possible set so that a new set partition is created; see Figure~\ref{fig:setp}~(e).
In the corresponding RGS listing, some strings differ in more than one position (in the listing for $n=6$ shown in the figure, there is precisely one such transition from the 73rd to the 74th string, which differ in the last two positions).
Williams' paper claims that this listing is the same as Kaye's, but this is not the case.
In particular, Williams' listing always ends with the partition $\{1\},\{2\},\ldots,\{n\}$, unlike Kaye's.
Yet another genlex greedy Gray code was communicated to the author by Joe Sawada and his student Kevin Kerr in~2021~\cite{kerr_2021}.
It starts with~$0^n$ and then repeatedly and greedily changes the rightmost possible entry (and only that entry) to the largest possible value so that a new RGS is created; see Figure~\ref{fig:setp}~(f).

\subsection{Fixed number of parts}

Ehrlich's paper~\cite{MR0366085} also describes a Gray code for set partitions with a fixed number of parts; see Figure~\ref{fig:setp}~(g2)--(g4).
The sets differ in moving one or two elements to another subset, and the corresponding RGS listings are genlex and differ in at most two positions in each step.
His construction works more generally for set partitions whose number of parts is in any fixed interval.
Ruskey~\cite{MR1288775} presents another Gray code for this problem in which the RGSs differ in exactly one position in each step, but this listing is not genlex; see Figure~\ref{fig:setp}~(h2)--(h4).
Ruskey and Savage~\cite{MR1296942} showed that a Gray code with only single $\pm 1$ transitions in the RGS representation does not exist in general.

Sabri and Vajnovszki~\cite{MR3654140} consider Gray codes for set partitions with a fixed number of parts or an upper bound on the number of parts.
They constructed a 5-Gray (in the RGS representation) for the first family and a 3-Gray code for the second family, based on sublists of generalizations of the BRGC, implying that both of these Gray codes are genlex; see Figure~\ref{fig:setp}~(i3)--(i5).

\subsection{Non-crossing, non-nesting and sparse partitions}
\label{sec:setp-geom}

Huemer, Hurtado, Noy and Oma\~na-Pulido~\cite{MR2510231} developed Gray codes for non-crossing set partitions.
A set partition is \emph{non-crossing} if for any two sets~$A$ and~$B$ of the partition and any four elements~$a,a'\in A$ and $b,b'\in B$ we do not have $a<b<a'<b'$.
The authors construct a cyclic Gray code for non-crossing partitions, valid for all~$n\geq 3$, that moves one element of a subset to another in each step.
A non-cyclic Gray code for this problem can be obtained easily, by following Kaye's construction described in Section~\ref{sec:setp-all} and skipping set partitions that have crossings, which yields a genlex listing; see Figure~\ref{fig:setp}~(j).
This paper also construct a cyclic Gray code for non-crossing partitions of fixed size~$k$, for all~$n>k\geq 3$, allowing however one additional flip operation, namely to exchange two elements between two subsets, i.e., $A\cup\{a\},B\cup\{b\}\longrightarrow A\cup\{b\},B\cup\{a\}$.
\openproblem{43}{The authors conjectured that the latter problem can be solved for all $n>k\geq 3$ without using this extra type of flip operation.}

Complementing the aforementioned results, Conflitti and Mamede~\cite{MR3367583} consider \emph{non-nesting} set partitions, which are set partitions such that for any two sets~$A$ and~$B$ of the partition and any four elements~$a,a'\in A$ and $b,b'\in B$ with $[b,b']\cap A=\emptyset$ we do not have $a<b<b'<a'$.
They present a Gray code for non-nesting partitions that moves one element of a subset to another in each step; see Figure~\ref{fig:setp}~(k).
The construction is a straightforward variation of Kaye's construction, and it also yields a genlex listing.
The authors used the same idea to construct a genlex Gray code for \emph{sparse} set partitions, which are partitions such that no subset contains two consecutive elements; see Figure~\ref{fig:setp}~(l).
Interestingly, there is a simple bijection between sparse set partitions of~$[n]$ and unrestricted set partitions of~$[n-1]$.
In particular, sparse set partitions are also counted by the Bell numbers.
In addition to these simple non-cyclic constructions, Conflitti and Mamede also provide constructions of cyclic Gray codes for both problems.
In a later paper, Conflitti and Mamede~\cite{MR3690404} present similar Gray codes for non-crossing and non-nesting partitions in Weyl groups of type~B and~D.

\openproblem{44}{For all the classes of restricted set partitions discussed in this section, it is open whether they admit Gray codes that change only few entries in the RGS representation in each step.}

\subsection{Generalizations}

For any integer $k\geq 0$, a \emph{restricted growth tail} is a string $a_1,\ldots,a_n$ such that $a_1\leq k$ and $a_i\leq 1+\max\{a_1,\ldots,a_{i-1},k-1\}$ for all $i=2,\ldots,n$.
Obviously, for $k=0$ these are the same as RGSs.
Ruskey and Savage~\cite{MR1296942} present Gray codes for restricted growth tails, valid for any $n\geq 1$ and $k\geq 0$, which satisfy Ehrlich's conditions mentioned in Section~\ref{sec:setp-all}, i.e., a single position changes in each step, either by $\pm 1$ or the value in that position changes from $0$ to the maximum value of at least~$k$ or vice versa.
These Gray codes are moreover cyclic.

Mansour, Nassar and Vajnovszki~\cite{MR2816662} consider another generalization of RGSs, called \emph{restricted growth functions}, for which there is a vector $b=(b_2,\ldots,b_n)$ with $b_i\geq 1$ for all $i=2,\ldots,n$, and the requirement for the string $a_1,\ldots,a_n$ is that $a_1=0$ and $a_i\leq b_i+\max\{a_1,\ldots,a_{i-1}\}$ for all $i=2,\ldots,n$.
Clearly, for $b=(1,\ldots,1)$ we recover standard RGSs.
They present genlex Gray codes that are valid for any vector~$b$, and those Gray codes change only one entry of the RGS representation at a time; see Figure~\ref{fig:setp}~(m)+(m').
For $b=(1,\ldots,1)$ the obtained listings coincide with the Gray code of Kerr-Sawada mentioned before, and they can also be constructed greedily, by starting with~$0^n$ and then repeatedly and greedily changing the rightmost possible entry to the largest possible symbol.

Subsequently, Mansour and Vajnovszki~\cite{MR3070506} introduced an even more general concept, called \emph{restricted growth words}.
They consider strings $a_1,\ldots,a_n$ such that $a_1=0$ and $a_i\leq 1+f(a_1,\ldots,a_{i-1})$, where $f$ is any non-negative function of $a_1,\ldots,a_{i-1}$.
Their paper describes a genlex Gray code for any such family of strings, which changes only one entry at a time, either by~$\pm 1$ or~$\pm 2$.
For $f=\max$ we obtain RGSs.
For $f=i-2$ we obtain \emph{subexcedant sequences}, which are counted by the factorial numbers; see Figure~\ref{fig:setp}~(n).
For $f=a_{i-1}$ we obtain \emph{staircase words}, which are another Catalan family; see Figure~\ref{fig:setp}~(n').
For $f=0$ we obtain bitstrings with first bit~0, which are counted by powers of~2, and the resulting listing is the (first half of the) BRGC.
For $f=|\{1\leq j<i-1\mid a_j<a_{j+1}\}|$ we obtain \emph{ascent sequences}, which are counted by the Fishburn numbers~[OEIS~A022493]; see Figure~\ref{fig:setp}~(n'').
Interestingly, the Kerr-Sawada greedy approach also works in this more general setting, yielding different listings, however.

Sabri and Vajnovszki~\cite{DBLP:journals/cj/SabriV15} showed that restricting generalizations of the BRGC to RGSs, subexcedant sequences, staircase words or ascent sequences yields genlex 3-Gray codes.

\section{Integer partitions and compositions}

In this section we discuss Gray codes for integer partitions and integer compositions.

\subsection{Integer partitions}

A \emph{partition} of~$n$ is a sequence $a_1,\ldots,a_k$ of integers such that $a_1\geq a_2\geq \cdots\geq a_k\geq 1$ and $a_1+\cdots+a_k=n$.
Integer partitions are counted by the partition numbers~[OEIS~A000041].

Wilf asked whether integer partitions of~$n$ admit a Gray code with transitions of the form $a_i,a_j\leftrightarrow a_i-1,a_j+1$ in each step, which includes transitions $a_i(,0)\leftrightarrow a_i-1,1$ that change the number of parts.
This problem was solved by Savage~\cite{MR1022113}, who showed more generally that for all $n\geq k\geq 1$, all partitions of~$n$ with parts of size at most~$k$ admit such a Gray code.
Her results also yield Gray codes for all partitions of~$n$ with largest part equal to~$k$, all partitions of~$n$ into at most $k$ parts, and all partitions of~$n$ into exactly~$k$ parts.

Rasmussen, Savage and West~\cite{MR1329389} constructed a Gray code for partitions of~$n$ into parts of size at most~$k$ and congruent to~1 modulo~$d$, where $d\geq 1$ is some integer parameter.
For $d=1$ this is the same problem as before.
For $d=2$ these are partitions into odd parts of size at most~$k$.
Their Gray code uses transitions $a_i,a_j\leftrightarrow a_i-d,a_j+d$ and $a_i(,0^d)\leftrightarrow a_i-d,1^d$.
Furthermore, the authors present a Gray code for partitions of~$n$ into distinct parts of size at most~$k$, using transitions $a_i,a_j\leftrightarrow a_i-1,a_j+1$ as before.
It is a well-known result of Euler that the number of partitions of~$n$ into odd parts is the same as the number of partitions of~$n$ into distinct parts~[OEIS~A000009].

\openproblem{45}{The paper~{\cite{MR1329389}} also raised the problem to construct Gray codes for the following five families of restricted partitions of~$n$: exactly $k$ distinct parts; distinct odd parts; distinct parts congruent to~1 modulo~$d$; at most $c$ copies of each part; parts congruent to~1 modulo~$d$ and at most $c$ copies of each part.}
Savage~\cite{MR1491049} conjectured that such codes exist for all $n\geq k\geq 1$ and $c,d\geq 1$.

A compact representation of an integer partition~$a_1,\ldots,a_k$ is obtained by grouping together runs of~$a_i$ that have the same value, and recording their multiplicities, i.e., we obtain a list of pairs $(a_1',b_1),\ldots,(a_\ell',b_\ell)$ such that $a_1'>\cdots>a_\ell'\geq 1$, $b_i\geq 1$ for $i=1,\ldots,\ell$ and $a_1'\cdot b_1+\cdots+a_\ell'\cdot b_\ell=n$.
Ruskey~\cite{ruskey_2003} observed that ordering partitions lexicographically, the corresponding compact representations differ at most in the last three pairs of integers.

\subsection{Integer compositions}
\label{sec:comp}

A \emph{composition} of~$n$ is a sequence of integers $a_1,\ldots,a_k$ such that $a_1,\ldots,a_k\geq 1$ and $a_1+\cdots+a_k=n$, i.e., in contrast to partitions, the order of the summands matters.
The number of compositions of~$n$ is $2^{n-1}$, and there is a simple bijection between bitstrings of length~$n-1$ and compositions:
Specifically, for any bitstring~$x$, we consider the lengths of all maximal non-overlapping substrings of~$x1$ that end with~1.
For example, for the string~$x=1101000110000$ of length~13 we have $x1=1|1|01|0001|1|00001$, where the vertical bars indicate the boundaries between maximal non-overlapping substrings ending with~1, and therefore $x$ is mapped to the composition~$1,1,2,4,1,5$ of~14.
Observe that the BRGC has the interesting property that for every flip~$1\leftrightarrow 0$, the flipped bit is the last one or at least one of the bits next to it is also a~1.
Applying the aforementioned bijection to the BRGC, we thus obtain a cyclic genlex Gray code for compositions that uses only transitions $a_i,1\leftrightarrow a_i+1$, where the~1 is at positions~$i+1$ or~$i-1$, respectively.

If we consider compositions of~$n$ with a fixed number~$k$ of parts, then the above bijection shows that they correspond to $(n-1,k-1)$-combinations.
Using the Gray codes for combinations discussed in Section~\ref{sec:comb-trans} that use transpositions $01\leftrightarrow 10$ or $001\leftrightarrow 100$, we obtain Gray codes for compositions with transitions $a_i,a_{i+1}\leftrightarrow a_i-d,a_{i+1}+d$ with $d=\pm 1$ or $d=\pm 2$, respectively (recall that using only the former requires some conditions on~$n$ and~$k$).

Klingsberg~\cite{MR646890} describes a genlex Gray code for compositions into a fixed number of parts that uses only transitions $a_i,a_j\leftrightarrow a_i-1,a_j+1$ (cf. also~\cite{MR2738477}).
In Klingsberg's paper and in~\cite{MR993775}, the idea for this construction is attributed to Knuth, who solved the corresponding problem stated in~\cite[p.~292, Problem~34]{MR0396274}.

Ruskey and Savage~\cite{MR1397157} consider compositions~$a_1,\ldots,a_k$ of~$n$ with a fixed number~$k$ of parts, with the additional requirement that the parts are bounded, i.e., there is a vector $b=(b_1,\ldots,b_k)$ such that $a_i\leq b_i$ for all $i=1,\ldots,k$.
For $b_i\geq n$ this restriction becomes irrelevant, so we obtain compositions of~$n$ into~$k$ parts.
Moreover, for $b_1=\cdots=b_k=2$ these are $(k,n-k)$-combinations (by changing the alphabet $1,2\rightarrow 0,1$).
The authors construct a Gray code that uses transitions $a_i,a_j\leftrightarrow a_i-1,a_j+1$.
\openproblem{46}{They raise the problem to determine when an adjacent transition Gray code for this problem is possible, i.e., if we also require that $|j-i|=1$ (in some cases, there are balancedness problems).}
Earlier, Ehrlich~\cite{DBLP:journals/cacm/Ehrlich73,MR0366085} presented an adjacent transition Gray code, but with changes that can be more than~$\pm 1$.

Walsh~\cite{MR1772772} considered the same setting of bounded compositions as Ruskey and Savage, and he showed that for any vector~$b$, the sublist of Klingsberg's ordering obtained by selecting only the compositions that satisfy the upper bounds given by~$b$ forms a genlex Gray code with transitions $a_i,a_j\leftrightarrow a_i-1,a_j+1$.
By choosing an appropriate~$b$ and interpreting every composition as an inversion vector of a permutation, he used this to derive a Gray code for permutations with a fixed number of inversions; see Section~\ref{sec:perm-restr}.

Vajnovszki and Vernay~\cite{MR2817083} showed more generally that we may start with a straightforward generalization of the BRGC for tuples $a_1,\ldots,a_k$ with $1\leq a_i\leq b_i$ (recall Section~\ref{sec:larger}), and consider the sublist of all tuples that satisfy $s\leq a_1+\cdots+a_k\leq t$, and this sublist will be a genlex Gray code with transitions $a_i\leftrightarrow a_i\pm 1$ or $a_i,a_j\leftrightarrow a_i-1,a_j+1$.
In particular, Walsh's result for compositions follows from this by setting $s=t=n$.
Furthermore, this also generalizes Tang and Liu's result about combinations discussed in Section~\ref{sec:comb-trans} (for $b_1=\cdots=b_k=2$ and $s=t$), who obtained their Gray code for combinations by restricting the ordinary BRGC to bitstrings with fixed weight.

\section{Geometric objects}

In this section we consider Gray codes for geometric objects equipped with flip operations that are constrained by geometry, e.g., from a planarity requirement, and this makes the topic particularly interesting.
For results on flips in planar graphs beyond Gray codes, we refer to the survey of Bose and Hurtado~\cite{MR2455502}.
In addition to the material in this section, non-crossing and non-nesting set partitions are discussed in Section~\ref{sec:setp-geom}.

\subsection{Triangulations and matchings}
\label{sec:triang}

We first consider the flip graph~$G_n$ that has as vertices all triangulations of a convex $n$-gon, and edges connecting triangulations that differ by a \emph{flip}, i.e., a diagonal in a quadrilateral formed by two triangles is replaced by the other diagonal.
The number of such triangulations is the $(n-2)$nd Catalan number.
It is well known that $G_n$ is the cover graph of a lattice, namely the Tamari lattice, and the skeleton of a polytope, namely the $(n-3)$-dimensional associahedron.
By considering the dual graphs of the triangulations, we see that triangulations are isomorphic to binary trees with $n-2$ nodes, and edge flips correspond to tree rotations; see Figure~\ref{fig:flip}~(f).
Lucas~\cite{MR920505} proved that $G_n$ has a Hamilton cycle, and her proof yields an efficient generation algorithm.
Subsequently, Lucas, Roelants van Baronaigien and Ruskey~\cite{MR1239499} provided a simpler proof for the existence of a Hamilton path in~$G_n$, which yields a loopless algorithm.
Possibly the simplest description of their path was given by Williams~\cite{MR3126386} via a greedy algorithm, which is to repeatedly rotate the edge with the largest inorder label so that a new tree is generated; see also~\cite{MR4415121} and Figure~\ref{fig:231cat}.
Another proof for the existence of a Hamilton cycle in~$G_n$ was obtained by Hurtado and Noy~\cite{MR1723053} via a generation tree.
Huemer, Hurtado and Pfeifle~\cite{MR2474724} later generalized this result to a flip graph of dissections of a convex polygon into $(b+1)$-gons under edge flips, whose dual graphs are $b$-ary trees under tree rotations.

Hernando, Hurtado and Noy~\cite{MR1939072} consider the set of all non-crossing perfect matchings on a set of $n=2m$ points in convex position, with the flip operation that replaces edges and non-edges along an alternating 4-cycle.
They proved that the corresponding flip graph has a Hamilton cycle for every even~$m\geq 4$, and no Hamilton path for any odd~$m\geq 5$.
Note that non-crossing matchings are also counted by the Catalan numbers.

Aichholzer, Aurenhammer, Huemer and Vogtenhuber~\cite{MR2346418} exhibit a set of six points in non-convex position for which the flip graph of triangulations is a star with three rays, so it does not admit a Hamilton path.
\openproblem{47}{They speculate whether Hamiltonicity holds for triangulations of every sufficiently large point set.
Moreover, they conjectured that the flip graph of pseudo-triangulations under edge removals/insertions and exchanges has a Hamilton cycle, where a pseudo-triangulation has faces that are bounded by three concave chains.}

Certain flip graphs for triangulations and non-crossing matchings on general point sets were considered by~Houle, Hurtado, Noy and Rivera-Campo~\cite{MR2190792}, and they established that these graphs are connected, which implies 3-Gray codes by Sekanina's~\cite{MR0140095} theorem discussed in Section~\ref{sec:general}.

\openproblem{66}{Gregor, Merino and~M\"utze~\cite{MR4747482} asked for the largest integer~$k$ for which there is a $k$-symmetric Hamilton cycle as defined in Section~\ref{sec:symm} in the associahedron.}
As a stepping stone, are there Gray codes for necklaces and bracelets of triangulations on a convex point set, i.e., for equivalence classes of these objects under cyclic rotation and/or reflection?
\openproblem{67}{The same two questions can naturally be asked for the flip graph of non-crossing matchings on a convex point set.}

\subsection{Ordered trees}
\label{sec:otrees}

An \emph{ordered} tree with $n$ nodes has a root and a specified left-to-right ordering of the children of each node.
These objects are counted by the $(n-1)$st Catalan number.
Lapey and Williams~\cite{MR4526706} provide a loopless Gray code algorithm for ordered trees that uses so-called pull operations, at most 2 in each step, and their method is derived from cool-lex order.
Nakano~\cite{DBLP:conf/tamc/Nakano24} gives a Gray code for ordered trees that uses single leaf moves, i.e., in each step a leaf of the current tree is removed and reattached elsewhere.

\subsection{Dissections, plane graphs and spanning trees}

Huemer, Hurtado, Noy and Oma\~na-Pulido~\cite{MR2510231} consider Gray codes for dissections of a convex polygon.
A \emph{dissection} of a convex polygon with $n$ vertices is a subdivision of the polygon by $k$ non-crossing diagonals, where $0\leq k\leq n-3$.
We let $G_n$ and $G_{n,k}$ be the graphs on the set of all dissections and those with exactly $k$ diagonals, respectively.
In $G_{n,k}$ two dissections are connected by an edge if they differ by removing one diagonal and adding another diagonal in the subpolygon that contains the end vertices of the removed diagonal.
The graph~$G_n$ is the union over all~$G_{n,k}$ for $k=0,\ldots,n-3$ plus edges that correspond to either removing or adding one diagonal in the dissection.
The authors prove that $G_n$ has a Hamilton cycle for all $n\geq 4$, and $G_{n,k}$ has a Hamilton cycle for all $k\geq 1$ and $n\geq 4$.
Note that the vertices of~$G_{n,k}$ are precisely the $(n-3-k)$-dimensional faces of the associahedron (as a polytope), so this result generalizes the aforementioned result about the Hamiltonicity of the associahedron (the special case $k=n-3$).

Hernando, Hurtado, Marquez, Mora and Noy~\cite{MR1684903} proved that all plane spanning trees on a set of $n\geq 3$ points in convex position can be listed cyclically so that consecutive trees differ by exchanging one edge.
Rivera-Campo and Urrutia-Galicia~\cite{MR1826999} proved the same result for plane spanning paths.

Aichholzer, Aurenhammer, Huemer and Vogtenhuber~\cite{MR2346418} consider point sets in general position, with plane graphs whose edges connect certain pairs of points by a straight line.
In this setting, they obtained Gray codes via generation trees for three different families of graphs: all plane graphs by edge removals/insertions, all plane spanning trees by edge exchanges (generalizing the result from~\cite{MR1684903} for convex position), and all connected plane graphs by edge removals/insertions or exchanges.

\subsection{Rectangulations and rainbow cycles}
\label{sec:rect}

\begin{figure}
\includegraphics{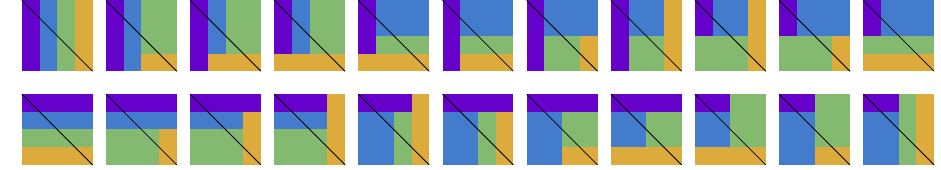}
\caption{Gray code for diagonal rectangulations with $n=4$ rectangles.}
\label{fig:rect}
\end{figure}

Merino and M\"utze~\cite{MR4598046} present Gray codes for various families of \emph{generic rectangulations}, which are subdivisions of a rectangle into $n$ smaller rectangles, subject to the constraint that no four rectangles meet in a point.
In particular, they constructed a Gray code for all generic rectangulations with $n$ rectangles by flip operations called wall slides, simple flips and T-flips.
They also obtained a Gray code for diagonal rectangulations, which are those generic rectangulations for which every rectangle intersects the diagonal from the top-left to the bottom-right corner of the outer rectangle, and this Gray code uses only simple and T-flips (diagonal rectangulations are counted by the Baxter numbers~[OEIS~A001181]); see Figure~\ref{fig:rect}.
More generally, their approach yields Gray codes for many subfamilies of generic rectangulations that are characterized by pattern avoidance, such as guillotine rectangulations or 1-sided rectangulations, by so-called jumps, which generalize wall slides, simple flips and T-flips.
\openproblem{48}{They raised the problem of constructing Gray codes for rectangulations of point sets as introduced in~{\cite{MR2244135}} (see also~{\cite{DBLP:journals/ieicet/YamanakaRN18}}).}
Downing, Einstein, Hartung and Williams~\cite{downing_einstein_hartung_williams_2023} provide Gray codes for Catalan squares and staircases, which are a Catalan family of geometric objects that can be viewed as T-avoiding rectangulations.

Felsner, Kleist, M\"utze and Sering~\cite{MR4046775} constructed cycles in the flip graphs of different geometric configurations that are balanced in the sense that each flip operation occurs equally often along the cycle.
They refer to these cycles as \emph{rainbow cycles}, and establish existence and non-existence results for the flip graphs of triangulations of a convex polygon, non-crossing perfect matchings on a convex set of points, and plane spanning trees on a general point set.
\openproblem{49}{They conjectured that the flip graph of non-crossing perfect matchings on a convex $2m$-gon under flips along alternating 4-cycles does not contain a cycle along which every possible edge on the point set appears and disappears exactly once (such a cycle has length~$m^2/2$), whenever $m\geq 6$ is even, which is confirmed for $m\in\{6,8,10,12\}$ (see also~{\cite{MR4184338}}).}

\section{Objects on graphs}
\label{sec:graphs}

In this section we discuss Gray codes for combinatorial objects that are defined with respect to a fixed graph~$H$.
We write $K_n$, $P_n$ and~$C_n$ for the complete graph, path or cycle on $n$ vertices, respectively, and $K_{a_1,a_2,\ldots,a_r}$ for the complete $r$-partite graph with partition classes of sizes~$a_1,\ldots,a_r$.

\subsection{Spanning trees, cycles, and generalizations to matroids}
\label{sec:span}

Cummins~\cite{MR245357} showed that the spanning trees of any fixed graph~$H$ can be listed cyclically by edge exchanges, i.e., in each step one edge is removed from the current spanning tree and another edge is inserted; see Figure~\ref{fig:flip}~(g).
Equivalently, his result asserts that the flip graph of the spanning trees of~$H$ under edge exchanges has a Hamilton cycle.
Cummins proved more strongly that this flip graph is \emph{edge-Hamiltonian}, i.e., for each of its edges, there is a Hamilton cycle containing this edge.
A short proof of this result was given by Shank~\cite{1082765}.
Smith~\cite{smith_1997} devised an algorithm for listing all spanning trees of a graph by edge exchanges, which requires only constant time on average per generated spanning tree.
Details of his algorithm are described in~\cite[Sec.~7.2.1.6]{MR3444818}.

Holzmann and Harary~\cite{MR307952} generalized the aforementioned results to matroids, proving that the flip graph of the bases of a matroid under element exchanges has a Hamilton cycle.
They showed more generally that if the flip graph has at least two cycles, then it is \emph{Hamilton-uniform}, i.e., for every edge of the flip graph, there is a Hamilton cycle containing this edge, and a Hamilton cycle avoiding this edge.
Furthermore, by the results from~\cite{MR638286} discussed in Section~\ref{sec:01}, the matroid base flip graph is either Hamilton-connected or a hypercube, which is known to be Hamilton-laceable.
Alspach and Liu~\cite{MR1027700} showed that for simple matroids, this flip graph is also \emph{edge-pancyclic}, i.e., for every edge and every integer~$\ell$ at least~3 and at most the number of vertices of the flip graph, there is a cycle of length~$\ell$ through this edge; see also~\cite{MR401539}.
Note that for the uniform matroid, the matroid base flip graph is the Johnson graph, and in this case a Hamilton cycle corresponds to a transposition Gray code for combinations.

Many of the aforementioned Gray codes for spanning trees of a graph~$H$ are based on a straightforward recursion, already observed by Feussner~\cite{feussner_1902} (see also~\cite{1082385}), namely to consider the spanning trees that do not contain a fixed edge~$e$, which are the spanning trees of~$H-e$, and those that do contain~$e$, which are the spanning trees of~$H/e$.
The graph~$H-e$ is obtained from~$H$ by removing~$e$, but none of its end vertices, whereas $H/e$ is obtained from~$H$ by contracting the end vertices of~$e$ to a single vertex.
Based on this recursion, Merino, M\"utze and Williams~\cite{DBLP:conf/fun/MerinoMW22} discovered a strikingly simple greedy algorithm to list the spanning trees of~$H$ by edge exchanges, by considering their characteristic vectors:
An edge exchange corresponds to a transposition of a~0 and~1 in the characteristic vector, and the rule is to perform a transposition that changes the shortest possible prefix so as to generate a previously unseen spanning tree.
This amazing rule works for any graph~$H$, for any edge ordering that determines the indices of the characteristic vector, for any initial spanning tree, and for any tie-breaking rule, it always yields a genlex listing, and it also generalizes straightforwardly to bases of a matroid (however, this Gray code is in general not cyclic).

\openproblem{50}{Cameron, Grubb and Sawada~{\cite{DBLP:conf/cocoon/CameronGS21}} asked whether there is a Gray code for spanning trees of a graph~$H$ that uses edge exchanges, with the additional requirement that the removed and added edge have a vertex of~$H$ in common.}
They refer to a listing satisfying this stronger requirement as a \emph{pivot Gray code}, and they constructed such a Gray code for the family of fan graphs, which are obtained by joining a single vertex to all vertices of a path.
\openproblem{51}{Knuth~\mbox{\cite[Ex.~102 in Sec.~7.2.1.6]{MR3444818}} considers a directed graph and all of its spanning trees in which all arcs are oriented away from a fixed root vertex, also known as \emph{out-trees} or \emph{arborescences}.
He asked for a Gray code to list those directed trees in Gray code order, exchanging one arc at a time.}
Note that the removed and added arc must point to the same vertex.
We conclude that the existence of Knuth's Gray code would imply a pivot Gray code for undirected graphs (simply replace each undirected edge with two arcs in opposite directions).

\openproblem{52}{Furthermore, Knuth~\mbox{\cite[Ex.~101 in Sec.~7.2.1.6]{MR3444818}} asked for a simple listing of the spanning trees of the complete graph~$K_n$ by edge exchanges.}
The precise definition of `simple' is not specified, but it could mean fast ranking and unranking algorithms, and possibly also a more fine-grained explanation of Cayley's formula~$n^{n-2}$ for the number of all spanning trees.

Li and Liu~\cite{MR2387612} consider the flip graph of circuits of a matroid, whose vertices are the circuits and two of them are connected by an edge if their intersection is nonempty.
They proved that the circuit flip graph of a connected matroid is Hamilton-uniform, provided that there are at least four circuits; see also~\cite{MR3057220}.
In the case of graphical matroids, a Hamilton cycle in the flip graph corresponds to a cyclic listing of the cycles of a graph so that any two consecutive cycles share at least one edge.
Li and Liu~\cite{MR2328946} proved that the circuit flip graph of a connected matroid is also edge-pancyclic; see also~\cite{MR2412260,li_liu_2012}.

\subsection{Perfect matchings}
\label{sec:match}

For a graph~$H$, we consider the flip graph~$G(H)$ that has as vertices all perfect matchings of~$H$, and two of them are connected by an edge in~$G(H)$ if they differ by an alternating cycle.
By the results from~\cite{MR638286} discussed in Section~\ref{sec:01}, $G(H)$ is either Hamilton-connected or a hypercube, which is known to be Hamilton-laceable; see also~\cite{MR762893} and~\cite{MR919493}.
Naddef~\cite{MR762890} showed that the flip graph~$G(H)$ is also edge-pancyclic, provided that the set of edges of~$H$ that are contained in some perfect matching form a connected subgraph of~$H$.

\subsection{Hamilton cycles and Eulerian tours}

Fleischner, Hor\'ak and \v{S}ir\'a\v{n}~\cite{MR1270502} consider the flip graph~$G(H)$ of Hamilton cycles of~$H$, where any two cycles~$C$ and~$C'$ are connected by an edge in~$G(H)$ if the symmetric difference of their edge sets is a 4-cycle in~$H$, and the vertices of this 4-cycle appear consecutively on~$C$ and~$C'$.
They showed that~$G(K_n)$ for $n\geq 4$ is edge-Hamiltonian (see also Sierksma's papers~\cite{MR1218043,MR1492012}).
Note that a Hamilton cycle in~$G(K_n)$ corresponds to a cyclic listing of rosary permutations of~$[n]$ by circularly adjacent transpositions.
\openproblem{53}{It would be very interesting to investigate the Hamiltonicity of~$G(H)$ for other graphs~$H$, for example $H=K_{a,b}$.}

Let~$H$ be a graph that admits a Eulerian tour, i.e., $H$ is connected and each vertex has even degree.
The graph~$H$ may also have multiple edges, but no loops.
Zhang and Guo~\cite{MR830589} consider the flip graph~$G(H)$ that has as vertices all Eulerian tours of~$H$ without considering their orientation, i.e., a forward or backward traversal of the same tour is represented as a single vertex in~$G(H)$.
Moreover, two Eulerian tours are connected by an edge in~$G(H)$ if they differ by reversing the orientation of traversing a closed proper subtour.
This flip operation is known as \emph{Kotzig transformation}.
They proved that if~$H$ has at least three Eulerian tours, then~$G(H)$ is edge-Hamiltonian.
In~\cite{MR887368}, the same authors also consider an analogous question for directed graphs.
Let~$H$ be a digraph that admits a Eulerian tour, i.e., $H$ is connected and for each vertex, the in-degree equals the out-degree.
As before, multiple edges in~$H$ are allowed, but no loops.
The flip graph~$G(H)$ has as vertices all directed Eulerian tours of~$H$, and two of them are connected if they differ by exchanging two directed subtrails between the same pair of vertices.
The authors showed that if~$H$ has at least three directed Eulerian tours, then~$G(H)$ is edge-Hamiltonian.

Guan and Pulleyblank~\cite{MR800996} consider orientations of an undirected graph~$H$ that admit a directed Eulerian tour.
The corresponding flip graph~$G(H)$ has as vertices all Eulerian orientations of~$H$, and two of them are connected in~$G(H)$ if they differ by reversing an oriented (simple) cycle.
They proved that~$G(H)$ is either Hamilton-connected or a hypercube, which is known to be Hamilton-laceable.

\subsection{Elimination orders and trees}
\label{sec:elim}

A \emph{perfect elimination ordering} (PEO) of a graph is an ordering $v_1,\ldots,v_n$ of its vertices such that for every $i=1,\ldots,n$, the neighbors of $v_i$ among $v_1,\ldots,v_{i-1}$ form a clique.
It is well known that a graph~$H$ admits a PEO if and only if $H$ is chordal.
Chandran, Ibarra, Ruskey and Sawada~\cite{MR2022580} present a Gray code for the PEOs of any chordal graph that performs one or a product of two adjacent transpositions in each step.
Their construction is based on the antimatroid Gray code from~\cite{MR1267190} (see also~\cite{MR2159645}).

An \emph{elimination tree} of a connected graph~$H$ is a rooted unordered tree, obtained by removing a vertex~$v$ of~$H$, which becomes the root of the tree, and recursing on the connected components of~$H-v$ that become the children of~$v$ in the tree.
An \emph{elimination forest} of a disconnected graph~$H$ is a collection of elimination trees, one for each connected component of~$H$.
Elimination trees generalize many combinatorial objects, such as binary trees (if $H$ is a path), permutations (if $H$ is a complete graph) and bitstrings (if $H$ is a disjoint union of edges).
Manneville and Pilaud~\cite{MR3383157} showed that if~$H$ has at least two edges, then its elimination trees can be generated cyclically by tree rotations, which naturally generalize classical rotations in binary trees.
The corresponding flip graphs arise as skeleta of \emph{graph associahedra}, which are a class of polytopes introduced by Carr and Devadoss~\cite{MR2239078,MR2479448} and Postnikov~\cite{MR2487491}.
Cardinal, Merino and M\"utze~\cite{MR4415121} present a simple and efficient algorithm to compute a Gray code for the elimination forests of~$H$ by tree rotations for the case when $H$ is a chordal graph, and this Gray code is cyclic if $H$ is also 2-connected.
This Gray code generalizes the Gray code for binary trees due to Lucas, Roelants van Baronaigien and Ruskey~\cite{MR1239499}, the SJT algorithm for permutations, and the BRGC for bitstrings.

\subsection{Acyclic orientations}
\label{sec:acyclic}

Savage, Squire and West~\cite{MR1267311} consider a graph~$H$ and all of its acyclic orientations.
The flip graph~$G(H)$ has these acyclic orientations as vertices, and two of them are adjacent if the differ by the reversal of a single arc.
Observe that for a graph~$H$ with $m$ edges, $G(H)$ is a subgraph of the $m$-cube, in particular, $G(H)$ is bipartite.
Moreover, if $H$ is a forest, then $G(H)$ is the full $m$-cube, i.e., acyclic orientations are in 1-to-1 correspondence to bitstrings of length~$m$, and reversing an arc corresponds to flipping a single bit.
Observe also that if $H=K_n$, then acyclic orientations are in 1-to-1 correspondence to permutations of~$[n]$, and reversing an arc corresponds to an adjacent transposition.
For these two cases we immediately obtain Hamilton cycles in~$G(H)$ by the BRGC and SJT algorithm, respectively.
The aforementioned paper of Savage, Squire and West showed more generally that~$G(H)$ has a Hamilton cycle for any chordal graph~$H$, and the resulting orderings generalize the BRGC and SJT.
To describe the other results, we introduce two more families of graphs:
The \emph{ladder}~$L_n$ is defined as the Cartesian product $L_n:=P_n\times P_2$, and the \emph{wheel}~$W_n$ is obtained from~$C_n$ by adding a new vertex adjacent to all others.
The authors constructed Hamilton cycles in~$G(C_n)$ and~$G(L_n)$ for odd~$n$, and a Hamilton path in~$G(W_n)$ for odd~$n$.
On the negative side, they proved that the following flip graphs have no Hamilton path because of balancedness problems: $G(C_n)$, $G(L_n)$ and $G(W_n)$ for even~$n$, $G(K_{a,b})$ for $a,b\geq 2$ if $a$ or~$b$ is even.
\openproblem{54}{They conjectured that $G(W_n)$ has a Hamilton cycle for odd~$n$ (not just a Hamilton path), and they also raise the problem of finding a Hamilton path or cycle in~$G(K_{a,b})$ when $a$ and~$b$ are both odd.
More generally, they conjectured that $G(H)$ has a Hamilton cycle if and only if its partition classes have the same size.}
This question was reiterated by Savage and Zhang~\cite{MR1609342} and by Edelman (see~\cite{MR1491049}).
Towards this question, Squire~\cite{squire_1994} showed that the square~$G^2(H)$ has a Hamilton cycle for any graph~$H$.
Pruesse and Ruskey~\cite{MR1321068} strengthened this result and proved that the prism $G(H)\times P_2$, i.e., the Cartesian product of~$G(H)$ with a single edge, admits a Hamilton cycle.
As $G(H)$ is bipartite, this easily implies that $G^2(H)$ has a Hamilton cycle.
Squire~\cite{MR1606508} gives an efficient algorithm for computing a Hamilton cycle in the square~$G^2(H)$, i.e., for listing all acyclic orientations of~$H$ by reversing at most two arcs in each step.

The problem of listing the acyclic orientations of~$H$ by arc reversals can be generalized so that it simultaneously generalizes the problem of listing the linear extensions of a poset by adjacent transpositions discussed in Section~\ref{sec:posets}.
Specifically, for a graph~$H$, we consider a subset~$F$ of edges of~$H$ to which we assign a fixed acyclic orientation, and we let $G(H,F)$ be the subgraph of~$G(H)$ induced by acyclic orientations that preserve the orientation of the edges in~$F$.
Clearly, we have $G(H,\emptyset)=G(H)$, so this captures acyclic orientations as before.
On the other hand, if $F$ are the cover edges of a poset~$P$ on~$[n]$, oriented upwards, then the acyclic orientations of~$G(K_n,F)$ (i.e., we take $H$ as the complete graph on~$[n]$) are in 1-to-1 correspondence with the linear extensions of~$P$, and reversing an arc corresponds to an adjacent transposition.
Savage~\cite{MR1491049} asked whether $G(H,F)$ admits a Hamilton cycle whenever its partition classes have the same size, but a counterexample was constructed by Naatz~\cite{MR1814541}.
Pruesse and Ruskey~\cite{MR1321068} also give an example showing that the prism of~$G(H,F)$ need not be Hamiltonian, and Squire~\cite{squire_1994} strengthened this by an example showing that the square $G^2(H,F)$ need not be Hamiltonian.
\openproblem{55}{In general it is open under which conditions $G(H,F)$ or its prism or square admit a Hamilton path or cycle.
As this generalizes Ruskey's conjecture discussed in Section~\ref{sec:posets}, this is arguably a hard problem.}

Pilaud~\cite{MR4490908} equipped~$G(H)$ with a poset structure, and found conditions on~$H$ for when this poset is a lattice with various additional properties.
\openproblem{56}{He also describes the corresponding lattice congruences and realizations of their quotients as polytopes, and he raises the question whether they admit a Hamilton cycle.}
Substantial progress towards this problem was made by Cardinal, Hung, Merino, Mi\v{c}ka and M\"utze~\cite{MR4614413}, who answered a slightly weaker version of this question, namely for Hamilton paths instead of cycles, and those can be computed with a simple greedy algorithm.

\subsection{Vertex colorings}

Choo and MacGillivray~\cite{MR2785821} consider the flip graph~$G(H,k)$ of $k$-colorings of a fixed graph~$H$, where two colorings are adjacent if they differ in the color of a single vertex of~$H$.
They showed that for every~$H$ there is a minimum integer~$k_0(H)$, such that for all $k\geq k_0(H)$, the flip graph~$G(H,k)$ admits a Hamilton cycle.
They also proved that $k_0(K_n)=n+1$ for $n\geq 2$, $k_0(C_n)=4$ for $n\geq 3$, and for any tree~$H$ we have $k_0(H)=3$, unless $H$ is a star with an odd number of vertices, in which case $k_0(H)=4$.
Celaya, Choo, MacGillivray and Seyffarth~\cite{MR3578754} proved that $k_0(K_{a,b})=3$ if $a$ and $b$ are both odd and $k_0(K_{a,b})=4$ otherwise.
Bard~\cite{bard_2014} established the upper bound $k_0(K_{a_1,a_2,\ldots,a_r})\leq 2r$.
He also asked whether $G(H,k)$ having a Hamilton cycle implies that $G(H,k+1)$ also has one.
Stijn Cambie showed that the answer is `No', by taking~$H$ as $K_{5,5}$ minus a perfect matching (cf.~\cite{MR2378926}).
The resulting graph~$G(H,3)$ has a Hamilton cycle, whereas $G(H,5)$ is disconnected and hence has none.
Cavers and Seyffarth~\cite{MR4045409} showed that any 2-tree $H$ satisfies $k_0(H)\in\{4,5\}$, and they characterized precisely which 2-trees attain which of these two possible values.
A \emph{2-tree} is obtained from a triangle by repeatedly adding a new vertex that is joined to the end vertices of an existing edge.
McDonald~\cite{mcdonald_2015} considers a variation that allows recoloring more than one vertex at a time, namely at most~$\ell$ vertices at a time, and he asked for the minimum~$\ell$ for which the resulting flip graph has a Hamilton cycle.

Haas~\cite{MR2875615} introduced another kind of flip graph for colorings.
Specifically, we consider isomorphism classes of $k$-colorings of~$H$ that differ only in permuting the colors.
For a fixed ordering~$\pi$ of the vertices of~$H$, the lexicographically smallest coloring is chosen as the representative of each isomorphism class.
The flip graph~$G(H,\pi,k)$ has those representative colorings as vertices, and an edge between any two that differ in the color of a single vertex of~$H$.
She showed that if $H$ is a tree on at least 4 vertices and $k\geq 3$, then there is an ordering~$\pi$ such that $G(H,\pi,k)$ has a Hamilton cycle.
\openproblem{57}{She also proved that $G(C_n,\pi,k)$ is connected for all $n\geq 3$, $k\geq 4$ and some~$\pi$, and she mentions computer experiments suggesting that this flip graph has a Hamilton cycle for $n>5$.}
Haas and MacGillivray~\cite{MR3798891} showed that $G(K_{a,b},\pi,k)$ has a Hamilton path for all $a,b\geq 2$, $k\geq 3$, and some~$\pi$.
They also report some positive results for graphs~$H$ obtained via unions and joins, and some negative results about $H=K_{a_1,a_2,\ldots,a_r}$ for $r\geq 3$, $a_i\geq 2$, $k\geq r+1$.

Finbow and MacGillivray~\cite{finbow_macgillivray_2014} consider the flip graph~$G(H,k)$ of vertex partitions of~$H$ into at most $k$ independent sets, where two of them are adjacent if they differ by moving a single vertex from one class to another.
Note that if $H$ is an independent set of size~$n$ (=complement of~$K_n$), then $G(H,k)$ is the flip graph of set partitions with at most~$k$ parts under element exchanges.
The authors showed that $G(H,n)$ has a Hamilton cycle for every $n$-vertex graph~$H$ except for $K_n$ and $K_n-e$, and that $G(H,k)$ has a Hamilton cycle for any tree~$H$ with at least 4 vertices and~$k\geq 3$.

\subsection{Vertex labelings}

Knor~\cite{MR1301949} considers all non-isomorphic vertex labelings of an $n$-vertex graph~$H$ with integers $1,\ldots,n$, and two of these vertex-labeled graphs are adjacent if one can remove and reinsert $k$ edges to transform one graph into the other.
For example, if $H=P_n$, then there are $n!/2$ non-isomorphic labelings of~$H$, and if $k=1$ then the flip operation is the well-known P\'osa rotation.
In general, we write $G(H,k)$ for the corresponding flip graph.
It is easy to see that $G(H,k)$ is vertex-transitive, as the neighbors of some vertex labeling of~$H$ only depend on~$H$ and not on the labeling, which makes this problem an interesting test case for Lov\'asz' conjecture (provided that $G(H,k)$ is connected, which is not always the case).
Knor proved that $G(P_n,1)$ has Hamilton cycle for $n\geq 3$, $G(C_n,2)$ has a Hamilton cycle for $n\geq 4$, and $G(F_n,1)$ is Hamilton-connected for $n\geq 7$, where $F_n$ is the \emph{fork}, i.e., the $n$-vertex graph obtained by appending two vertices to the same end vertex of the path~$P_{n-2}$.
He further showed that $G(K_{a,b},a+b-2)$ for $a\neq b$ has the Johnson graph~$J(a+b,a)$ as a spanning subgraph, so it is Hamilton-connected by the results discussed in Section~\ref{sec:comb-trans} (for $a=b$ an analogous statement holds for $J(2a-1,a)$).
Another interesting special case is $G(K_{a+1,a},a)$, which is isomorphic to the odd graph~$K(2a+1,a)=O_a$.

\subsection{Dominating sets}

Adaricheva et al.~\cite{adaricheva_et_al_2021} consider the flip graph~$G(H)$ of dominating sets of~$H$, where two dominating sets are connected by an edge if they differ by adding or removing a single vertex of~$H$ from the set.
As the number of dominating sets of any graph~$H$ is odd~\cite{brouwer_csorba_schrijver_2009}, the bipartite graph~$G(H)$ has two partition classes of different sizes and therefore no Hamilton cycle, so we only need to to focus on finding Hamilton paths.
Note that every nonempty subset of vertices of~$K_n$ is a dominating set, so $G(K_n)$ is isomorphic to the $n$-cube $Q_n$ minus a single vertex.
Similarly, $G(K_{1,n})$ is isomorphic to~$Q_n$ plus a pending edge attached to one vertex.
As $Q_n$ has a Hamilton cycle, it follows immediately that $G(K_n)$ and $G(K_{1,n})$ have a Hamilton path.
The aforementioned paper shows moreover that $G(P_n)$ has a Hamilton path, and $G(K_{a,b})$ has a Hamilton path if and only if $a$ or $b$ is odd (if both are even there are balancedness problems).

Adaricheva et al.~\cite{MR4496479} proved that $G(H)$ has a Hamilton path for any tree~$H$.
They also consider the case of cycles $H=C_n$.
Note that characteristic vectors of dominating sets of~$C_n$ are precisely bitstrings of length~$n$ that cyclically avoid the factor~$000$ (three consecutive vertices on~$C_n$ that are not in the dominating set leave the middle vertex undominated), i.e., they are complements of Lucas words.
Therefore, by the results of Baril and Vajnovszki~\cite{MR2187406} discussed in~Section~\ref{sec:factor-avoid}, there is a Hamilton path in~$G(C_n)$ if and only if $n$ is not a multiple of~4.
Similarly, characteristic vectors of dominating sets of~$P_n$ are precisely bitstrings of length~$n$ that avoid the factor~$000$ and also~$00$ at the first or last two positions, which are a proper subset of complemented Fibonacci words.

\openproblem{58}{Haas and Seyffarth~{\cite{MR3195801}} asked whether the flip graph~$G(H,k)$ of dominating sets of size at most~$k$ in~$H$ has a Hamilton path or cycle.}
More results and open problems about flip graphs on colorings and dominating sets of graphs are surveyed by Mynhardt and Nasserasr~\cite{mynhardt_nasserasr_2019}.

\subsection{Degree sequences and 0/1-matrices}

Arikati and Peled~\cite{MR1672981} consider degree sequences of graphs, which are integer sequences $d=(d_1,\ldots,d_n)$ for which there exists a simple graph with vertex set~$[n]$ such that the degree of the vertex~$i$ is~$d_i$ for all $i=1,\ldots,n$.
The flip graph~$G(d)$ has as vertices all graphs that realize the degree sequence~$d$, and two of them are connected by an edge in~$G(d)$ if they differ by removing two edges and adding two edges that together form a 4-cycle on~$[n]$, an operation that clearly preserves all degrees.
The authors proved that $G(d)$ has a Hamilton cycle provided that $d$ has majorization gap~1 (see the paper for definitions).
\openproblem{59}{In general, the question whether $G(d)$ admits a Hamilton path or cycle is open.}

\setlength{\tabcolsep}{0.3pt}
\renewcommand{\arraystretch}{0.5}

Brualdi~\cite{MR585770} considers $m\times n$ matrices of 0s and 1s, with fixed row sum vector~$r=(r_1,\ldots,r_m)$ and column sum vector~$s=(s_1,\ldots,s_n)$, where $0<r_i<n$ for $i=1,\ldots,m$ and $0<s_j<m$ for $j=1,\ldots,n$.
The flip graph $G(r,s)$ has as vertices all 0/1-matrices whose row and column sums are~$r$ and~$s$, respectively, and two of them are connected by an edge in~$G(r,s)$ if the matrices agree in all but four entries that form submatrices $\begin{pmatrix} 1 & \!\! 0 \\ 0 & \!\! 1 \end{pmatrix}$ and $\begin{pmatrix} 0 & \!\! 1 \\ 1 & \!\! 0 \end{pmatrix}$, respectively.
\openproblem{60}{Brualdi raised the question whether the flip graph~$G(r,s)$ admits a Hamilton path or cycle.}
Li and Zhang~\cite{MR1279625} showed that the answer is `yes' if $r=(1,\ldots,1)$ or $s=(1,\ldots,1)$, proving more strongly that $G(r,s)$ is edge-Hamiltonian.
By taking complements, the same result follows if $r=(n-1,\ldots,n-1)$ or $s=(m-1,\ldots,m-1)$.

In the aforementioned paper~\cite{MR1672981}, Arkati and Peled point out that Brualdi's problem can be rephrased equivalently using degree sequences as follows:
We can interpret an $m\times n$ matrix with row and column sums~$r$ and~$s$, respectively, as the biadjacency matrix of a bipartite graph with partition classes of sizes~$m$ and~$n$.
We then turn the bipartite graph into a split graph by adding all edges in one of the partition classes, making this partition class a clique.
Note that the degree sequence~$d$ of the resulting split graph yields a flip graph~$G(d)$ which is isomorphic to~$G(r,s)$.
This is because in a split graph, every flip operation removes and adds two edges that go between the independent set and the clique.
Therefore, Brualdi's question is equivalent to asking whether the flip graph~$G(d)$ for the degree sequence~$d$ of any split graph admits a Hamilton path or cycle.

\section{Greedy Gray codes}
\label{sec:greedy}

Throughout this survey, we encountered a number of versatile methods for constructing Gray code listings for several different combinatorial objects, by a suitable general encoding of the objects.
These include flip-swap languages (Section~\ref{sec:flip-swap}), bubble languages (Section~\ref{sec:bubble}), pattern-avoiding permutations (Section~\ref{sec:avoid}), linear extensions of posets (Section~\ref{sec:posets}), elimination trees (Section~\ref{sec:elim}) and acyclic orientations (Section~\ref{sec:acyclic}).

Williams~\cite{MR3126386} discovered another powerful general paradigm for constructing Gray codes, which he dubbed the \emph{greedy approach}.
The idea is to describe the listing implicitly by an algorithm, in which the next object in the list is defined by applying a flip operation to the preceding object, and this operation is selected greedily from a set of applicable flip operations so as to generate an object which has not been encountered in the list before.
Williams showed in his paper that many classical Gray codes can be described by such a greedy rule:
\begin{itemize}[leftmargin=5mm, noitemsep, topsep=3pt plus 3pt]
\item the BRGC for bitstrings (Figure~\ref{fig:bits1}~(a)): flip the rightmost possible bit;
\item the homogeneous transposition Gray code for combinations due to Eades and McKay (Figure~\ref{fig:comb1}~(b)): transpose the rightmost possible~1 homogeneously with the rightmost possible~0;
\item the SJT algorithm for permutations (Figure~\ref{fig:perm1}~(d)): transpose the largest possible entry with its neighbor;
\item the Ord-Smith/Zaks ordering of permutations (Figure~\ref{fig:perm1}~(k)): reverse the shortest possible prefix;
\item the Lucas, Roelants van Baronaigien and Ruskey Gray code for binary trees (Figure~\ref{fig:231cat}): rotate the edge with the largest possible inorder label.
\end{itemize}
Besides describing known Gray codes in a unified way, the greedy approach has proven to be very useful for discovering new Gray codes.
One major application of the greedy method are the Gray codes for zigzag languages of permutations introduced by Hartung, Hoang, M\"utze and Williams~\cite{MR4391718}, which include many classes of pattern-avoiding permutations (recall Section~\ref{sec:avoid}), and which encode several classes of rectangulations~\cite{MR4598046} (recall Section~\ref{sec:rect}), elimination trees of chordal graphs~\cite{MR4415121} (recall Section~\ref{sec:elim}), and acyclic orientations of chordal graphs (recall Section~\ref{sec:acyclic}).

Another major application area where the greedy method thrives are 0/1-polytopes.
As was shown by Merino and M\"utze~\cite{MR4720318}, the greedy approach can be used to compute a Hamilton path on the skeleton of any 0/1-polytope (recall Section~\ref{sec:01}).
This vastly extends the method to all families of combinatorial objects that can be naturally encoded by 0/1-strings, including spanning trees of graphs (recall Section~\ref{sec:span} and~\cite{DBLP:conf/fun/MerinoMW22}), matchings of graphs (recall Section~\ref{sec:match}), etc. 

The Gray codes derived from the greedy approach all have an underlying hierarchical generation tree, instances of which can be found in several predecessor papers, including Hurtado and Noy's tree of triangulations~\cite{MR1723053}, Takaoka's twisted lexico tree~\cite{DBLP:journals/cj/Takaoka99}, the production trees of the ECO method~\cite{MR2074330}, and in Li and Sawada's reflectable languages~\cite{MR2483809} (see also~\cite{DBLP:conf/iccsa/XiangCU10}).

The key to understanding the success of the greedy method was recently revealed in~\cite{MR4720318}, by exhibiting an underlying connection to Gray code listings that are genlex.
Specifically, if there is a genlex listing of combinatorial objects using a particular flip operation, where objects with the same prefix (or suffix) appear consecutively, then the greedy rule to perform a flip that modifies the shortest possible suffix (or prefix, respectively) can generate this listing, when using a suitable tie-breaking rule.
This explains why the greedy approach works in all the aforementioned cases.
More importantly, it extends the applicability of the greedy method to any of the other genlex Gray codes mentioned in this survey.

\section*{Acknowledgements}

The author thanks Jean-Luc Baril, Clemens Huemer, Dragan Maru\v{s}i\v{c}, Arturo Merino, Ond\v{r}ej Mi\v{c}ka, Namrata, Carla Savage, Joe Sawada and Terry Wang for numerous corrections and remarks that helped improving this survey.

\bibliographystyle{alpha}
\bibliography{refs}

\end{document}